\documentclass{amsart}
\usepackage[all]{xy}
\usepackage{amssymb}
\usepackage{amsmath}

\swapnumbers
\newtheorem{thm}[subsection]{Theorem}
\newtheorem{prop}[subsection]{Proposition}
\newtheorem{lem}[subsection]{Lemma}
\newtheorem{cor}[subsection]{Corollary}
\newtheorem{defn}[subsection]{Definition}
\newtheorem{notn}[subsection]{Notation}
\newtheorem{rem}[subsection]{Remark}
\newtheorem{exmpls}[subsection]{Examples}

\newcommand{\Irr}{{\mbox{\it Irr}_{S_n}}}
\newcommand{\Sh}{{\mbox{\it Sh}(n,m)}}
\newcommand{\Shm}{{\mbox{\it Sh}(n_{1},\dots ,n_{r})}}

\def\noi{\noindent}
\def\ms{\medskip}
\def\bs{\bigskip}

\def\m{\, \bot\,  }

\def\P{\noindent \emph{Proof.} }
\def\n{\underline {\bold n}}
\def\m{\underline {\bold m}}

\def\cc{\mathfrak{c}}

\def\arboluno{\smallmatrix \searrow &&\swarrow \\ &\bullet &\\&\downarrow &\endsmallmatrix }

\begin{document}

\author[M. Ronco]{Mar\'\i a Ronco}

\address{Depto. de Matem\'atica, Fac. de Ciencias\\
Universidad de Valpara\'\i so\\
Avda. Gran Breta\~na 1091\\
Valpara\'\i so. Chile}
\email{maria.ronco@uv.cl}
\thanks{Research supported by FONDECYT 1060224}
\title{Shuffle bialgebras}


\begin{abstract} The goal of our work is to study the spaces of primitive elements of the Hopf algebras associated to the permutahedra and the associahedra. We introduce 
the notion of shuffle bialgebras, and compute the subpaces of primitive elements associated to these algebras. These spaces of primitive elements have natural structure of 
some type of algebras which we describe in terms of generators and relations. Applying these results we are able to compute primitive elements of other combinatorial 
Hopf algebras, and describe the algebraic theories associated to them. 
\end{abstract}
\maketitle

\section*{Introduction} \label{S:int} 

The aim of this paper is to compute the subspaces of primitive elements of some combinatorial Hopf algebras and to describe them as free objects for some new types 
of algebras.

The main examples of combinatorial Hopf algebras studied are the Malvenuto-Reutenauer Hopf algebra (see \cite{MR}), the algebra spanned by the faces of the 
permutahedra (see \cite{Ch},\cite{Be}), the algebra of functions between finite sets, and the algebra of planar rooted trees (see \cite{LR1}). In all cases we describe
 them as free objects for some algebraic theories, and compute the subspaces of their primitive elements. Although the primitive elements of the Malvenuto-Reutenauer 
algebra and of the algebra of planar binary rooted trees have been previously computed in \cite{AS1}, \cite{AS2} and \cite{HNT}, our description has the advantage of 
showing them as free objects for some algebraic theories described in terms of generators and relations. These results give in each case a Cartier-Milnor-Moore type theorem.

In the general case, there does not exist a standard method to compute the space of primitive elements of a non-cocommutative coalgebra. The examples studied in this paper 
have one point in common: they are equipped with an associative product $\times $, called the concatenation product, which verifies a nonunital infinitesimal relation with 
the coproduct. Nonunital infinitesimal bialgebras were introduced in \cite{LR3}, where we proved that any connected nonunital infinitesimal bialgebra is isomorphic to the cofree 
coalgebra spanned by the space of its primitive elements. This result is the main tool used in the present work to compute the primitive elements.

We deal first with shuffle algebras, whose free objects are closely related to the Malvenuto-Reutenauer and Solomon-Tits Hopf algebras. 
Shuffle algebras are a particular case of monoids in the category of ${\mathcal S}$-modules, as described in \cite{PS} and \cite{Liv}, where the operations do not preserve
 the action of the symmetric group. Afterwards, we extend our results to other algebraic structures: the preshuffle algebras and the grafting algebras, the last ones 
 are given by the underlying spaces of non-symmetric operads. The way we study them is 
largely inspired by the treatment given by J.-L. Loday to the so-called triples of operads, see \cite{Lo1}. Let us describe briefly the method employed: 
\begin{enumerate}\item Given a linear algebraic theory ${\mathcal T}$, we introduce the notion of 
${\mathcal T}$ bialgebra, in such a way that any free ${\mathcal T}$ algebra has a natural structure of ${\mathcal T}$ bialgebra.
\item In a second step, we identify the theory ${\mathcal T}$ with the space of all the operations of the theory, and compute a basis for a subspace ${\mathcal Prim}_{\mathcal T}$ 
of ${\mathcal T}$, such that the space of primitive elements of any ${\mathcal T}$ bialgebra is closed under the action of the elements of ${\mathcal Prim}_{\mathcal T}$.
\item Afterwards, we prove that any free ${\mathcal T}$ algebra ${\mathcal T}(X)$ is isomorphic, as a coalgebra, to the cofree coalgebra spanned by the space 
${\mathcal Prim}_{\mathcal T}(X)$, generated by the operations of ${\mathcal Prim}_{\mathcal T}$ on the elements of $X$. Applying results proved in \cite{LR3}, we get that 
${\mathcal Prim}_{\mathcal T}(X)$ is the space of primitive elements of ${\mathcal T}(X)$.
\item Finally, we describe the algebraic theory associated to ${\mathcal Prim}_{\mathcal T}$ in terms of generators and relations; and prove that the category of connected 
${\mathcal T}$ bialgebras is equivalent to the category of ${\mathcal Prim}_{\mathcal T}$ algebras.
\end{enumerate}

It is quite easy to compute the theory ${\mathcal Prim}_{sh}$ when ${\mathcal T}$ is the theory of shuffle algebras. The other examples follow 
easily from this case.
\ms

The paper is organised as follows:

The first section of the paper recalls some contructions on planar rooted trees and on permutations, needed in the following sections.

In Section 2 we give the definition of a shuffle algebra, describe the free objects for this theory in terms of permutations, and give the main examples. 
Shuffle bialgebras are introduced in Section 3. We also show that there exist natural functors between the categories of graded infinitesimal bialgebras and the category of shuffle
bialgebras.

In Section 4 we compute the primitive elements of a shuffle  bialgebra, and prove that they are given by some new algebraic objects called 
${\mathcal Prim}_{sh}$ algebras. We prove a Cartier-Milnor-Moore Theorem in this context, showing that the category of connected shuffle
 bialgebras is equivalent to the category of ${\mathcal Prim}_{sh}$ algebras.

The next Section is devoted to introduce the relationship between the algebraic theories introduced in the paper and the operads of dendriforn, infinitesimal and $2$--
associative algebras. In Section 6 we show that any coproduct on a free shuffle algebra gives rise to a boundary map, and show that the classical boundary map of the 
permutohedra  is an example of this construction. 

In Section 7 the notion of preshuffle bialgebras is introduced. As a particular case of preshuffle algebras, 
we introduce grafting algebras, and prove that the free objects for this theory are given by the spaces spanned 
by planar coloured trees.  In fact, the notion of grafting algebra coincides with non symmetric algebraic operad, however since we study them as algebras we keep the name 
of grafting algebra to designate them. From? the definition of ${\mathcal Prim}_{sh}$ algebras, we compute the suspaces of primitive elements of preshuffle bialgebras and 
grafting bialgebras, and describe any connected preshuffle (respectively, grafting) bialgebra as an enveloping algebra over its primitive part. 
A variation of Cartier-Milnor-Moore Theorem is proved in this context, showing that the category of connected preshuffle 
(respectively, grafting) bialgebras is equivalent to the category of ${\mathcal Prim}_{psh}$ (respectively, ${\mathcal Prim}_{gr}$) algebras.

 The last section of the paper contains some applications of the computations of primitive elements made for preshuffle, shuffle and grafting bialgebras to some good triples of 
 operads (see \cite{Lo1}).
 
  \noindent {\bf Acknowledgement.}
The author wants to thank the support of Ecos-Conicyt Project E06C01 during the correction of this paper, as well as J.-L. Loday and T. Schedler for helpful comments.

\section{Preliminaries}
\ms

We introduce here some definitions and notations that are used in the paper.

Let $K$ be a field, $\otimes $ denotes the tensor product of vector spaces over $K$. Given a graded $K$-vector space $A$, $A_{+}$ is the space 
$A\oplus K$ equipped with the canonical maps $K \hookrightarrow A_{+}\longrightarrow K$. 

\noi Given a graded vector space $A=\displaystyle {\bigoplus _{n\geq 0}A_{n}}$, we denote the degree of an homogeneous element $x\in A_{n}$ by
$\vert x\vert =n$.

\noi For any set $X$, we denote by $K[X]$ the vector space spanned by $X$.

\noi Given a $K$-vector space $V$, the tensor space over $V$ is the graded vector space $T(V):=\displaystyle{\bigoplus _{n\geq 0} V^{\otimes n}}$. The reduced 
tensor space ${\overline T}(V)$ over $V$ is the subspace $\displaystyle {\bigoplus _{n\geq 1}V^{\otimes n}}$.

\noi The space ${\overline T}(V)$ with the concatenation product, given by:
$$(v_{1}\otimes \dots \otimes v_{n})\cdot (w_{1}\otimes \dots \otimes w_{m}):=v_{1}\otimes \dots \otimes v_{n}\otimes w_{1}\otimes 
\dots \otimes w_{m},$$
for $v_{1},\dots ,v_{n},w_{1},\dots ,w_{m}\in V$, is the free associative algebra spanned by $V$.

This product is extended to $T(V)$ in the unique way such that the unit $1_{K}$ of the field $K$ becomes the unit for the concatenation product.
\bs

\noi{\bf Coalgebras.} A coalgebra $C$ over $K$ is a vector space, equipped with a coproduct $\Delta :C\longrightarrow C\otimes C$, which is 
coassociative.
\ms

\noi We use Sweedler\rq s notation, and denotes  $\Delta (x)=\sum x_{(1)}\otimes x_{(2)}$, for $x\in C$.
\ms

\noi A coalgebra $C$ is counital if there exists a linear map $\epsilon: C\longrightarrow K$ such that 

\noi $(\epsilon \otimes Id_{C})\circ \Delta =Id _{C}=(Id_{C}\otimes \epsilon )\circ \Delta $, where we identify 
$K\otimes C$ and $C\otimes K$ with $C$, via the canonical isomorphism. 
\ms

\noi For a counital coalgebra $(C, \Delta ,\epsilon )$, we define the reduced coproduct 

\noi $\overline {\Delta }:=\Delta - Id_{C}\otimes \epsilon -\epsilon\otimes 
Id_{C}$ on $ Ker(\epsilon)$. Note that $\overline {\Delta }: Ker(\epsilon)\longrightarrow  Ker(\epsilon)\otimes  Ker(\epsilon)$ is coassociative too.

\noi Let $C=\displaystyle{\bigoplus _{n\geq 0}C_{n}}$ be a graded $K$-vector space. A graded coassociative coproduct on $C$ is a coassociative 
coproduct $\Delta$ such that $\Delta (A_{n})\subseteq \displaystyle{\bigoplus _{i=0}^{n}A_{i}\otimes A_{n-i}}$.
\ms

\noi Given a coassociative coproduct $\Delta $ on $C$ and $r\geq 1$, we denote by $\Delta ^{r}:C\longrightarrow C^{\otimes r+1}$ the homomorphism 
defined recursively as $\Delta ^{1}:=\Delta$ and $\Delta ^{r+1}:=(\Delta ^{r}\otimes Id_{C})\circ \Delta $, for $r\geq 1$.
\ms

Most of the coalgebras we deal with in the paper are not counital. Given such a coalgebra $(C,\Delta )$
we define an associate counital coalgebra $(C_{+},\Delta _{+})$, where $C_{+}:=K\oplus C$ and 
$$\Delta _{+}(x)=\begin{cases} 1\otimes x +x\otimes 1+\Delta (x),&{\rm for}\ x\in C\\
x 1_K\otimes 1_K,&{\rm for}\ x\in K.\end{cases}$$

\noi The space $C_{+}$ with the coproduct $\Delta_{+}$ and the counit $\epsilon$ given by:
$$\epsilon (x):=\begin{cases}0,&{\rm for}\ x\in C,\\
x&{\rm for}\  x\in K,\end{cases}$$
is a counital coassociative coalgebra.

\noi Let $V$ be a vector space, the deconcatenation coproduct on $\overline {T}(V)$ is given by:
$$\overline{\Delta }^{c}(v_{1}\otimes \dots \otimes v_{n}):=\sum _{i=1}^{n-1}(v_{1}\otimes \dots \otimes v_{i})\otimes 
(v_{i+1}\otimes \dots \otimes v_{n}).$$
The graded space $\overline {T}(V)=\overline {T}(V)_+$ equipped with the coproduct $\overline{\Delta }^{c}$ is a coassociative coalgebra.
\ms

\begin{defn}{\rm Let $C={\displaystyle \bigoplus _{n\geq 1}}C_{n}$ be a positively graded $K$-vector space, equipped with a coassociative coproduct $\Delta $. 
An element $x\in C$ is called {\it primitive} if $\Delta (x)=0$. The subspace of primitive elements of $C$ is denoted by Prim$(C)$.}
\end{defn}
\ms

\begin{defn}{\rm Let $(C, \Delta )$ be a coassociative coalgebra, and let $F_{p}C$ be the following filtration on $C$:
$$\displaylines { F_{1}C:= {\rm Prim}(C)$$\hfill \cr
F_{p}C:=\{ x\in C \ \vert\  \Delta (x)\in F_{p-1}C\otimes F_{p-1}C\} .\hfill\cr }$$
Following the definition of D. Quillen (see \cite{Qu}), we say that $C$ is connected if 

\noindent $C={\displaystyle \bigcup _{p\geq 1}}F_{p}C.$ }
\end{defn}
\ms

\noindent The definition of primitive element for the counital coalgebra $C_+$ becomes 

\noindent $x\in {\rm Prim}(C_{+})$ if $\Delta _{+}(x)=x\otimes 1+1\otimes x$. In this case, ${\rm Prim}(C_{+})={\rm Prim}(C)$.
\bs

\bs

\bs

The main purpose of this work is to study bialgebra structures on spaces spanned by (coloured) functions between finite sets, permutations and trees. 
The rest of this section is devoted to introduce definitions and elementary results on these objects.
\ms

\noi{\bf Permutations and shuffles.} Let $S_n$ be the group of permutations on $n$ elements. A permutation $\sigma $ is denoted by its image 
$(\sigma (1),\dots ,\sigma (n))$. 

\noi The element $1_{n}=(1,2,\dots ,n)$ denotes the identity of $S_{n}$. The set 
$S_{\infty }:=\displaystyle{\bigcup _{n\geq 1}S_{n}}$ is the graded set of all permutations.

\begin{defn} {\rm Given $1\leq r\leq n$, a {\it composition} ${\underline {\bold n}}$ of $n$ of length $r$ is a family of positive integers $(n_1,\dots ,n_r)$ such that $\displaystyle
{\sum _{i=1}^rn_i}=n$. The number $r$ is called the length of the composition $\n$.}
\end{defn}

For any  composition $\n= (n_1,\dots ,n_r)$ of $n$, there exists a homomorphism 

\noi ${\mbox {\it S}_{n_1}\times \dots \times S_{n_r}\hookrightarrow S_n}$ given by
$(\sigma _1,\dots ,\sigma _r)\mapsto \sigma _1\times \dots \times \sigma _r$, where
$$(\sigma _1\times \dots \times \sigma _r)(i):=\sigma _k(i-n_1-\dots -n_{k-1})+n_1+\dots +n_{k-1},$$
for $n_1+\dots +n_{k-1}<i\leq n_1+\dots +n_k.$
For any composition $\n$, $S_{\n}$ denotes the subgroup of $S_n$ which is the image of ${\mbox {\it S}_{n_1}\times \dots \times S_{n_r}}$ under
 this embedding.
 
The operation $\times :{\mbox {\it S}_{n}\times S_{m}\longrightarrow S_{n+m}}$ defined previously is an associative product on $S_{\infty }$, called the {\it concatenation}.
\ms

\noi In \cite{MR}, C. Malvenuto and C. Reutenauer introduce the following definition.

\begin{defn}  {\rm A permutation $\sigma \in S_{n}$ is {\it irreducible}  if $\sigma \notin \displaystyle{\bigcup _{i=1}^{n-1}{\mbox {\it S}_i\times S_{n-i}}}$.}
{\rm We denote by $\Irr$ the set of irreducible permutations of $S_{n}$.}
\end{defn}

\noi Note that the graded vector space $K[S_{\infty}]:={\displaystyle \bigoplus _{n\geq 1}}K[S_{n}]$, where $K[S_{n}]$ denotes the vector space spanned by the set of permutations, 
equipped with the concatenation product, is the free associative algebra generated by $\displaystyle{\bigcup _{n\geq 1}\Irr}$.
\ms

\begin{defn} {\rm \begin{enumerate}\item A subgroup $W$ of $S_{n}$ is called a {\it parabolic standard subgroup} if $W=S_{n_1\times \dots \times 
n_r}$, for some composition $(n_{1},\dots ,n_{r})$ of $n$.
\item Given a composition $\n=(n_{1},\dots ,n_{r})$ of $n$, a $(n_{1},\dots ,n_{r})${\it -shuffle}, or $\n${\it -shuffle}, is an element $\sigma $ of $S_{n}$ such that:
$$\sigma^{-1} (n_{1}+\dots +n_{k-1}+1)<\dots<\sigma^{-1} (n_{1}+\dots +n_{k}),\ {\rm for}\ 1\leq k\leq r-1.$$
\end{enumerate}
The set of all $(n_{1},\dots ,n_{r})$-shuffles is denoted indistinctly $\Shm$ or ${\mbox {\it Sh}(\n)}$.}
\end{defn}

\noi Note that\begin{enumerate}
\item [\upshape (i)] ${\mbox{\it Sh}(n_{1},\dots ,n_{i},0,n_{i+1},\dots ,n_{r}})={\mbox {\it Sh}(\n)}$ , for any composition $\n=(n_{1},\dots ,n_{r})$ and any
 $0\leq i\leq r$,  
\item [\upshape (ii)] ${\mbox{\it Sh}(n)}=\{ 1_{n}\}$, 
\item [\upshape (iii)]${\mbox{\it Sh}({\underline {\bold 1}})}={\mbox{\it Sh}(1,\dots ,1)}=S_{n}$, for $n\geq 1$.
\end{enumerate}
\ms

\noi Given positive integers $n,m$, the permutation $\epsilon _{n,m}:=(n+1,\dots ,n+m,1,\dots ,n)$ belongs  to $\Sh$.
\ms

The following results about Coxeter groups are well-known. For the first assertion  see for instance \cite{Sol}, the second one is proved, in a more general 
context, in \cite{BBHT}. 
\ms

\begin{prop}\label{shuff} {\rm (1)} Given a  permutation $\sigma \in S_{n}$ and an integer $0\leq i\leq n$ there exists 
unique elements $\sigma _{(1)}^{i}\in S_{i}$, $\sigma _{(2)}^{n-i}\in S_{n-i}$ and $\gamma \in {\mbox{\it Sh}(i,n-i)}$ such that 
$\sigma = (\sigma _{(1)}^{i}\times \sigma _{(2)}^{n-i})\cdot \gamma $.

\noi {\rm (2)} Given compositions $\n$ of $n$ and ${\mbox{\it Sh}({\underline {\bold m}})}$ of $m$, 
we have that:
$$({\mbox{\it Sh}({\underline {\bold n}})} \times {\mbox{\it Sh}({\underline {\bold m}})})\cdot  \Sh = {\mbox{\it Sh}(\n\cup\m)},$$
where $\n\cup\m:=(n_1,\dots ,n_r,m_1,\dots ,m_p)$.
\end{prop}
\bs

\noi{\bf Coxeter poset of the symmetric group.} The results stated in the previous section imply that, for any composition $\n$ of $n$, and any right coset
$S_{n_{1}\times \dots\times n_{r}}\cdot \sigma $, there exists  a unique element $\delta \in {\mbox{\it Sh}(\n)}$ such that 

\noi $S_{n_{1}\times \dots\times n_{r}}\cdot \sigma = S_{n_{1}\times \dots\times n_{r}}\cdot \delta $.
\ms

\begin{defn} {\rm For $n\geq 1$, the Coxeter poset ${\mathcal P}_{n}$ of $S_{n}$ is the set of cosets:
$${\mathcal P}_{n}:=\{ S_{n_{1}\times \dots \times n_{r}}\cdot \delta\ \mid \ {\rm with}\ \n=(n_{1},\dots ,n_{r})\ {\rm a\ composition\ of}\ n\ 
{\rm and}\ 
\delta \in {\mbox{\it Sh}({\underline {\bold n}})}\},$$ ordered by the inclusion relation.}
\end{defn}

\noi The maximal element of ${\mathcal P}_{n}$ is the unique coset modulo $S_{n}$, that is ${\mbox {\it S}_{n}\cdot 1_n}$ which will be denoted by $\xi _{n}$, and the 
minimal elements are the cosets $S_{1\times \dots \times 1}\cdot \sigma $ which are denoted simply by $\sigma $, for $\sigma \in S_{n}$.
 
 \noi Note that ${\mathcal P}_{n}$ is the disjoint union of the subsets 
 $${\mathcal P}_{n}^{r}:=\{ S_{n_{1}\times \dots \times n_{r}}\cdot \delta \ \mid  \delta \in 
{\mbox{\it Sh}({\underline {\bold n}})}\},\ {\rm for}\ 1\leq r\leq n.$$
 \bs
 
 \noi{\bf Functions on finite sets.} Given positive integers $n$ and $r$, let ${\mathcal F}_{n}^{r}$ be the set of all maps $f:\{ 1,\dots ,n\} \longrightarrow 
 \{ 1,\dots ,r\}$. For $f\in {\mathcal F}_{n}^{r}$, we denote it by its image $(f(1),\dots ,f(n))$. The constant function $(1,\dots ,1)\in {\mathcal F}_{n}^1$ is 
 denoted by $\xi_{n}$. For $n\geq 1$, let ${\mathcal F}_{n}:=\displaystyle{\bigcup _{r=1}^{n}{\mathcal F}_{n}^{r}}$.
 \ms
 
\noi For any $n,m,r$ and $k$, there exists an embedding ${\mathcal F}_{n}^{r}\times {\mathcal F}_{m}^{k}\longrightarrow {\mathcal F}_{n+m}^{r+k}$,
  given by $$f\times g:=(f(1),\dots ,f(n),g(1)+r,\dots ,g(m)+r),\ {\rm for}\ f\in {\mathcal F}_{n}^{r}\ {\rm and}\ g\in{\mathcal F}_{m}^{k}.$$
 \ms

 \begin{rem}\label{funfin}{\rm \begin{enumerate}\item The set of permutations $S_{n}$ is a subset of ${\mathcal F}_{n}^{n}$ and the embedding 
 ${\mathcal F}_{n}^{r}\times {\mathcal F}_{m}^{k}\hookrightarrow {\mathcal F}_{n+m}^{r+k}$, restricted to $S_{n}\times S_{m}$, coincides with the 
 concatenation $\times$, previously defined.
 \item For any element $f\in {\mathcal F}_{n}^{r}$ there exists a unique non-decreasing function $f^{\uparrow}\in {\mathcal F}_{n}^{r}$ and 
 a unique permutation $\sigma _{f}\in \Shm$ such that 
 $$f=f^{\uparrow}\cdot \sigma _{f},$$
 where $n_i=\vert f^{-1}(i)\vert$ for $1\leq i\leq r$, and $\cdot$ denotes the composition of functions.
 \item For $n\geq 1$ and $1\leq r\leq n$, there exists a natural bijection between ${\mathcal P}_{n}^{r}$ and the set of all surjective maps from 
 $\{ 1,\dots ,n\} $ to $\{1,\dots ,r\}$. Given a composition $\n=(n_{1},\dots ,n_{r})$ of $n$ and $\delta \in {\mbox{\it Sh}({\underline {\bold n}})}$, the element 
 $S_{n_{1}\times \dots \times n_{r}}\cdot \delta $ maps to the function $\xi _{\n}\cdot \delta $, where:
 $$\xi _{\n}(j):=(\xi_{n_{1}}\times \dots \times \xi_{n_{r}})(j)= k,$$
  for $n_{1}+\dots +n_{k-1}<j\leq n_{1}+\dots +n_{k}$.
 
 \noindent For example, the element $S_{2,3,1}\cdot (3,1,2,4,6,5)$ maps to the function 
 
 \noi $(2,1,1,2,3,2)\in {\mathcal F}_{6}^{3}$.
 
\noi For $n,r\geq 1$, we may identify ${\mathcal P}_{n}^{r}$ with the subset $\{ f\in {\mathcal F}_{n}\ \mid\ f\ {\rm is\ 
 surjective}\}$ of ${\mathcal F}_{n}$. The map $\times: {\mathcal F}_{n}\times {\mathcal F}_{m}\longrightarrow {\mathcal F}_{n+m}$ restricts to
 $\times :{\mathcal P}_{n}\times {\mathcal P}_{m}\longrightarrow {\mathcal P}_{n+m}$.
 \end{enumerate}}
 \end{rem}
 \ms
 
\noi From now on, we shall use Remark \ref{funfin}, and denote the elements of ${\mathcal P}_{n}$ as surjective maps by their image $f=(f(1),\dots ,f(n))$.
 \ms
 
 Let ${\mathcal K}_{n}$ be the set of all maps $f\in {\mathcal P}_{n}$ verifying the following condition:
 $${\rm if}\ f(i)=f(j),\ {\rm for\ some}\ i<j,\ {\rm then}\ f(k)\leq f(i)\ {\rm for\ all}\ i\leq k\leq j.$$
 
 Is is immediate to check that $f\times g\in {\mathcal K}_{n+m}$, for any $f\in {\mathcal K}_{n}$ and $g\in {\mathcal K}_{m}$.
 \ms
 
We extend the definition of irreducible permutation to $\displaystyle{\bigcup _{n\geq 1}{\mathcal F}_{n}}$ as follows:

\begin{defn} {\rm An element $f\in {\mathcal F}_{n}$ is called {\it irreducible} if 
$f\notin \displaystyle{\bigcup _{i=1}^{n-1}{\mathcal F}_{i}\times {\mathcal F}_{n-i}}$. 
The set of irreducible elements of ${\mathcal F}_{n}$ is denoted ${\mbox {\it Irr}_{{\mathcal F}_{n}}}$. In a similar way, the sets of irreducible elements of 
${\mathcal P}_{n}$ and ${\mathcal K}_{n}$ are the sets ${\mbox {\it Irr}_{{\mathcal P}_{n}}}:={\mathcal P}_{n}\cap {\mbox {\it Irr}_{{\mathcal F}_{n}}}$ and  
${\mbox {\it Irr}_{{\mathcal K}_{n}}}:={\mathcal K}_{n}\cap {\mbox {\it Irr}_{{\mathcal F}_{n}}}$, respectively. }
\end{defn}

Again, the graded space $K[{\mathcal F}_{\infty}]:=\displaystyle{\bigoplus _{n\geq 1}K[{\mathcal F}_{n}]}$ equipped with the concatenation product is the free 
associative algebra spanned by $\displaystyle{\bigcup _{n\geq 1}{\mbox {\it Irr}_{{\mathcal F}_{n}}}}$. Analogous results hold for the spaces 
$K[{\mathcal P}_{\infty}]:=\displaystyle{\bigoplus _{n\geq 1}K[{\mathcal P}_{n}]}$ 
and  $K[{\mathcal K}_{\infty}]:=\displaystyle{\bigoplus _{n\geq 1}K[{\mathcal K}_{n}]}$
\bs

\bs

\bs

\noi{\bf Planar rooted trees}
\ms

\begin{defn} {\rm A {\it planar rooted tree} is a non-empty oriented connected planar graph such that any vertex has at
least two input edges and one output edge, equipped with a final vertex called the {\it root}. For $n\geq 2$, a planar $n$-ary tree is a planar rooted tree such that 
any vertex has exactly $n$ input edges.} 
\end{defn}

\noi Note that in a planar tree the set of input edges of any vertex is totally ordered.

\noi All trees we deal with are reduced planar rooted ones. From now on, we shall use the term planar tree instead of planar rooted tree.

\begin{notn} {\rm \begin{enumerate}\item We denote by $Y_m$ the set of planar binary trees with $m+1$ leaves, and by $T_m$ the set of
 all planar trees with $m+1$ leaves. Clearly, $Y_m$ is a subset of $T_m$.
\item Given a tree $t\in T_m$, we denote by $v(t)$ the number of internal vertices of $t$ and by $\vert t\vert =m$ the degree of $t$.  The set $T_{m}$ is the 
disjoint union $\displaystyle{\bigcup _{r=1}^{m}T_{m}^{r}}$, where $T_{m}^{r}$ is the set of planar rooted trees with $m+1$ leaves and $r$ internal vertices.
\end{enumerate}}
\end{notn}

Let $t$ be an element of $T_m$, the leaves of $t$ are numbered from left to right, beginning with $0$ up to $m$. We denote by
$\cc _m$ the unique element of $T_{m}$, which has $m+1$ leaves and only one vertex (the corolla). For any $t\in Y_{m}$, the number of internal vertices 
of $t$ is $v(t)=m$. 

Let $X=\displaystyle{\bigcup _{n\geq 1}X_{n}}$ be a positively graded set. The set $T_{m,X}$ is the set of planar binary trees with the internal vertices coloured by the 
elements of $X$ in such a way that any vertex with $k$ input edges is coloured by an element of $X_{k-1}$. 

\begin{defn}\label{injerto} {\rm Given coloured trees $t$ and $w$, for any $0\leq i\leq \vert w\vert$, define $t\circ _iw$ to be the coloured tree obtained by attaching 
the root of $t$ to the $i$-th leaf of $w$.}
\end{defn} 

For instance $${\smallmatrix \searrow &&\swarrow \\ &x &\\&\downarrow &\endsmallmatrix } \circ _2{\smallmatrix \searrow &&\swarrow 
&&\\&y &&&\\ &&\searrow &&\swarrow \\ &&&z &\\ &&&\downarrow &\endsmallmatrix } = {\smallmatrix \searrow &&\swarrow &&\searrow&
&\swarrow\\&y &&&&x &\\&&\searrow &&\swarrow &&\\&&&z &&&\\&&&\downarrow&&&\endsmallmatrix } .$$

\begin{notn} {\rm \begin{enumerate}\item  Given coloured trees $t^0, t^1,\dots t^{\mid w\mid}$ and $w$, we denote by 

\noi $(t^0,\dots ,t^{\vert w\vert })\circ w$ the tree obtained as follows:
$$(t^0,\dots ,t^{\vert w\vert})\circ w:=t^0\circ _0 (t^1\circ _1 (\dots t^{\vert w\vert -1}\circ _{\vert w\vert -1}(t^{\vert w\vert }
\circ _{\vert w\vert } w))).$$
\item Given two coloured trees $t$ and $w$ and $x\in X_{1}$, we denote by $t\vee _{x} w$ the tree obtained by joining the roots of $t$ and $w$
to a new root, coloured by $x$. More generally, we denote by $\bigvee _{x}(t^0,\dots,t^r)$ the tree $(t^0,\dots,t^r)\circ (\cc _r,x)$, for $x\in X_{r}$.
\end{enumerate}}
\end{notn}
\ms

Any coloured tree $t$ may be written in a unique way as $t = \bigvee _{x}(t^0,\dots ,t^r)$, with $\vert t\vert  = \displaystyle{\sum _{i=0}^r\vert t^i\vert + r -1}$ 
and $x\in X_{r}$.
\ms

Recall that the number of elements of the set of planar binary rooted trees with $n+1$-leaves is the Catalan number 
$c_n=  {(2n)!\over n! (n+1)!}$ (see \cite{BB}). The number of planar rooted trees with $n+1$ leaves is called the super Catalan
number $C_m$, and it is given by the following recursive formula:
$$C_m= C_{m-1} + 2 \sum _{i=1}^{n-1}C_iC_{n-(i+1)}.$$
\bs

\bs

\bs

\section{Shuffle algebras}
\ms

Our goal is to describe the spaces spanned by coloured permutations and coloured elements of ${\mathcal P}_{\infty}$ as free objects 
for some type of algebraic structure.

In order to achieve our task we introduce the notion of shuffle algebras. 
\ms

\begin{defn} \label{SH} {\rm A {\it shuffle} algebra over $K$ is a graded $K$-vector space $A=\displaystyle{\bigoplus _{n\geq 0}A_{n}}$ equipped with 
linear maps
$$\bullet _{\gamma }:A_{n}\otimes _{K}A_{m}\rightarrow A,\ {\rm for}\ \gamma \in \Sh,$$
verifying that:}
 $$x\bullet _{\gamma }(y\bullet _{\delta }z)=(x\bullet _{\sigma}y)\bullet _{\lambda }z,$$
{\rm whenever } $(1_{n}\times \delta )\cdot \gamma= (\sigma \times 1_{r})\cdot \lambda\in {\mbox {\it Sh}(n,m,r)}$.
\end{defn}
\ms

Let ${\mbox{\it $\mathbb S$--Mod}}$ be the category whose objects are infinite families $M=\{ M(n)\}_{n \geq 0}$ of $K$-modules, such that each $M(n)$ is a right 
$K[S_n]$-module, for $n\geq 1$. 
A homomorphism $f$ from $M$ to $N$ in ${\mbox{\it $\mathbb S$--Mod}}$ is a family of linear maps 

\noi $f(n):M(n)\longrightarrow N(n)$ such that
 each $f(n)$ is a morphism of $K[S_n]$-modules.

\noi The category ${\mbox{\it $\mathbb S$--Mod}}$ is endowed with a symmetric monoidal structure $\otimes _{\mathbb S}$  given by:
$$({\mbox {\it M$\otimes _{\mathbb S}$N}})(n)=\bigoplus _{i=0}^n({\mbox {\it M$(i)\otimes$N$(n-i)$}})\otimes _{K[S_i\times S_{n-i}]}K[S_n],$$
where ${\mbox {\it M$(i)\otimes$N$(n-i)$}}$ has the natural structure of  right $K[{\mbox {\it S$_i\times $S$_{n-i}$}}]$-module.

By Proposition \ref{shuff}, the tensor product $({\mbox {\it M$(i)\otimes$N$(n-i)$}})\otimes _{K[S_i\times S_{n-i}]}K[S_n]$ is isomorphic to ${\mbox {\it M$(i)\otimes$N$(n-i)$}} \otimes K[{\mbox {\it Sh}(i,n-i)}]$. 

\noi Moreover, the associativity  and symmetry of  
$\otimes _{\mathbb S}$ are given by the isomorphisms:\begin{enumerate}
\item $a _{MNR}:({\mbox {\it M$\otimes _{\mathbb S}$N$)\otimes _{\mathbb S}$R}}\longrightarrow {\mbox {\it M$\otimes _{\mathbb S}($N$\otimes _{\mathbb S}$R}}),$
with $$a _{MNR}((x\otimes y\otimes \sigma)\otimes z\otimes \delta ):=x\otimes (y\otimes z\otimes \gamma )\otimes \tau ,$$ whenever 
$(\sigma \times 1_r)\cdot \delta=(1_n\times \gamma )\cdot \tau $ in ${\mbox {\it Sh}(m,n,r)}$, for ${\mbox {\it x}\in M(m)}$, $y\in N(n)$ and $z\in R(r)$.
\item $c_{MN}:{\mbox {\it M$\otimes _{\mathbb S}$N}}\longrightarrow {\mbox {\it N$\otimes _{\mathbb S}$M}},$ with 
$$c_{MN}(x\otimes y\otimes \sigma ):=y\otimes x\otimes (\epsilon _{n,m}\cdot \sigma ),$$
where $\epsilon _{n,m}=(n+1,\dots ,n+m,1,\dots ,n)\in \Sh$, for $x\in M(m)$ and $y\in N(n)$.
\end{enumerate}
\ms

Let $(M,\circ )$ be a monoid in $({\mbox{\it $\mathbb S$--Mod}}, \otimes _{\mathbb S})$,  the space $M=\displaystyle{\bigoplus _{n\geq 0}M(n)}$ has a 
natural structure of shuffle algebra, given by:
$$x\bullet _{\gamma }y:= \circ (x \otimes y\otimes \gamma ),$$
for $x\in M(n)$ and $y\in M(m)$.
The associativity of $\circ $ implies that the products $\bullet _{\gamma}$ fullfill the conditions of Definition \ref{SH}. 

\noi  In \cite{PS} and \cite{Liv}, an associative monoid in $({\mbox{\it $\mathbb S$--Mod}}, \otimes _{\mathbb S})$ is called a twisted associative algebra or 
an $As$-algebra in the category ${\mbox{\it $\mathbb S$--Mod}}$. 
\ms

\noi Given an associative graded algebra $(A=\displaystyle{\bigoplus _{n\geq 0}A_n},\bullet )$,  consider ${\overline A}=
\{ A_n\otimes K[S_n]\}_{n\geq 0}$. The ${\mathbb S}$-module ${\overline A}$ has a natural structure of monoid in $({\mbox{\it $\mathbb S$--Mod}}, 
\otimes _{\mathbb S})$, given by:
$$\circ ((x,\sigma )\otimes (y,\tau )\otimes \gamma):=(x\bullet y)\otimes (\sigma\times \tau)\cdot \gamma)\in A_{n+m}\otimes K[S_{n+m}],$$
for $x\in A_n$, $y\in A_m$, $\sigma \in S_n$, $\tau \in S_m$ and $\gamma \in \Sh$.
\ms

\begin{exmpls}\label{bish} {\rm {\bf a) The tensor space} For any vector space $V$ the tensor space ${\overline T}(V):=\bigoplus _{n\geq 1}V^{\otimes n}$, with the products $\bullet _{\gamma }$ given by:}
$$(v_1\otimes \dots \otimes v_n)\bullet _{\gamma }(v_{n+1}\otimes \dots \otimes v_{n+m}):=v_{\gamma (1)}\otimes \dots \otimes 
v_{\gamma (n+m)},$$
{\rm for $v_{1},\dots ,v_{n+m}\in V$, is a shuffle algebra.}
\ms

\noi {\rm {\bf b) Free shuffle algebras.} Define, on the graded vector space $K[S_{\infty }]:=\displaystyle{\bigoplus _{n\geq 1}K [S_{n}]}$, the operations 
$\bullet _{\gamma}:K[S_{n}]\otimes K [S_{m}]\longrightarrow K[S_{n+m}]$ as follows:}
$$\sigma \bullet _{\gamma}\tau :=(\sigma \times \tau)\cdot \gamma ,$$
{\rm for} $\sigma \in S_{n}$, $\tau \in S_{m}$ {\rm and} $\gamma \in \Sh$.
{\rm It is immediate to check that $(K[S_{\infty }],\bullet _{\gamma})$ is a shuffle algebra.
\ms

\noi In general, for any set $E$, the space $K [S_{\infty},E] :=\displaystyle{\bigoplus _{\geq 1}K[S_{n}\times E^{n}]}$ equipped with the operations 
$$(\sigma ,x_{1},\dots ,x_{n})\bullet _{\gamma}(\tau ,x_{n+1},\dots ,x_{n+m}):=(\sigma \bullet _{\gamma}\tau ,x_{\gamma (1)},\dots ,
x_{\gamma (n+m)})$$ 
is a shuffle algebra.}

{\rm Note that the map $E\hookrightarrow K[S_{\infty},E]$ maps $e\in E$ to the element $((1);e)\in S_{1}\times E$. So, the degree of any element $e\in E$ is one. }
\ms

{\rm In general, given a positively graded set $X=\displaystyle{\bigcup _{n\geq 1}X_{n}}$, consider the vector space $K [{\mathcal F}_{\infty},X]$ spanned by the elements 
$(f,x_{1},\dots ,x_{r})\in {\mathcal F}_{n}^{r}\times X^{r}$ such that $\vert x_{i}\vert =\vert f^{-1}(i)\vert $ for $1\leq i\leq r$, with the operations given by:
$$(f;x_{1},\dots ,x_{r}) \bullet _{\gamma}(g;y_{1},\dots ,y_{k}) :=((f \times g )\cdot \gamma ;x_{1},\dots ,x_{r},y_{1},\dots ,y_{k}),$$
for $(f;x_{1},\dots ,x_{r}) \in {\mathcal F}_{n,X}$, $(g;y_{1},\dots ,y_{k})\in {\mathcal F}_{m,X}$ and $\gamma \in \Sh$, is a shuffle algebra. }
 \ms
 
 \noi {\rm The subspace $K [{\mathcal P}_{\infty},X]$ of $K [{\mathcal F}_{\infty},X]$ is closed under the products $\bullet _{\gamma}$. 
 So, $K [{\mathcal P}_{\infty},X]$ has also a natural structure of shuffle algebra.}
\ms

\subsubsection {\bf Proposition.} Given a positively graded set $X$, the algebra $(K [{\mathcal P}_{\infty},X], \bullet _{\gamma })$ is the free shuffle 
algebra spanned by $X$.

\P  {\rm From the definition of shuffle algebra and Proposition \ref{shuff}, one has that any element in the free shuffle algebra spanned by $X$
is a sum of elements $x$, with}
$$x=x_{1}\bullet _{\gamma _{1}}(x_{2}\bullet _{\gamma _{2}}(\dots (x_{k-1}\bullet _{\gamma _{k-1}}x_{k}))),$$
{\rm for unique elements $x_{i}\in X$ and unique shuffles $\gamma _{i}$, for $1\leq i\leq k$.}
{\rm Let $\psi $ be the homomorphism from the free shuffle algebra spanned by $X$ to the space $K [{\mathcal P}_{\infty},X]$, such that:}
$$\psi (x_{1}\bullet _{\gamma _{1}}(x_{2}\bullet _{\gamma _{2}}(\dots (x_{k-1}\bullet _{\gamma _{k-1}}x_{k})))):= 
(\xi _{\n}\cdot \gamma; x_{1}, \dots , x_{k}), $$
where\begin{enumerate}\item $n_{i}=\vert x_{i}\vert $, for $1\leq i\leq k$,
\item $\gamma =(1_{n_{1}+\dots +n_{k-2}}\times \gamma _{k-1})\cdot \dots \cdot (1_{n_{1}}\times \gamma _{2})\cdot \gamma _{1}.$
\end{enumerate}
\ms

{\rm Conversely, let $f:\{1,\dots ,n\}\rightarrow \{ 1,\dots ,r\}$ be a surjective map, and let $n_{i}:=\vert f^{-1}(i)\vert $, for $1\leq i\leq r$. 
There exists a unique permutation $\gamma \in {\mbox {\it Sh}(\n)}$ such that $f=\xi _{\n}\cdot \gamma$.}

\noi {\rm Moreover, there exist unique permutations 

\noi $\gamma _{i}\in {\mbox {\it Sh}(n_{i},n_{i+1}+\dots +n_{k})}$ such that:}
$$\gamma = (1_{n_{1}+\dots +n_{k-2}}\times \gamma _{k})\dots \cdot\cdot(1_{n_{1}}\times \gamma _{2})\cdot \gamma _{1}.$$

\noi {\rm The inverse of $\psi $ is given by:} $$\displaylines {\psi ^{-1}(f; x_{1}, \dots , x_{k})=
x_{1}\bullet _{\gamma _{1}}(x_{2}\bullet _{\gamma _{2}}(\dots (x_{k-1}\bullet _{\gamma _{k-1}}x_{k}))).\hfill \diamondsuit }$$

{\rm Clearly, if $E$ is concentrated in degree $0$, the free shuffle algebra spanned by $E$ is $K [S_{\infty},E]$.}
\ms

\noi {\rm {\bf c) Non-unital infinitesimal bialgebras.}  Suppose that $(A,\cdot )$ is a graded $K$-algebra, equipped with a coassociative coproduct 
$\Delta :A\otimes A\rightarrow A$ such that:}
$$\Delta (x\cdot y)=\sum x\cdot y_{(1)}\otimes y_{(2)} + \sum x_{(1)}\otimes x_{(2)}\cdot y + x\otimes y,\ {\rm for}\ x,y\in A,$$
{\rm where $\Delta (z)=\sum z_{(1)}\otimes z_{(2)}$, for $z\in A$. The triple $(A,\cdot ,\Delta )$ is called a {\it nonunital infinitesimal bialgebra} (see \cite{LR3}).}
\ms  

{\rm It is easy to see that the reduced tensor space ${\overline T}(V)$, equipped with the concatenation product and the deconcatenation coproduct, 
is a graded unital infinitesimal bialgebra which is denoted ${\overline T}^{c}(V)$.} 

\subsubsection {\bf Remark.}\label{util} {\rm Given a permutation $\gamma \in \Sh$ there exists unique integers $n_{1},\dots ,n_{r}$ and $m_{1},\dots, m_{r}$ such that:
$$\gamma =(1,\dots ,n_{1},n+1,\dots n+m_{1},n_{1}+1,\dots ,n_{1}+n_{2},\dots ,m_{1}+\dots +m_{r-1}+1,\dots ,m),$$
where $\displaystyle{\sum _{i=1}^{r}n_i=n}$, $\displaystyle{\sum _{j=1}^{r}m_j=m}$, $n_{1}\geq 0$, $n_i\geq 1$ for $i>2$, $m_j\geq 1$ for 
$j<r$, and $m_{r}\geq 0$.}
\bs

{\rm Let $A=\displaystyle {\bigoplus _{n\geq 1}A_{n}}$ be a positively graded nonunital infinitesimal bialgebra. The map $\Delta _{n_{1},\dots ,n_{r}}:A_{n}\longrightarrow A_{n_{1}}\otimes \dots \otimes A_{n_{r}}$ is given by the composition of 
$\Delta ^{r-1}$ with the projection $p_{n_{1}\dots n_{r}}:A^{\otimes n}\longrightarrow A_{n_{1}}\otimes \dots \otimes A_{n_{r}}$. 

\noi For any $x\in A_{n}$, let 
$\Delta _{n_{1},\dots ,n_{r}}(x)=\sum x_{(1)}^{n_{1}}\otimes \dots \otimes x_{(r)}^{n_{r}}$.}
\ms

{\rm The proof of the following result is given for a general case in Theorem \ref{grftbi}.}

\subsubsection {\bf Lemma.} Let $(\displaystyle {\bigoplus _{n\geq 1}A_{n}},\cdot ,\Delta )$ be a graded nonunital infinitesimal bialgebra. 
The graded space $A$, equipped with the operations:
$$x\bullet _{\gamma}y=\sum _{\vert y_{(1)}\vert =i}x_{(1)}^{n_{1}}\cdot y_{(1)}^{m_{1}}\cdot x_{(2)}^{n_{2}}\cdot \dots \cdot y^{m_{r}}_{(r)},\ {\rm for}\ 
x\in A_{n},y\in A_{n},\ {\rm and}\ \gamma \in \Sh,$$
where $n_{1},\dots ,n_{r}$ and $m_{1},\dots ,m_{r}$ are the integers which determine $\gamma $, as pointed out in Remark \ref{util}
is a shuffle algebra.
\ms

\noi {\rm {\bf d) The algebra of parking functions.} (see \cite{NT1} and \cite{NT2}) Let ${\mbox{\it PF}_{n}}$ be the subset of all functions $f$ in 
${\mathcal F}_{n}^{n}$ which may be written as a composition $f=f^{\uparrow}\cdot \sigma $, with $f_{1}\in {\mathcal F}_{n}^{n}$ such that 
$f^{\uparrow}(i)\leq i$ for all $1\leq i\leq n$, and $\sigma \in S_{n}$.  Such a function is called a {\it parking} function. }
\ms

{\rm Applying Remark \ref{funfin}, we get that for any parking function $f\in {\mbox{\it PF}_{n}}$ there exist unique elements $f^{\uparrow}
\in {\mbox{\it PF}_{n}}$ and $\sigma \in {\mbox {\it Sh}(r_1,\dots ,r_n)}$ such that $f^{\uparrow}$ is a non-decreasing parking function and 
$f=f^{\uparrow}\cdot \sigma$, where} $r_i=\vert f^{-1}(i)\vert $.
\ms

{\rm Define the concatenation map $\times: {\mbox{\it PF}_{n}}\times {\mbox{\it PF}_{m}}\longrightarrow {\mbox{\it PF}_{n+m}}$  as the restriction of 
the concatenation product ${\mathcal F}_{n}^{n}\times {\mathcal F}_{m}^{m}\longrightarrow {\mathcal F}_{n+m}^{n+m}$, that is:

$$f\times g:=(f(1),\dots ,f(n),g(1)+n,\dots ,g(m)+n).$$ 

\noi Note that $f\times g$ is also a parking function.
 Moreover, for any functions $f\in {\mbox{\it PF}_{n}}$, $g\in {\mbox{\it PF}_{m}}$ and $\gamma \in \Sh$, the product 
 $f\bullet _{\gamma }g= (f\times g)\cdot \gamma $ belongs to ${\mbox{\it PF}_{n+m}}$. Let $\mbox{\bf {PQSym}}_{n}$ denote the 
 $K$-vector space spanned by the set ${\mbox{\it PF}_{n}}$ for $n\geq 1$, the space spanned by all parking 
 functions ${\mbox {\bf {PQSym}}}:=\displaystyle{\bigoplus _{n\geq 1}{\mbox{\bf {PQSym}}_n}}$ is a shuffle subalgebra of
  $K[{\mathcal F}_{\infty}]$.}
\ms

\noi {\rm Following {\bf 2.2.3} of \cite{NT2}, for a parking function $f\in {\mbox{\it PF}_{n}}$, an integer $b\in \{0,1,\dots ,n\}$ is called a {\it breakpoint } 
of $f$ if $\vert \{ i\mid F(i)\leq b\}\vert =b$. 
\ms

 I. Gessel defined a {\it primitive parking function} as an element  $f\in {\mbox{\it PF}_{n}}$ such that its unique breakpoints are 
the trivial ones: $0$ and $n$.  Let $\mbox {\it PPF}_n$ be the subset of prime parking functions of $\mbox{\it PF}_{n}$. It is immediate to check that 
$f\in {\mbox{\it PF}_{n}}$ if its associated non-decreasing parking function cannot be written a a concatenation of parking functions of smaller degree.}

\noi {\rm Note that the definition of breakpoint implies that if for any parking function $f\in P_{n}$
and any permutation $\sigma \in S_{n}$ the sets of breakpoints of $f$ and of $f\cdot \sigma$ are the same. So, the subset $\mbox {\it PPF}_{n}$ is invariant under the 
right action of $S_{n}$}. 

\subsubsection {\bf Remark.}\label{park} {\rm (see {\bf 2.2.3} of \cite{NT2}) A element in $\mbox {\it PPF}_n$ is a parking function which cannot be described as 
$f\bullet _{\gamma } g$ for some $f\in {\mbox {\it PF}_{k}}$, $g\in {\mbox{\it PF}_{n-k}}$ and $\gamma \in {\mbox {\it Sh}(k,n-k)}$.}
\ms

\noi {\rm Remark \ref{park} implies the following result.}

\subsubsection {\bf Lemma.}\label{parki} The shuffle algebra $\mbox{\bf {PQSym}}$ is the free shuffle algebra spanned by the set ${\mbox{\it PPF}}:=\displaystyle
{\bigcup _{n\geq 1}{\mbox {\it PPF}_n}}$ of all prime parking functions.
\ms

\P {\rm As is pointed out in \cite{NT2}, a parking function $f\in {\mbox{\it PF}_{n}}$ has a breakpoint at $0<b<n$ if and only if there exist unique functions 
$f_{1}\in {\mbox{\it PF}_{b}}$, 
$f_{2}\in {\mbox{\it PF}_{n-b}}$ and a shuffle $\sigma \in {\mbox {\it Sh}(b,n-b)}$ (not necesarily unique) such that $f=(f_{1}\times f_{2})\cdot \sigma$. 
By a recursive argument on the number of breakpoints of $f$, it is immediate to check that the set ${\mbox{\it PPF}}:=\displaystyle {\bigcup _{n\geq 1}
{\mbox {\it PPF}_n}}$ spans ${\bold {PQSym}}$ as a shuffle algebra.}

\noi{\rm To see that ${\bf {PQSym}}$ is free as a shuffle algebra, it suffices to note that for any function $f\in {\mbox{\it PF}_{n}}$ with breakpoints $0<b_1<\dots 
<b_r<n$, there exist unique elements $f_1\in {\mbox{\it PPF}_{b_1}},f_2\in {\mbox{\it PPF}_{b_2-b_1}}, \dots , f_{r+1}\in {\mbox{\it PPF}_{n-b_r}}$ such that} 
$$f=(f_1\times f_2\times\dots 
\times f_{r+1})\cdot \sigma ,\ {\rm with}\ \sigma \in {\mbox {\it Sh}(b_{1},b_{2}-b_{1},\dots ,n-b_{r})}.\hfill \diamondsuit $$
\ms

\noi {\rm Note that the group $S_{n}$ acts on the right on the set $\mbox {\it PPF}_{n}$, for $n\geq 1$. 

\noi So, ${\mbox {\it PPF}}=\{ K[{\mbox {\it PPF}_{n}}]\}_{n\geq 1}$ is an object in the category 
\mbox{$\mathbb S$--Mod}. Applying Lemma \ref{parki}, it is immediate to check that ${\mbox {\bf PQSym}}=
{\overline T}_{\mathbb S}({\mbox {\it PPF}})=
\displaystyle{\bigoplus _{n\geq 1}{\mbox {\it PPF}}^{\otimes _{\mathbb S}n}}$ in the category \mbox{$\mathbb S$--Mod}, which means that as a shuffle algebra 
${\bf {PQSym}}$ is the free monoid 
spanned by $\mbox {\it PPF}$ in the monoidal category $(\mbox{$\mathbb S$--Mod} ,\otimes _{\mathbb S})$.}

\end{exmpls}
\bs

\bs

\begin{defn}\label{prods}{\rm Given graded spaces $V=\displaystyle{\bigoplus _{n\geq 0}V_{n}}$ and $W=\displaystyle{\bigoplus _{m\geq 0}W_{m}}$ 
there exist two ways to obtain the product of both spaces:\begin{enumerate}
\item The {\it Hadamard product} of $V$ and $W$, denoted by $V\underset {H}{\otimes} W$, is the graded vector space such that 
$(V\underset {H}{\otimes} W)_{n}:= V_{n}\otimes W_{n}$, for $n\geq 0$.
\item The {\it tensor product} of $V$ and $W$, denoted by $V\otimes W$, is the graded vector space such that $(V\otimes W)_{n}:= \displaystyle{\bigoplus _{i=0}^{n}
V_{i}\otimes W_{n-i}}$, for $n\geq 0$.
\end{enumerate}}\end{defn}
\ms

\begin{rem}\label{diagonal} {\rm Let $0\leq r\leq n+m$ be an integer and let $\gamma $ be a ${\mbox {\it (n,m)}}$-shuffle. There exist a unique non negative integer
$0\leq n_{1}\leq r$ and permutations $\gamma _{(1)}^{r}\in S_{r}$ and $\gamma _{(2)}^{n+m-r}\in S_{n+m-r}$ such that 
$\gamma =(1_{n_{1}}\times \epsilon _{n-n_{1},m_{1}}\times 1_{m-m_1})\cdot (\gamma _{(1)}^{r}\times \gamma _{(2)}^{n+m-r})$,
where $n_1:=\vert \gamma ^{-1}(\{1,\dots ,n\})\cap \{1,\dots r\}\vert $ and $m_1:=r-n_1$. Moreover, the permutation $\gamma _{(1)}^{r}$ 
belongs to ${\mbox{\it Sh}(n_{1},m_{1})}$ and $\gamma _{(2)}^{n+m-r}$ belongs to ${\mbox{\it Sh}(n-n_{1},m-m_1)}$.}
\end{rem}
\ms

\noi The proof of the following result is immediate.

\begin{lem}\label{hadsh}  Let $(A, \bullet _{\gamma })$ and $(B,\circ _{\delta })$ be two shuffle algebras.
\begin{enumerate}\item The Hadamard product $A\underset {H}{\otimes}B$ has a natural structure of shuffle algebra, given by the operations:
$$(x\otimes  y)\bullet _{\gamma}(x\rq \otimes y\rq ):= (x\bullet _{\gamma }x\rq )\otimes (y\circ _{\gamma }y\rq),$$ 
for $x\in A_{n}, y\in B_{n}, x\rq \in A_{m}, y\rq \in B_{m}$ and $\gamma \in \Sh$.
\item The tensor product $A\otimes B$  has a natural structure of shuffle algebra, given by the operations:
$$(x\otimes y)\bullet _{\gamma }(x\rq \otimes y\rq ):=\begin{cases} 
(x\bullet _{\gamma _{(1)}^{n+n\rq}}x\rq )\otimes (y\circ _{\gamma _{(2)}^{m+m\rq}}y\rq ),& \ {\rm for}\ 
\ n = (n+n\rq )_{1},\\
0,&\ {\rm otherwise,}\end{cases}$$
where $x\in A_{n}$, $x\rq \in A_{n\rq }$, $y\in A_{m}$, $y\rq \in A_{m\rq}$, $\gamma _{(1)}^{n+n\rq }\in {\mbox {\it Sh}(n,n\rq )}$ and $\gamma _{(2)}^{m+m\rq }\in 
{\mbox {\it Sh}(m,m\rq )}$ are the permutations defined in Remark \ref{diagonal}, and 

\noi $(n+n\rq )_1:=\vert \gamma ^{-1}(\{1,\dots n+m\}\cap \{1,\dots ,n+n\rq \}\vert $.
\end{enumerate}
\end{lem}
\bs

\noi For any shuffle algebra $(A,\bullet _{\gamma})$, define the products $\bullet _{0}$ and $\bullet _{top}$ as follows:
$$\displaylines {
x\bullet _{0}y:=x\bullet _{1_{n+m}}y,\hfill \cr
x\bullet _{top} y:=y\bullet _{\epsilon _{m,n} } x,\hfill \cr }$$
for $x\in A_{n}$ and $y\in A_{m}$.

\begin{rem}\label{assprod} {\rm The products $\bullet _0$ and $\bullet _{top}$ are associative.}
\end{rem}
\ms

\noi Moreover, Proposition \ref{shuff} implies the following result.
\ms

\begin{lem}\label{shuffprod} Let $(A,\bullet _{\gamma})$ be a shuffle algebra. The product $*:A\otimes A\rightarrow A$ defined as:
$$x * y:=\sum _{\gamma \in \Sh}x\bullet _{\gamma }y,\ {\rm for}\ x\in A_{n}\ {\rm and}\ y\in A_{m},$$
is associative.
\end{lem}
\bs

\bs

\bs

\section{Shuffle bialgebras.}
\ms

We study coproducts on shuffle algebras that turn them into Hopf algebras.
\ms

\begin{defn}\label{defbi} {\rm Let $(A,\bullet _{\gamma })$ be a positively graded shuffle algebra, such that $A$ is equipped with a graded coassociative coproduct $\Delta $. 
We say that $(A,\bullet _{\gamma } , \Delta  )$ is a {\it shuffle bialgebra} if it verifies:
$$\Delta (x\bullet _{\gamma } y)=\sum _{r=1}^{n+m-1}\bigl (\sum (x_{(1)}\bullet _{\gamma _{(1)}^{r}}y_{(1)})\otimes (x_{(2)}
\bullet _{\gamma _{(2)}^{n+m-r}} y_{(2)})\bigr ),$$
where  $\gamma _{(1)}^{r}$ and $\gamma _{(2)}^{n+m-r}$ are defined in Remark \ref{diagonal}, the second sum is taken over all $\vert x_{(1)}\vert =n_{1}$ and 
$\vert y_{(1)}\vert =m_{1}$, and we fix 
$$x_{(1)}\bullet _{\gamma _{(1)}^{r}}y_{(1)}:=\begin{cases}x,&{\rm for}\ n_{1}=n\\
y,&{\rm for}\ n_{1}=0,\end{cases}$$
$$x_{(2)}\bullet _{\gamma _{(2)}^{n+m-r}} y_{(2)}:=\begin{cases}x,&{\rm for}\ n_{1}=0\\
y,&{\rm for}\ n_{1}=n,\end{cases}$$ }
\end{defn}
\ms

\begin{prop} Let $(A,\bullet _{\gamma }, \Delta _{A})$ and $(B,\circ _{\delta }, \Delta_{B} )$ be shuffle bialgebras. The Hadamard product 
$A\underset {H}{\otimes}B$ with the operations $\bullet _{\gamma}$ given in Lemma \ref{hadsh} and the coproduct given by:
$$\Delta _{A\underset {H}{\otimes}B}(x\otimes y)=\sum _{\vert x_{(1)}\vert =\vert y_{(1)}\vert }(x_{(1)}\otimes y_{(1)})\otimes (x_{(2)}\otimes y_{(2)}),$$ 
is a shuffle bialgebra.
\end{prop}
\ms

\P  Let $x\in A_{n}$, $y\in B_{n}$, $z \in A_{m}$, $w \in B_{m}$ and $\gamma \in \Sh$. 

\noi For $1\leq r,s \leq n+m$, we have that 
$\vert x_{(1)}^{n_{1}}\bullet _{\gamma _{(1)}^{r}}z _{(2)}^{m_{1}}\vert =r$ and $\vert y_{(1)}^{k_{1}}\bullet _{\gamma _{(1)}^{s}}w _{(2)}^{l_{1}}
\vert =s$, for $n_1+m_1 =r$ and $k_1+l_1=s$. So,
$\vert x_{(1)}^{n_{1}}\bullet _{\gamma _{(1)}^{r}}z _{(2)}^{m_{1}}\vert =\vert y_{(1)}^{k_{1}}\bullet _{\gamma _{(1)}^{s}}w _{(2)}^{l_{1}}\vert $
if, and only if $r=s$, which implies that 

 $$\vert x_{(1)}^{n_{1}}\vert =n_{1}=\vert y_{(1)}^{k_{1}}\vert\ {\rm and}\ \vert z _{(2)}^{m_{1}}\vert =m_{1}
=\vert w _{(2)}^{l_{1}}\vert .$$

\noi The argument above implies the result.\hfill $\diamondsuit$
\bs

The proof of the following Lemma is straightforward.

\begin{lem} $(1)$ Let $(A,\bullet _{i},\Delta )$ be a shuffle bialgebra. The relationship between $\Delta $ and the associative products 
$\bullet _{0}$ and $\bullet _{top}$ defined in the previous subsection, is given by the following equalities:
\begin{align}
\Delta (x\bullet _{0}y)&=\sum (x\bullet _{0}y_{(1)})\otimes y_{(2)} +\sum x_{(1)}\otimes (x_{(2)}\bullet _{0}y)+ x\otimes y,\notag\\
\Delta (x\bullet _{top}y)&=\sum (x\bullet _{top}y_{(1)})\otimes y_{(2)}+\sum x_{(1)}\otimes (x_{(2)}\bullet _{top}y) + x\otimes y,\notag\end{align} 
for $x,y\in A$.

\noi $(2)$ Let $(A,\bullet _{\gamma },\Delta )$ be a shuffle bialgebra and let $*$ be the associative product defined on $A_{+}=A\oplus K$ by:
$$x*y=\begin{cases}
\displaystyle{\sum _{\gamma \in \Sh}x\bullet _{\gamma }y},&\ {\rm for}\ x\in A_{n}\ {\rm and}\ y\in A_{m}\\
\qquad xy,&\ {\rm if}\ x\in K\ {\rm or}\ y\in K .\end{cases}$$
Then $(A_{+},*,\Delta _{+})$ is a Hopf in the usual sense, which means that:
$$\Delta _{+}(x*y)=\sum (x_{(1)}*y_{(1)})\otimes (x_{(2)}*y_{(2)}),\ {\rm for}\ x,y\in A_{+},$$
where $\Delta _{+}(\lambda )=\lambda 1_{K}\otimes 1_{K}$, for $\lambda \in K$, and $\Delta _{+}(x):=x\otimes 1_{K}+1_{K}\otimes x
+\Delta (x)$, for $x\in A$.
\end{lem}
\ms

\begin{cor}\label{hs} $(1)$ If $(A,\bullet _{i},\Delta )$ is a shuffle bialgebra, then $(A,\bullet _{0},\Delta )$ and 
$(A, \bullet _{top},\Delta )$ are nonunital infinitesimal bialgebras. 

\noi $(2)$ If $(A,\bullet _{\gamma },\Delta )$ is a shuffle bialgebra, then $(A_{+},*,\Delta _{+})$ is a Hopf algebra.
\end{cor}

The previous result implies that there exists two functors, $H_{0}$ and $H_{top}$, from the category of shuffle bialgebras to the category of graded  
nonunital infinitesimal bialgebras.
\bs

\bs

We prove that all the examples of shuffle algebras given in the previous Section can be equipped with a structure of shuffle bialgebra.

\begin{exmpls}\label{exbial} {\rm {\bf a) The Malvenuto-Reutenauer bialgebra} (see \cite{MR}) On the vector space $K[S_{\infty}]$, let $\Delta _{MR}$ be the 
unique coproduct such that:}
$$\Delta _{MR}(\sigma ):= \sum _{r=1}^{n-1}\sigma _{(1)}^{r}\otimes \sigma _{(2)}^{n-r},$$
{\rm for} $\sigma \in S_{n}${\rm , where} $\sigma =\delta _{r}\cdot (\sigma _{(1)}^{r}\times \sigma _{(2)}^{n-r})$, with $\delta _{r}^{-1}\in 
{\mbox {\it Sh}(r,n-r)}${\rm , for} $1\leq r\leq n-1$.
\ms

\subsubsection {\bf Proposition} \label{propiedades} The shuffle algebra  $(K[S_{\infty}], \bullet _{\gamma })$, equipped with the coproduct $\Delta _{MR}$ 
is a shuffle bialgebra. 
\ms

\P {\rm Let $\gamma $ be a {\mbox {\it (n,m)}}-shuffle, and let $0\leq r\leq n+m$ be an integer. }

\noi {\rm From Remark \ref{diagonal}, we have that there exists $0\leq n_1,m_1\leq r$ such that:} 
  $$\gamma  = (1_{n_1}\times \epsilon _{ n-n_1,m_1}\times 1_{m-m_1})\cdot (\gamma _{(1)}^r\times \gamma _{(2)}^{n+m-r}).$$
{\rm Suppose that, for $\sigma \in S_n$ and $\tau \in S_m$,}
$$\sigma =\alpha \cdot (\sigma _{(1)}^{n_1}\times \sigma _{(2)}^{n-n_1})\ {\rm and}\ \tau =\beta \cdot (\tau _{(1)}^{m_1}\times \tau _{(2)}
^{m-m_1}).$$

\noi {\rm Note that} $$(\sigma _{(2)}^{n-n_1}\times \tau _{(1)}^{m_1})\cdot \epsilon_{ n-n_1,m_1}=
 \epsilon _{n-n_1,m_1}\cdot (\tau _{(1)}^{m_1}\times \sigma _{(2)}^{n-n_1}).$$

\noi {\rm So, the following equality holds:}
$$\displaylines {
\sigma \bullet _{\gamma }\tau =\hfill \cr 
 (\alpha \times \beta )\cdot (\sigma_{(1)}^{n_1}\times  \epsilon _{ n-n_1,m_1})\cdot (\tau _{(1)}^{m_1}
 \times \sigma _{(2)}^{n-n_1})\times \tau _{(2)}^{m-m_1})\cdot (\gamma_{(1)}^r\times \gamma _{(2)}^{n+m-r})=\hfill \cr
\hfill (\alpha \times \beta )\cdot (1_{n_1}\times  \epsilon _{ n-n_1,m_1}\times 1_{m-m_1})\cdot ((\sigma _{(1)}^{n_1}
\bullet _{\gamma _{(1)}^r} \tau _{(1)}^{m_1})\times (\sigma _{(2)}^{n-n_1}
 \bullet _{\gamma _{(2)}^{n+m-r}}\tau _{(2)}^{m-m_1})).\cr } $$
 
{\rm To end the proof it suffices to note that} $ (\alpha \times \beta )\cdot (1_{n_{1}}\times  \epsilon _{m_1,n-n_1}\times 1_{m-m_1})$ 
{\rm belongs to} ${\mbox {\it Sh}(r,n+m-r)}$.\hfill  $\diamondsuit $
\bs

{\rm Let $E$ be a set, the coproduct of $K[S_{\infty}]$ extends to $K[S_{\infty},E]$ as follows:}
$$\Delta _{MR}(\sigma ,e_{1},\dots ,e_{n})=\sum _{i=1}^{n-1}(\sigma _{(1)}^{i},e_{\delta _{i}(1)},\dots ,e_{\delta _{i}(i)})\otimes (\sigma _{(2)}^{n-i},
e_{\delta _{i}(i+1)},\dots ,e_{\delta _{i}(n)}),$$
{\rm where} $\sigma =\delta _{i}\cdot (\sigma _{(1)}^{i}\times \sigma _{(2)}^{n-1})${\rm , with} $\delta _{i}\in Sh(i,n-i)$ for $1\leq i\leq n-1$.
\ms

\noi {\rm Using the Proposition above, is not difficult to check that $(K[S_{\infty},E],\bullet _{i},\Delta _{MR})$ is a shuffle bialgebra.}
\bs

\bs

{\rm {\bf b) Nonunital infinitesimal bialgebras. } Let $(A,\cdot ,\Delta )$ be a graded nonunital infinitesimal bialgebra.} 

\subsubsection {\bf Lemma } The associated shuffle algebra $(A,\bullet _{\gamma})$, with the coproduct $\Delta $ is a shuffle bialgebra.
\ms

\P  {\rm Let $\gamma \in \Sh$ be the permutation given by the sequences $n_{1},\dots ,n_{r}$ and 
$m_{1},\dots ,m_{r}$, as described in Remark \ref{util}. 

\noi For any $1\leq s\leq n+m-1$, note that the decomposition 
$$\gamma =( 1_{p_{1}}\times \epsilon _{n-p_1,q_1}\times 1_{m-q_1})\cdot 
(\gamma _{(1)}^{s}\times \gamma _{(2)}^{n+m-s})$$
of Remark \ref{diagonal}, is such that there exists $1\leq k\leq r$, and $1\leq n\rq _{k}\leq n_{k}$ or $1\leq m\rq _{k}\leq m_{k}$ with:}
$$p_{1}=\begin{cases}n_{1}+\dots +n\rq_{k},&\ {\rm if}\ s=\displaystyle{\sum _{i=1}^{k-1}(n_{i}+m_{i})+n\rq _{k}}\\
n_{1}+\dots +n_{k},&\ {\rm if}\ s=\displaystyle{\sum _{i=1}^{k-1}(n_{i}+m_{i})+n_{k}+m\rq _{k}}.\end{cases}$$
{\rm In the first case $\gamma _{(1)}^{s}$ is given by the sequence $n_{1},\dots ,n_{k-1},n\rq _{k}$ and $m_{1},\dots ,m_{k-1},0$, while in the second one 
$\gamma _{(1)}^{s}$ is given by the sequence $n_{1},\dots ,n _{k}$ and $m_{1},\dots ,m_{k-1},m\rq _{k}$. }

\noi {\rm Given elements $x\in A_n$, $y\in A_m$, the coassociativity of $\Delta$ and the relation between $\cdot$ and $\Delta$ state that:}
$$\displaylines {
\Delta (x\bullet _{\gamma}y)=\sum \Delta (x_{(1)}^{n_{1}}\cdot y_{(1)}^{m_{1}}\cdot \dots \cdot x_{(r)}^{n_{r}}\cdot y_{(r)}^{m_{r}})=\hfill \cr
\sum  (\sum _{1\leq k\leq r}(\sum _{1\leq n\rq _{k}\leq n_{k}}(x_{(1)}^{n_{1}}\cdot y_{(1)}^{m_{1}}\cdot \dots \cdot x_{(k)}^{n\rq _{k}})
\otimes (x_{(k+1)}^{n_{k}-n\rq _{k}}\cdot y_{(k)}^{n_{k}}\cdot \dots \cdot y_{(r)}^{m_{r}}))+\hfill\cr
\hfill (\sum _{1\leq m\rq _{k}\leq m_{k}}(x_{(1)}^{n_{1}}\cdot \dots \cdot x_{(k)}^{n_{k}}\cdot y_{k}^{m\rq _{k}})\otimes (y_{(k+1)}^{m_{k}-m
\rq _{k}}\cdot x_{(k+1)}^{n_{k+1}}\cdot \dots \cdot y_{(r)}^{m_{r}})) ).\cr }$$

\noi {\rm But} $(x_{(1)}\bullet _{\gamma _{(1)}^{s}}y_{(1)})\otimes (x_{(2)}\bullet _{\gamma _{(2)}^{n+m-s}}y_{(2)})=$
$$\begin{cases}(x_{(1)}^{n_{1}}\cdot y_{(1)}^{m_{1}}\dots x_{(k)}^{n\rq _{k}})
\otimes (x_{(k+1)}^{n_{k}-n\rq _{k}}\cdot y_{(k)}^{n_{k}}\dots y_{(r)}^{m_{r}}),&{\rm for}\ s=\displaystyle{\sum _{i=1}^{k-1}(n_{i}+m_{i})+n\rq _{k}}\\
(x_{(1)}^{n_{1}} \dots x_{(k)}^{n_{k}}\cdot y_{(k)}^{m\rq _{k}})\otimes (y_{(k+1)}^{m_{k}-m
\rq _{k}}\cdot x_{(k+1)}^{n_{k+1}}\cdot \dots \cdot y_{(r)}^{m_{r}}),&{\rm for}\ s=\displaystyle{\sum _{i=1}^{k-1}(n_{i}+m_{i})+n_{k}+m\rq _{k}},\end{cases}$$
{\rm which implies the result.}\hfill $\diamondsuit $
\ms

{\rm Let $G: {\mbox {\it Gr}\epsilon} \longrightarrow {\mbox {\it Sh}} $ be the functor which assigns to any graded nonunital infinitesimal bialgebra $(A,\cdot ,\Delta )$ the 
shuffle bialgebra $(A,\bullet _{\gamma},\Delta )$. It is easy to check that the compositions $H_{0}\circ G$ and $H_{top}\circ G$ are the identity functor, where $H_{0}$ 
and $H_{top}$ are the functors defined in Corollary \ref{hs}.}
\bs

\noi{\rm {\bf c) The free shuffle algebra over a graded set.} Let $X$ be a positively graded set, 
and let $\Theta :K[X]\longrightarrow K[X]\otimes K[X]$, be a graded coassociative coproduct on $K[X]$.

\noi The coproduct $\Delta _{\theta}$ on $K[{\mathcal P}_{\infty},X]$ is defined as follows:

\noi given $f=\xi_{\n}\cdot \sigma$, with $\n=(n_1,\dots, n_r)$, $\sigma \in {\mbox {\it Sh}(\n)}$, and elements $x_{1},\dots ,x_{r}\in X$, with $x_i\in X_{n_i}$:}
$$\displaylines {\Delta _{\theta}(f;x_{1},\dots ,x_{r}):=\hfill\cr
\hfill\sum _{i=0}^{n}\bigl(\sum _{\vert x_{j(1)}\vert =m_{j}^{i}}(\xi _{\m^i}\cdot \sigma _{(1)}^{i};x_{1(1)},\dots ,
x_{r(1)})\otimes (\xi _{\n-\m^i}\cdot \sigma _{(2)}^{n-i};x_{1(2)},\dots ,x_{r(2)})\bigr ),\cr }$$
{\rm where}\begin{enumerate}
\item $\sigma =\delta _{i}\cdot (\sigma _{(1)}^{i}\otimes \sigma _{(2)}^{n-i})$, with $\delta _{i}\in {\mbox {\it Sh}(i,n-i)}$,
\item {\rm for each} $1\leq i\leq n-1$, $$m_{j}^{i}:=\vert \delta _{i}^{-1}\{ 1,\dots ,i\}\cap \{ n_{1}+\dots +n_{j-1}+1,\dots ,n_{1}+\dots +n_{j}\}\vert,$$
\item $\m^i:=(m_1^i,\dots ,m_r^i)$ and $\n-\m^i:=(n_1-m_1^i,\dots ,n_r-m_r^i),$
\item $\Theta (x_{j})=\sum x_{j(1)}\otimes x_{j(2)}$, for $1\leq j\leq r$.
\end{enumerate}
\medskip

{\rm For example, suppose that $X_n=\{\xi _n\}$, for $n\geq 1$, if $\Theta $ is the unique coassociative coproduct on $K[X]$ such that $\Theta (\xi _n)={\displaystyle\sum _{i=0}^n\xi _i\otimes \xi_{n-i}}$ then for $f=(2,3,3,5,4,1,4,3)$, we get that:}
$$\displaylines {
\Delta_{\theta} (f)= (1)\otimes (2,2,5,3,1,3,2) + (1,2)\otimes (2,4,3,1,3,2) + \hfill\cr
(1,2,2)\otimes (4,3,1,3,2) +(1,2,2,3)\otimes (3,1,3,2) +(1,2,2,4,3)\otimes (1,3,2) +\cr
\hfill (2,3,3,5,4,1)\otimes (2,1) + (2,3,3,5,4,1,4)\otimes (1).\cr }$$
\ms

 \subsubsection {\bf Proposition.} For any positively graded set $X$ and any coassociative coproduct $\Theta $ defined on $K[X]$, the coproduct $\Delta _{\theta}$ 
 defines a shuffle bialgebra structure on $(K[{\mathcal P}_{\infty},X],\bullet _{\gamma })$. 
 \ms

 \P {\rm Let} $f=\xi _{\n}\cdot \sigma \in {\mathcal P}_{n}^{r}$, $g=\xi _{\m}\cdot \tau \in {\mathcal P}_{m}^{k}$, 
 $x_{1},\dots ,x_{r},y_{1}, \dots ,y_{k}\in X$ {\rm and } $\gamma \in \Sh${\rm , be such that} $\sigma \in {\mbox {\it Sh}(\n)}$, 
 $\tau \in {\mbox {\it Sh}(\m)}$,  $\vert x_{i}\vert = n_{i}$ and $\vert y_{j}\vert =m_{j}${\rm , for }$1\leq i\leq r$ and $1\leq j\leq k$.
\ms

\noi {\rm Suppose that, for} $0\leq r\leq n+m$, 
$\gamma =\delta _{r}\cdot (\gamma _{(1)}^{r}\times \gamma _{(2)}^{n+m-r})$, with $$\delta _{r}= 1_{n_1}\times \epsilon_{n-n_1,m_1}\times 1_{m-m_1}\in 
{\mbox {\it Sh}(r,n+m-r)}.$$
 
 \noi {\rm In Example {\bf a)} we prove that:}
 $$ \sigma \bullet _{\gamma}\tau=(\alpha \times \beta )\cdot \delta _{r}\cdot ((\sigma _{(1)}^{n_1}\bullet _{\gamma _{(1)}^{r}} 
 \tau _{(1)}^{m_1})\times (\sigma _{(2)}^{n-n_1}\bullet _{\gamma _{(2)}^{n+m-r}}\tau _{(2)}^{m-m_1})),$$
{\rm with }$\alpha \in {\mbox {\it Sh}(n_1,n-n_1)}$ {\rm and} $\beta \in {\mbox {\it Sh}(m_1,m-m_1)}$. {\rm So, }
 $$\sigma \bullet _{\gamma}\tau = ((\alpha \times \beta )\cdot \delta _{i})\cdot ((\sigma \bullet _{\gamma }\tau)_{(1)}^{r}\times 
 (\sigma \bullet _{\gamma}\tau)_{(2)}^{n+m-r}).$$
{\rm We have that: }$$\displaylines{
\Delta _{\theta}((f;x_{1},\dots ,x_{r})\bullet _{\gamma }(g;y_{1},\dots ,y_{k}))=\hfill\cr
\sum _{i}\bigl (\sum _{\mid x_{(1)j}\mid =l_{j}\atop 
\mid y_{(1)j}\mid =s_{j}}(\xi _{{\underline{\bold  l}},{\underline {\bold h}}}\cdot (\sigma \bullet _{\gamma }\tau )_{(1)}^{i};x_{1(1)},\dots ,x_{r(1)},y_{1(1)},
\dots ,y_{k(1)})\otimes\hfill \cr
\hfill  (\xi _{{\underline {{\bold n}-{\bold l}}} ,{\underline {{\bold m}-{\bold h}}}}\cdot (\sigma \bullet _{\gamma}\tau)_{(2)}^{n+m-i};x_{1(2)},\dots ,x_{r(2)},
y_{1(2)},\dots ,y_{k(2)})\bigr ),\cr}$$
{\rm where} $l_{j}:=\vert \alpha _{p_{i}}^{-1}\{ 1,\dots ,p_{i}\}\cap \{ n_{1}+\dots +n_{j-1}+1,\dots ,n_{1}+\dots +n_{j}\}\vert$, 

\noi $h_{j}:=\vert \beta _{i-p_{i}}^{-1}\{ 1,\dots ,i-p_{i}\}\cap \{ m_{1}+\dots +m_{j-1}+1,\dots ,m_{1}+\dots ,m_{j}\}\vert $, ${\underline {\bold l}}=(l_1,\dots , l_r)$ and
${\underline {\bold h}}=(h_1,\dots ,h_k)$. 
\ms

\noi {\rm So, we get that } $ \Delta _{\theta }((f;x_{1},\dots ,x_{r})\bullet _{\gamma }(g;y_{1},\dots ,y_{k}))=$ 
$$\sum _{i} (f;x_{1},\dots ,x_{r})_{(1)}
\bullet _{\gamma _{(1)}^{i}}(g;y_{1},\dots ,y_{k})_{(1)}\otimes (f;x_{1},\dots ,x_{r})_{(2)}\bullet _{\gamma _{(2)}^{n+m-i}}
(g;y_{1},\dots ,y_{k})_{(2)}.\hfill \diamondsuit $$ 
\bs

\noi {\rm {\bf d) Monoids in $({\mbox{\it $\mathbb S$--Mod}},\otimes _{\mathbb S})$}. For an ${\mathbb S}$-module $M$, a coproduct on $M$ is a family of homomorphisms of 
$K[S_n]$-modules $\Omega _n: M(n)\longrightarrow \displaystyle {\bigoplus_{i=0}^n M(i)\otimes M(n-i)\otimes K[{\mbox {\it Sh}(i,n-i)}]}$, for each $n\geq 0$. 
For $x\in M(n)$, we have that $$\Omega (x)=\sum _{i=0}^n(\sum _{\sigma \in {\mbox {\it Sh}(i,n-i)}}x_{(1)}^{\sigma}\otimes x_{(2)}^{\sigma }\otimes \sigma).$$
The coproduct $\Omega $ is coassociative if for any $\sigma \in {\mbox {\it Sh}(n,m+r)}$, $\tau \in {\mbox {\it Sh}(m,r)}$, 

\noi $\delta \in {\mbox {\it Sh}(n+m,r)}$ and 
$\omega \in {\mbox {\it Sh}(n,m)}$, such that $(1_n\times \tau)\cdot \sigma =(\omega \times 1_r)\cdot \delta $, it verifies the equality:}
$$\sum x_{(1)}^{\sigma}\otimes (x_{(2)}^{\sigma })_{(1)}^{\tau}\otimes (x_{(2)}^{\sigma })_{(2)}^{\tau }=\sum (x_{(1)}^{\delta })_{(1)}^{\omega}\otimes (x_{(1)}^{\delta })_{(2)}^{\omega }
\otimes x_{(2)}^{\delta }.$$
\medskip

{\rm A monoid $(M, \circ)$ in the category $({\mbox{\it $\mathbb S$--Mod}},\otimes _{\mathbb S})$ is  a bialgebra if it is equipped with a coassociative coproduct verifying the condition:
$$\Omega (\circ (x\otimes y\otimes \gamma))=
\sum (x_{(1)}^{\delta}\otimes y_{(1)}^{\tau}\otimes \alpha _{1})\otimes (x_{(2)}^{\delta}\otimes 
y_{(2)}^{\tau}\otimes \alpha _{2})\otimes \rho ,\eqno (*)$$
for $x\in M(n)$, $y\in M(m)$ and $\gamma\in {\mbox {\it Sh}(n,m)}$, where $$({\mbox 1_{n_{1}}\times  \epsilon _{m_1,n_2}\times 1_{m_2}})\cdot (\delta \times \tau)\cdot \gamma= (\alpha _1\times \alpha _2)
\cdot \rho \ {\rm in}\ {\mbox {\it Sh}(n_1,m_1,n_2,m_2)},$$
with $\alpha _{1}\in {\mbox {\it Sh}(n_1,m_1)}$, $\alpha_2\in {\mbox {\it Sh}(n_2,m_2)}$ and $\rho\in {\mbox {\it Sh}(r,n+m-r)}$, with $r=n_1+m_1$.
\medskip

Note that if $(M,\circ )$ is an algebra in $( {\mbox{\it $\mathbb S$--Mod}}, \otimes _{\mathbb S})$, then $\bigoplus M(n)$ is a shuffle algebra with the products 
$x\bullet _{\gamma}y=\circ (x\otimes y\otimes \gamma )$. However, even if $(M,\circ ,\Omega )$ is a bialgebra in $({\mathbb S}-Mod, \otimes _{\mathbb S})$, the space 
$\bigoplus M(n)$ with $\bullet _{\gamma}$ and $\Omega $ is not necessarily a shuffle bialgebra.  But it is possible to obtain two shuffle bialgebras from it, as we describe above.}

\subsubsection {\bf Proposition.} Let $(M,\circ ,\Omega )$ is a bialgebra in $({\mbox{\it $\mathbb S$--Mod}}, \otimes _{\mathbb S})$, then:
\begin{enumerate}\item The shuffle algebra $(M=\bigoplus _{n\geq 0}M(n),\bullet _{\gamma})$ and the coproduct $\Omega _0$ given by:
$$\Omega _0(x):=\sum _{i=0}^nx_{(1)}^{1_n}\otimes x_{(2)}^{1_n},$$
where $1_n$ is considered as a ${\mbox (i,n-i)}$-shuffle for $0\leq i\leq n$, is a shuffle bialgebra.
\item The shuffle algebra $(M=\bigoplus _{n\geq 0}M(n),\bullet _{\gamma})$ and the coproduct $\Omega _{top}$ given by:
$$\Omega _{top}(x):=\sum _{i=0}^nx_{(2)}^{\epsilon _{i,n-i}}\otimes x_{(1)}^ {\epsilon _{i,n-i}},$$
where ${\epsilon _{i,n-i}}$ is considered as a ${\mbox (i,n-i)}$-shuffle for $0\leq i\leq n$, is a shuffle bialgebra.
\end{enumerate}
\medskip

 \P {\rm For elements $x\in M(n)$ and $y\in M(m)$, a shuffle $\gamma \in {\mbox {\it Sh}(n,m)}$ and an integer $0\leq r\leq n+m$, we have that:
$$\gamma = ({\mbox 1_{n_1}\times \epsilon _{n-n_1,m_1}\times 1_{m-m_1}})\cdot (\gamma _{(1)}^r\times \gamma _{(2)}^{n+m-r}),$$
for $n_1=\vert \gamma ^{-1}(\{1,\dots ,r\})\cap \{1,\dots ,n\}\vert $ and $m_1=r-n_1$.

\noindent If $\delta =1_n$ and $\tau =1_m$,then by formula $(*)$ we get that:
$$\displaylines {
({\mbox 1_{n_{1}}\times \epsilon _{m_{1},n-n_{1}}\times 1_{m-m_{1}}})\cdot \gamma=\hfill \cr
({\mbox 1_{n_{1}}\times  \epsilon _{m_{1},n-n_{1}}\times 1_{m-m_{1}}})\cdot ({\mbox 1_{n_1}\times \epsilon _{n-n_1,m_1}\times 1_{m-m_1}})\cdot (\gamma _{(1)}^r\times \gamma _{(2)}^{n+m-r})=\cr
\hfill (\gamma _{(1)}^r\times \gamma _{(2)}^{n+m-r})\cdot 1_{n+m},\cr }$$
which implies the first statement.
\medskip

Suppose that $\delta =\epsilon _{n_1(n-n_1)}$ and $\tau =\epsilon _{m_1(m-m_1)}$. We have that $\epsilon _{r(n+m-r)}=$
$$(1_{n_1}\times \epsilon _{m_1(n-n_1)}\times 1_{m-m_1})\cdot (\epsilon_{n_1(n-n_1)}\times \epsilon_{m_1(m-m_1)})\cdot (1_{n-n_1}\times 
\epsilon _{n_1(m-m_1)}\times 1_{m_1}).$$
Moreover, if $\gamma = ( 1_{n-n_1}\times \epsilon _{n_1,m-m_1}\times 1_{m_1})\cdot (\gamma _{(1)}^{n+m-r}\times \gamma _{(2)}^r)$, then 
$$\epsilon_{r,n+m-r}\cdot (\gamma _{(1)}^{n+m-r}\times \gamma _{(2)}^{r})=(\gamma _{(2)}^r\times \gamma _{(1)}^{n+m-r})\cdot \epsilon _{n+m-r,r}.$$

So, the formula $(*)$ implies that $(x\bullet _{\gamma }y)_{(1)}^{\epsilon _{r,n+m-r}}=
x_{(1)}^{\epsilon _{n_1,n-n_1}}\bullet _{\gamma _{(2)}^r}y_{(1)}^{\epsilon _{m_1,m-m_1}}$ and
$(x\bullet _{\gamma }y)_{(2)}^{\epsilon _{r,n+m-r}}=x_{(2)}^{\epsilon _{n_1,n-n_1}}\bullet _{\gamma _{(1)}^{n+m-r}}y_{(2)}^{\epsilon _{m_1,m-m_1}}.$
We may conclude that:
$$\displaylines {
\Omega _{top}(x\bullet _{\gamma }y)=\sum _{r=0}^{n+m}(x\bullet _{\gamma }y)_{(2)}^{\epsilon _{r,n+m-r}}\otimes (x\bullet _{\gamma }y)_{(1)}^{\epsilon _{r,n+m-r}}=\hfill\cr
\hfill\sum _{r=0}^{n+m}x_{(2)}^{\epsilon _{n_1,n-n_1}}\bullet _{\gamma _{(1)}^{n+m-r}}y_{(2)}^{\epsilon _{m_1,m-m_1}}\otimes 
x_{(1)}^{\epsilon _{n_1,n-n_1}}\bullet _{\gamma _{(2)}^r}y_{(1)}^{\epsilon _{m_1,m-m_1}},\cr}$$
which ends the proof. }\hfill $\diamondsuit$
\bigskip

\noi {\rm {\bf e) The bialgebra of parking functions.} (see \cite{NT2}) Given a function $f\in {\mathcal F}_{n}^{k}$, there exist a 
unique non-decreasing function $f^{\uparrow}\in {\mathcal F}_{n}^{k}$, and 
a unique permutation $\sigma \in {\mbox {\it Sh}(\n)}$ such that
$$f=f^{\uparrow}\cdot \sigma ,$$
where $n_{i}:=\vert f^{-1}(k_{i})\vert $ for $Im(f)=\{ k_{1}<\dots <k_{r}\}$.}

\noi {\rm In \cite{NT2}, J.-C. Novelli and J.-Y. Thibon define a graded map 

\noi ${\mbox {\it Park}}:\bigcup _{n\geq 1}{\mathcal F}_{n}\longrightarrow \bigcup _{n\geq 1}{\mbox {\it PF}}_n$. We give a different description of ${\mbox {\it Park}}$, but it is 
easy to check that it coincides with the one defined in \cite{NT2}.}

\noi  {\rm Let $f^{\uparrow}\in {\mathcal F}_{n}$ be a non-decreasing function, the parking function ${\mbox {\it Park}}(f^{\uparrow})$ is defined as follows:}
$${\mbox {\it Park}}(f^{\uparrow})(j):=\begin{cases}1,&\ {\rm for}\ j=1,\\
{\rm Min}\{ {\mbox {\it Park}}(f^{\uparrow})(j-1))+f^{\uparrow}(j)-f^{\uparrow}(j-1), \ j\}, & \ {\rm for}\ j>1.\end{cases}$$

{\rm For $f=f^{\uparrow}\cdot \sigma$, define:}
$${\mbox {\it Park(f)}}:={\mbox {\it Park}}(f^{\uparrow})\cdot \sigma.$$

\subsubsection  {\bf Remark.}\label{parkin} {\rm \begin{enumerate}\item  Let $f\in PF_n$ be a parking function. It is easy to check that:
\begin{enumerate}
\item $f(i)=f(j)$  if, and only if ${\mbox {\it Park}}(f)(i)={\mbox {\it Park}}(f)(j)$.
\item $f(i)<f(j)$ if, and only if  ${\mbox {\it Park}}(f)(i)<{\mbox {\it Park}}(f)(j)$
\end{enumerate}
 for $1\leq i,j\leq n$.
\item  Let $f,g$ be a pair of parking functions ,
$${\mbox {\it Park}}(f\times g)={\mbox {\it Park}}(f)\times {\mbox {\it Park}}(g).$$
\item If $f\in {\mbox {\it PF}}_n$ is a parking function and $\gamma \in S_n$ is a permutation, then $${\mbox {\it Park}}(f\cdot \gamma )={\mbox {\it Park}}(f)\cdot \gamma .$$
\end{enumerate}}

\noi {\rm Let ${\mbox {\bf {PQSym}}}={\displaystyle \bigoplus _{n\geq 1}}{\mbox {\it PF}}_n$ be the graded space of all parking functions. The coproduct on $\mbox {\bf {PQSym}}$ is defined as follows (see \cite{NT2}). }

{\rm For $f\in {\mbox {\it PF}}_n$ and $0\leq r\leq n$,
$$\Delta _{PQSym}(f):=\sum _{r=0}^{n}{\mbox {\it Park}}(f_1^r)\otimes {\mbox {\it Park}}(f_2^{n-r}),$$
for $f_{1}^{r}:=(f(1),\dots ,f(r))$ and $f_{2}^{n-r}:=(f(r+1),\dots ,f(n)) $.}

\subsubsection {\bf Proposition.} The shuffle algebra $({\mbox {\bf {PQSym}}}, \bullet _{\gamma})$ equipped with the coproduct $\Delta _{PQSym}$ is a shuffle 
bialgebra.
\ms

\P {\rm Let $f\in {\mbox {\it PF}_n}$, $g\in {\mbox {\it PF}_m}$ be parking functions, and let $\gamma $ be a {\mbox {\it (n,m)}}-shuffle. For $0\leq r\leq n$, we want to check that:}
$$\displaylines{ 
{\mbox {\it Park}}((f\bullet _{\gamma }g)_1^r)\otimes {\mbox {\it Park}}((f\bullet _{\gamma }g)_2^{n+m-r})=\hfill\cr
\hfill ({\mbox {\it Park}}(f_1^{n_{1}})\bullet _{\gamma _{(1)}^r}{\mbox {\it Park}}(g_1^{m_{1}}))
\otimes ({\mbox {\it Park}}(f_2^{n-n_{1}})\bullet _{\gamma _{(2)}^{n+m-r}}{\mbox {\it Park}}(g_2^{m-m_1})),\cr }$$
{\rm where} $\gamma =(1_{n_{1}}\times  \epsilon _{m_{1},n-n_{1}}\times 1_{m-m_1})\cdot (\gamma _{(1)}^r\times \gamma _{(2)}^{n+m-r})$ {\rm is the decomposition described in Remark \ref{diagonal}.}

\noi {\rm Computing} $(f\times g)\cdot ({\mbox 1_{n_{1}}\times \epsilon _{m_{1},n-n_{1}}\times 1_{m-m_1}})$, we get that:
$$(f\bullet _{\gamma }g)_1^r=(f_1^{n_1}\times g_1^{m_1}[n-n_1])\cdot \gamma _{(1)}^r$$ {\rm and} 
$$(f\bullet _{\gamma }g)_2^{n-r}=
(f_2^{n-n_1}\times g_2^{m-m_1}[n_1])\cdot  \gamma _{(2)}^{n+m-r},$$
{\rm where} $f_1^{n_1}\times g_1^{m_1}[n-n_1]=(f(1),\dots , f(n_1),g(1)+n,\dots ,g(m_1)+n)$.

\noi {\rm By the Remark \ref{parkin}, we get that 
$$\begin{matrix}
{\mbox {\it Park}}((f\bullet _{\gamma }g)_1^r)&=& {\mbox {\it Park}}(f_{1}^{n_{1}})\bullet_{\gamma _{(1)}^r} {\mbox {\it Park}}(g_{1}^{m_{1}})) ,\\
{\mbox {\it Park}}((f\bullet _{\gamma }g)_2^{n+m-r})& =& {\mbox {\it Park}}(f_{2}^{n-n_{1}})\bullet_{\gamma _{(2)}^{n+m-r}} {\mbox {\it Park}}(g_{2}^{m-m_1})),\end{matrix} $$
which ends the proof. }$\diamondsuit$
\end{exmpls}

\bs

\bs

\bs

\section{Primitive elements of shuffle bialgebras} 
\ms

In this section we use the notations and definitions given in the preliminaries about coalgebras. 

\noi We recall some results proved in \cite{LR3} that we need to study primitive elements in shuffle algebras.
\ms

Following \cite{LR3}, let $(H,\cdot ,\Delta  )$ be a triple such that $( H,\Delta )$ is a connected coassociative 
 coalgebra  and  $( H,\cdot )$ is an associative algebra. Define the $K$-linear map $e\in End _{K}(H)$ as follows:
$$\displaylines {e(x):=\hfill\cr
x-\sum x_{(1)}\cdot x_{(2)}+\dots +(-1)^{r+1}\sum x_{(1)}\cdot x_{(2)}\cdot \dots \cdot x_{(r)}+\dots =
\cr
\hfill \sum _{r\geq 1}(-1)^{r+1}\cdot ^{r}\circ {\Delta}^{r}(x),\cr }$$
where $ {\Delta}^{r}(x)=\sum x_{(1)}\otimes x_{(2)}\otimes \dots \otimes x_{(r)}$.

\begin{prop} \label{idempotent} (see Proposition 2.5 of \cite{LR3}) Any connected infinitesimal bialgebra $(H,\cdot ,\Delta )$ verifies that:
\begin{enumerate}\item  the image of $e$ is ${\rm Prim}(H)$, 
\item the restriction $e\mid _{{\rm Prim}(H)}=Id_{{\rm Prim} (H)}$, and 
\item $e(x\cdot y)=0$ for all $x,y\in Ker(\epsilon )$.  
 \item any element $x$ of $Ker(\epsilon )$ verifies that
$$x=e(x)+\sum e(x_{(1)})\cdot e(x_{(2)})+\dots +\sum e(x_{(1)})\cdot \dots \cdot e(x_{(n)})+\dots ,$$
where ${\Delta}^{n}(x)=\sum x_{(1)}\otimes \dots \otimes x_{(n)}$.
\end{enumerate}
\end{prop}
\ms
 
\begin{thm} \label{cofree} (see Theorem 2.6 of \cite{LR3} ) Any connected infinitesimal bialgebra $H$ is isomorphic to 
$({\overline T}({\rm Prim}(H)):=(\bigoplus _{n\geq 1}{\rm Prim}(H)^{\otimes n},\nu ,\Delta )$, where $\nu$ is the concatenation product and 
$\Delta $ is the deconcatenation coproduct.
\end{thm}
\ms

This section is devoted to compute the subspace of primitive elements of a shuffle bialgebra, and describe it as an algebra for certain type of algebraic 
structure.
\bs

Let $(A,\bullet _{\gamma})$ be a shuffle algebra over $K$. For any pair of positive integers, the permutations $1_{n+m}$ and 
$\epsilon _{nm}:=(n+1,\dots ,m+n,1,\dots ,n)$ belong to $\Sh$. Given elements $x\in A_{n}$ and $y\in A_{m}$, we shall keep the notations $x\bullet _{0}y$ 
for the element $x\bullet _{1_{n+m}}y$, and $y\bullet _{top}x$ for the element $x\bullet _{\epsilon _{nm}}y$, 
in order to simplify notation. Recall that both operations are associative.
\ms

\noi The binary operation $\{- ,- \}:A\otimes A\longrightarrow A$ is given by the formula:
$$\{x,y\}:=x\bullet _{top}y-x\bullet _{0}y,\ {\rm for}\ x,y\in A.$$
\ms 

We want to prove that the subspace of primitive elements of a shuffle bialgebra is closed under the operations $\{- ,- \}$ and $\bullet _{\gamma}$, for $\gamma\in 
\Sh\setminus \{ 1_{n+m},\epsilon _{nm}\}$, with $n,m\geq 1$.
\ms 

\begin{prop}\label{elsrim} Let $(A={\displaystyle \bigoplus _{n\geq 1}}A_{n},\bullet _{\gamma }, \Delta )$ be a shuffle bialgebra. If the homogeneous elements $x\in A_{n}$ and $y\in A_{m}$ belong to ${\rm Prim}(A)$, then $x\bullet _{\gamma } y$ and $\{ x,y\}$ belong to ${\rm Prim}(A)$, for any 
$\gamma \in \Sh\setminus \{ 1_{n+m},\epsilon _{nm}\}$.
\end{prop}
\ms

\P  Note that: 
$$\Delta (x\bullet _{ top}y)=x\otimes y=\Delta (x\bullet _{0}y),\ {\rm for}\ x,y\in {\rm Prim}(A),$$
which implies that $\Delta (\{ x,y\})=0$, whenever $x,y\in {\rm Prim}(A)$.
\ms  

\noi For any permutation $\gamma \in \Sh \setminus \{ 1_{n+m},\epsilon _{nm}\}$, the coproduct  verifies that:
$$\Delta (x\bullet _{\gamma } y)=\sum _r\bigl( \sum (x_{(1)}\bullet _{\gamma _{(1)}^r}y_{(1)})\otimes (x_{(2)}\bullet _{\gamma _{(2)}^{n +
m -r}}y_{(2)})\bigr),$$
where $\gamma =(\gamma _{(1)}^r\times \gamma _{(2)}^{n+m -r})\cdot (1_{n_1}\times \epsilon _{(n-n_1)m_1}\times 1_{ m-m_1}$.

Since $\sum x_{(1)}\otimes x_{(2)}=0$ and $\sum y_{(1)}\otimes y_{(2)}=0$, we get that $(x_{(1)}\bullet _{\gamma _{(1)}^r}y_{(1)})\otimes (x_{(2)}\bullet _{\gamma _{(2)}^{\mbox n +m -r}}y_{(2)})$ if, and only if, $r=n_1=n$ and therefore $\gamma =1_{n+m}$, or $r=m_1=m$ and therefore $\gamma =\epsilon _{nm}$.

But $\gamma \notin \{1_{n+m},\epsilon_{nm}\}$, so $\Delta (x\bullet _{\gamma } y)=0$.\hfill $\diamondsuit$
\bs

In order to get a nice expression for the elements of the subspace of primitive elements of a shuffle bialgebra, let us introduce $(q+1)$-ary 
operations $B_q^{\gamma }$.

\begin{defn}{\rm Let $(A,\bullet _{\gamma })$ be a shuffle algebra over $K$. For $n\geq 1$, let $x\in A_n$, $y_{1}\in A_{m_1}$,\dots ,
$y_{q}\in A_{m_q}$ be elements of $A$, and let 
$\gamma  \in \Sh$ be such that $\gamma ^{-1} (1)<\gamma^{-1} (n+m_{1})$ and $\gamma^{-1} (n+m_{1}+\dots 
+m_{q-1}) +1<\gamma ^{-1}(n)$, where $m ={\displaystyle \sum _{i=1}^{q}}m_{i} $. The $(q+1)$-ary operation $B_{q}^{\gamma }: 
A^{\otimes {q+1}}\longrightarrow A$ is defined by the following formula:}
$$B_{q}^{\gamma }(x;y_{1},\dots ,y_{q}):=x\bullet _{\gamma }(y_{1}\bullet _{0}y_{2}\bullet _{0}\dots \bullet _{0}y_{q}).$$
\end{defn}
\ms

Note that, in the case $n=1$, the conditions on $\gamma $ imply that $\gamma \notin \{ 1_{n+m},\epsilon _{nm}\}$.
So, $B_{1}^{\gamma}(x,y)$ is simply $x\bullet _{\gamma }y$.
\ms

\begin{rem}\label{remsh} {\rm Let $x\in A_{n}$, $y\in A_{m}$, $y_{i}\in A_{m_{i}}$ for $1\leq i\leq q$, and $z\in A_{r}$.}
\begin{enumerate}\item {\rm For $q=1$ the following equalities are immediate to check:}
$$B_{1}^{\gamma}(B_{1}^{\tau}(x;y);z)=\begin{cases}B_{2}^{\delta }(x;y,z),&{\rm for}\ \sigma =1_{m+r }\\
B_{1}^{\delta}(x;B_{1}^{\sigma }(y;z)),&{\rm for}\ \sigma \notin \{1_{\mbox {\it m+r }},\epsilon _{mr}\}\\
B_{1}^{\delta }(x;\{z,y\})+B_{2}^{\delta}(x;z,y),&{\rm for}\ \sigma =\epsilon _{\mbox {\it mr}},\end{cases}$$
{\rm where} $(\tau \times 1_{r})\cdot \gamma=\delta \cdot (1_{n}\times \sigma )\cdot \delta $.
\item {\rm Given a permutation $\tau \in \Sh$ such that $\tau^{-1} (1)<\tau^{-1} (n+m_{1} )$ and 
$\tau^{-1} (n+m_{1} +\dots +m_{n-1})+1<\tau^{-1} (n)$, we get that:}
$$\displaylines {
B_{q}^{\tau}(x;y_{1},\dots ,y_{q})=B_{q-1}^{\tau _{1}}(B_{1}^{\sigma _{1}}(x;y_{1});y_{2},\dots ,y_{q})=\hfill\cr
\hfill B_{1}^{\tau _{q-1}}(B_{1}^{\sigma _{q-1}}(\dots B_{1}^{\sigma _{1}}(x;y_{1});\dots );y_{q-1});y_{q}),\cr }$$
{\rm for} $\tau =(\sigma _{1}\times 1_{\mbox {\it m}_{2} +\dots +m_{q} })\cdot \tau _{1}$ {\rm and} $\tau _{j}=(\sigma _{j+1}\times 
1_{m_{j+1} +\dots +m_{q} })\cdot \tau _{j+1}$, {\rm where} $1\leq j\leq q-2$.
\item {\rm Let $\sigma\in {\mbox {\it Sh}(n,m_1+\dots +m_r)}$, $\sigma\neq 1_{n+m_1+\dots +m_r}$, then there exists integers $0\leq j\leq k+1\leq r+1$ such that:}
$$\displaylines {
\sigma ^{-1}(n+m_1+\dots +m_{j-1})+1<\sigma ^{-1}(1)\leq \sigma^{-1} (n+m_1\dots +m_j)\hfill\cr
\sigma^{-1}(n+m_1+\dots +m_{k-1}+1)<\sigma ^{-1}(n)\sigma ^{-1}(n+m_1+\dots +m_k+1).\hfill\cr }$$
{\rm If $j=k+1$, then $\sigma = \epsilon _{n(m_1+\dots +m_{j-1})}\times 1_{m_j+\dots +m_r}$. In this case, }
$$\displaylines { 
x\bullet _{\sigma }(y_1\bullet _0\dots \bullet_0y_r)=\hfill \cr
\sum_{i=1}^k y_1\bullet _0\dots \bullet _0y_{i-1}\bullet _0
(\{ y_i;x\}\bullet _{\sigma _i}(y_{i+1}\bullet _0\dots \bullet _0 y_k))\bullet _0y_{k+1}\bullet _0\dots \bullet _0 y_r,\cr }$$
{\rm where $\sigma _i=1_{m_i}\times \epsilon _{n(m_{i+1}+\dots +m_k)}$, for $1\leq i\leq k$.}
\medskip

\noi {\rm If $j\leq k$, then} 
$$\sigma =(\epsilon_{n(m_1+\dots +m_{j-1})}\times 1_{m_j+\dots +m_r})\cdot (1_{m_1+\dots +m_{j-1}}\times {\tilde {\sigma}} \times 
1_{m_{k+1}+\dots +m_r}),$$
{\rm where ${\tilde {\sigma}}\in {\mbox {\it Sh}(n;m_j+\dots +m_k)}$. The permutation}
 $$\sigma _i:=(1_{m_i}\times \epsilon _{n(m_{i+1}+\dots +m_{j-1})}\times 1_{m_j+\dots +m_k})\cdot (1_{m_i+\dots +m_j}\times {\tilde {\sigma}},$$
{\rm belongs to ${\mbox {\it Sh}(n,m_i+\dots +m_k)}$ for $1\leq i\leq j$, and we have that:}
 $$\displaylines { 
x\bullet _{\sigma }(y_1\bullet _0\dots \bullet_0y_r)=\hfill \cr
\sum _{i=1}^{j-1}y_1\bullet _0\dots \bullet _0y_{i-1}\bullet _0B_{k-i}^{\sigma _i}(\{y_i;x\};y_{i+1},\dots ,y_k)\bullet _0y_{k+1}\bullet _0\dots \bullet_0y_r+\cr
\hfill y_1\bullet _0\dots \bullet _0y_{j-1}\bullet _0B_{k-j+1}^{\sigma _j}(x;y_j,\dots ,y_k)\bullet _0y_{k+1}\bullet _0\dots \bullet_0y_r.\cr }$$
 \end{enumerate}
\end{rem}
\ms

Point (2) of Remark \ref{remsh} states that given primitive elements  $x,y_1,\dots ,y_q$ in a shuffle bialgebra $A$, the element $B_q^{\gamma }(x;y_1,
\dots ,y_q)$ is primitive too.
\bs

The following Proposition describes the relationships between the operations $B_{q}^{\gamma}$ and $\{ -,-\}$.

\begin{prop}\label{lemsh} Let $x\in A_{n}$,$y\in A_{m}$, $z_{i}\in A_{r_i}$ for $1\leq i\leq q$ and $w\in A_s$ be elements of a shuffle algebra $A$. 
The operations $B_{q}^{\gamma }$ and $\{-,-\}$ defined above verify the following equalities:
\begin{multline} \{x,\{y,w\}\}=\{\{x,y\},w\}+B_{1}^{1_{n}\times \epsilon_{sm}}(\{x,w\};y)-B_{1}^
{\epsilon_{mn}\times 1_{s}}(y;\{ x,w\}).\hfill \\
{\rm (2)}\ \{ x;B_{1}^{\gamma }(y;w)\} = B_{1}^{\overline {\gamma}}(y;\{x,w\})+B_{1}^{\tilde {\gamma}}(\{ x,y\};w),\hfill\\
{\rm where}\ {\overline {\gamma}}:=(\epsilon_{mn}\times 1_s)\cdot (1_{n}\times \gamma )\
{\rm and}\ {\tilde {\gamma}}=1_{n}\times \gamma .\hfill\\
{\rm (3)}\ \{ B_{1}^{\gamma }(x;y),w\}=B_{1}^{\underline {\gamma}}(\{x;w\};y)+B_{1}^{\gamma \times 1_{s}}(x;\{y,w\}),\hfill\\
{\rm where}\ {\underline \gamma}:=(1_{m}\times \epsilon_{sn})\cdot (\gamma\times 1_{s}).\hfill\\
{\rm (4)}\ {\rm Let}\ \gamma \in {\mbox {\it Sh}(n +m, r )}\setminus \{1_{n +m +r}, \epsilon _{(n +m )r}\},\ {\rm and}\ \tau \in {\mbox {\it Sh}(n,m)}\setminus \{ 1_{n +m},\epsilon _{nm}\},\hfill\\
{\rm be\ permutations\ such\ that}\ (\tau \times 1_{r})\cdot \gamma =(1_{n}\times \sigma)\cdot \delta ,\hfill\\
{\rm with}\ \delta \in {\mbox {\it Sh}(n,m +r)}\ {\rm and}\ \sigma\in {\mbox {\it Sh}(m,r)}, \ {\rm where}\ r:=\sum_{i=1}^q r_i. \hfill\\
{\rm a)\ If}\ \sigma =1_{m+r},\ {\rm then}\hfill\\
B_{q}^{\gamma }(B_{1}^{\tau}(x;y);z_{1},\dots ,z_{q})=B_{q+1}^{\delta }(x;y,z_{1},\dots ,z_{q}).\hfill\\
{\rm b)\ If}\ \sigma \neq 1_{m+r},\ {\rm then\ there\ exist\ integers}\ 0\leq j\leq q\ {\rm and}\ 1\leq k\leq q,\ {\rm defined\ in\ the\ same\ way}\\
{\rm that\ in\ point \ (4)\ of\ Remark\ \ref{remsh},\ and\ we\ get\ that}:
{\rm i)\ For}\ j=k+1,\hfill\\
B_{q}^{\gamma }(B_{1}^{\tau}(x;y);z_{1},\dots ,z_{q})=\hfill\\
\sum _{i=1}^{k}B_{q-k+i}^{\delta}(x;z_{1},\dots ,z_{i-1},B_{k-i}^{\sigma _{i}}
(\{z_{i},y\};z_{i+1},\dots ,z_{k}),z_{k+1},\dots ,z_{q})+\hfill\\
\hfill  B_{q+1}^{\delta}(x;z_{1},\dots ,z_{k},y,z_{k+1},\dots ,z_{q}),l\\
{\rm where}\ \sigma _{i}:=1_{r_{i}}\times \epsilon _{(r_{i+1}+\dots +r_{k})m}{\rm ,\ for\ }1\leq i\leq k.\hfill\\
{\rm  ii)\ For}\ j\leq  k,\hfill\\
B_{q}^{\gamma }(B_{1}^{\tau}(x;y);z_{1},\dots ,z_{q})=\hfill\\
\sum _{i=1}^{j-1}B_{q-k+i}^{\delta }(x;z_{1},\dots ,z_{i-1},B_{k-i}^{\sigma _{i}}(\{z_{i},y\},z_{i+1},\dots ,z_k),z_{k+1},\dots ,z_{m})+\\
\hfill B_{q-k+j}^{\delta}(x;z_{1},\dots ,z_{j-1},B_{k-j}^{\sigma _{j}}(y;z_{j},\dots ,z_{k}),z_{k+1},\dots ,z_{m}),\\
{\rm where}\ \sigma_{i}\ {\rm is\ the\ permutation\ defined\ in\ Point\ (4)\ of\ Remark\ \ref{remsh}\ for}\ 1\leq i\leq j.\hfill \\
\end{multline}
\end{prop}
\ms

\P  First, note that if $\gamma \notin \{1_{n +m +r},\epsilon _{(n +m)r}\}$ and 
$\tau \notin \{1_{n+m},\epsilon _{n m}\}$, then $\delta \notin \{1_{n +m +r},
\epsilon _{n (m+r )}\}$.

Points {(1)}, {(2)} and {(3)} are easily verified, while point {(4)} is a straightforward consequence of Remark \ref{remsh}.
\hfill $\diamondsuit $
\bs

\begin{defn}{\rm A {\it ${\mathcal Prim}_{sh}$ algebra} is a graded vector space $V$ equipped with operations $B^{\gamma}:
V_{n}\otimes V_{m}\longrightarrow V$, for $\gamma \in Sh(n,m)\setminus\{1_{n+m},\epsilon _{nm}\}$, and a binary operation $\{ -,-\}$ 
 which verify the following relations: }\begin{enumerate}
  \item $\{ x,\{ y,w\}\}=\{\{x,y\},w\} -B^{1_{n}\times \epsilon _{s m}}(\{x,w\};y)+B^{\epsilon_{mn}\times 1_{s}}(y;\{x,w\});$
 \item $\{ x;B^{\gamma }(y;w)\} = B^{\overline {\gamma}}(y;\{x,w\})+B^{\tilde {\gamma}}(\{ x,y\};w),$
 
\noindent {\rm where} ${\overline {\gamma}}:=(\epsilon_{mn}\times 1_s)\cdot (1_{n}\times \gamma )$ {\rm and} ${\tilde {\gamma}}=1_{n}\times \gamma $;
 \item $\{ B^{\gamma }(x;y),w\}=B^{\underline {\gamma}}(\{x;w\};y)+B^{\gamma \times 1_{s}}(x;\{y,w\})$,

\noindent {\rm where} ${\underline \gamma}:=(1_{m}\times \epsilon_{sn})\cdot (\gamma\times 1_{s})$;
\item $B^{\gamma }(B^{\delta}(x;y);w)=B^{\tau}(x;B^{\sigma}(y;w)),$
 {\rm whenever} $(\delta \times 1_{s})\cdot \gamma=(1_{n}\times \sigma )\cdot \tau$, with $\sigma \neq 1_{m +s}$;
\end{enumerate}
{\rm for} $x\in A_{n}$, $y\in A_{m}$ {\rm and} $w\in A_{s}$.
\end{defn}

\noi Proposition \ref{lemsh} states that there exists a functor from the category ${\mbox  {\it Sh-alg}}$ of shuffle algebras to the category of ${\mathcal Prim}_{sh}$ algebras. Given a 
shuffle bialgebra $A$, the subspace of  primitive elements ${\rm Prim}(A)$ is a ${\mathcal Prim}_{sh}$ subalgebra of $A$.
\bs

Let $(A, \bullet _{\gamma}, \Delta) $ be a shuffle bialgebra, recall that  $(A_{+},\bullet _{0}, \Delta _{+})$ is a unital 
infinitesimal bialgebra. 
\ms

\noi Given a positively graded set $X$, we use Proposition \ref{idempotent} and Theorem \ref{cofree} to describe the free shuffle algebra
 $K[{\mathcal P}_{\infty},X]$ in terms of its primitive elements. We look at the $K$-linear map $e: K[{\mathcal P}_{\infty},X]_{+}\longrightarrow 
{\rm Prim}(K[{\mathcal P}_{\infty},X])$.

\noi For any $x\in X_{n}$, the map $e(\xi _{n},x)=$
$$ (\xi _{n};x)-\sum (\xi _{n_{1},n_{2}};x_{(1)},x_{(2)})+\dots 
+(-1)^{r+1}\sum (\xi _{n_{1},\dots ,n_{r}};x_{(1)},\dots ,x_{(r)})+\dots ,$$
 gives a bijection between $X$ and the subspace $e(X)$ of ${\rm Prim}(K[{\mathcal P}_{\infty},X])$. 
 In order to simplify notation, we denote by $\overline x$ the image under $e$ of the element $(\xi _{n};x)$.
 \ms

\noi Let  ${\mathcal Prim}_{sh}(X)$ be the subspace of $K[{\mathcal P}_{\infty},X]$ spanned by the set $e(X)$ with the operations $B^{\gamma}$ and $\{-,-\}$
, and let ${\mathcal Prim}_{sh}(X)^{\bullet _{0}n}$ be the space spanned by all the elements of the form $z_{1}\bullet _{0}z_{2}\bullet _{0}\dots 
\bullet _{0}z_{n}$, with each $z_{j}\in {\mathcal Prim}_{sh}(X)$, for $1\leq j\leq n$.

\noi Note that Proposition \ref{lemsh} implies that any homogeneous element in ${\mathcal Prim}_{sh}(X)$ is a sum of elements of type $B_{q}^{\gamma}(x;y_{1},
\dots ,y_{q})$, with $x=\{\dots \{x_{1},x_{2}\},x_{3}\},\dots \},x_{n}\} $, for $x_{1},\dots ,x_{n}\in X$, $y_{1},\dots ,y_{q}\in 
{\mathcal Prim}_{sh}(X)$, $q\geq 0$, $n\geq 1$ and $\mid x\mid + \sum \mid y_{j}\mid =n$.
\ms

\begin{prop}\label{todosh} Let $X$ be a positively graded set, equipped with a coassociative graded coproduct $\Theta $ on $K[X]$. Any element $z$ in 
$K[{\mathcal P}_{\infty},X]$ may be written as a sum  $z=\sum _{k}z_{1}^{k}\bullet _{0}z_{2}^{k}\bullet _{0}\dots \bullet _{0}z_{r_{k}}^{k}$, 
with $z_{i}^{k}\in {\mathcal Prim}_{sh}(X)$.
\end{prop}
\ms

\P Clearly, if $x\in X_{1}$ then ${\overline x}=e(\xi _{1},x)\in {\mathcal Prim}_{sh}(X)$. If $x\in X_{n}$, for $n\geq 2$, then
$$(\xi _{n},x)= {\overline x}+\sum (\xi _{n_{1}};x_{(1)})\bullet _{0}{\overline {x_{(2)}}},$$
with ${\overline x}$ and $\overline {x_{(2)}}$ in $e(X)\subseteq {\mathcal Prim}_{sh}(X)$. By a recursive argument
$$(\xi _{n_{1}};x_{(1)})=\sum _{k}\overline {y}_{1}^{k}\bullet _{0}\overline {y}_{2}^{k}\bullet _{0}\dots \bullet _{0}
\overline {y}_{r_{k}}^{k},$$ with $y_{l}^{k}\in X_{m_{l}^{k}}$, for $1\leq l\leq r_{k}$. 

\noi So, any element of the form $(\xi _{n};x)$ belongs to ${\displaystyle \bigoplus _{n\geq 1}}{\mathcal Prim}_{sh}(X)^{\bullet _{0}n}$.
\bs

Since $(K[{\mathcal P}_{\infty},X],\bullet _{\gamma})$ is the free shuffle algebra spanned by all the elements $(\xi _{n};x)$, with $x\in X_{n}$, one has that 
any homogeneous element $y\in K[{\mathcal P}_{\infty},X]_{n}$ may be written in a unique way as $y=\sum _{l}(\xi _{n_{l}};x_{l})\bullet _{\gamma_{l}}y_{l}\rq $, with
$x_{l}\in X_{n_{l}}$, $y_{l}\rq \in K[{\mathcal P}_{\infty},X]$ such that $\vert y_{l}\rq \vert <n$ and $\gamma_{l}\in {\mbox {\it Sh}(n_l,n-n_{l})}$.

\noi We have proved yet that
 $$(\xi _{n};x)=\sum _{l}\overline {x_{1}^{l}}\bullet _{0}\dots \bullet _{0}\overline {x_{r-l}^{l}},$$ 
so, we only need to check that an element $$y={\overline x}\bullet _{\gamma}(z_{1}\bullet _{0}z_{2}\bullet _{0}\dots \bullet _{0}z_{q}),$$ with 
$x\in X_{n}$ and $z_{i}\in {\mathcal Prim}_{sh}(X)$, belongs to $\displaystyle {\oplus _{n\geq 1}}{\mathcal Prim}_{sh}(X)^{\bullet _{0}n}$. 

\noi If $\gamma =1_{\vert y \vert}$, then $y=x\bullet _{0}z_{1}\bullet _{0}\dots \bullet _{0}z_{q}$, with $x$ and $z_{j}\in {\mathcal Prim}_{sh}(X)$.
\ms

\noi If $\gamma \neq 1_{\vert y\vert}$, then there exist $1\leq j\leq q+1$ and $1\leq k\leq q$ such that
$$\displaylines {
\gamma (n+r_{1}+\dots +r_{j-1})\leq \gamma (1)-1<\gamma (n +r_{1} +\dots +r_{j})\hfill\cr
\gamma (n +r_{1}+\dots +r_{k-1}+1)<\gamma (n )<\gamma (n +r_{1} +\dots +r_{k} +1), \hfill \cr}$$
where $\vert z_{i}\vert =r_{i}$.

\noi If $j=k+1$, then $\gamma =\epsilon _{n(r_1+\dots +r_k)}\times 1_{r_{k+1}+\dots +r_{q}}$, and we have that:
$$\displaylines {
y=x\bullet _{\gamma }(z_{1}\bullet _{0}\dots \bullet _{0}z_{q})=B_{k-1}^{\gamma _{1}}(\{ z_{1},x\};z_{2},\dots ,z_{k})\bullet _{0}z_{k+1}\bullet _{0}
\dots \bullet _{0}z_{q}+\hfill \cr
z_{1}\bullet _{0}B_{k-2}^{\gamma _{2}}(\{ z_{2},x\};z_{3},\dots ,z_{k})\bullet _{0}z_{k+1}\bullet _{0}\dots \bullet _{0}z_{q}+\dots +\cr
z_{1}\bullet _{0}\dots \bullet _{0}z_{k-3}\bullet _{0}B_{1}^{\gamma _{q-1}}(\{z_{k-1},x\};z_{k})\bullet _{0}z_{k+1}\bullet _{0}\dots 
\bullet _{0}z_{q}+\cr
z_{1}\bullet _{0}\dots \bullet _{0}z_{k-1}\bullet _{0}\{z_{k},x\}\bullet _{0}z_{k+1}\bullet _{0}\dots 
\bullet _{0}z_{q}+\cr
\hfill z_{1}\bullet _{0}\dots \bullet _{0}z_{k}\bullet _{0}x\bullet _{0}z_{k+1}\bullet _{0}\dots 
\bullet _{0}z_{q},\cr }$$
where $\gamma _{i}=1_{r_{i}}\times \epsilon _{n(r_{i+1}+\dots +r_{k})}$.
\ms

\noi If $j\leq k$, then $$\displaylines {
y=x\bullet _{\gamma }(z_{1}\bullet _{0}\dots \bullet _{0}z_{q})=B_{k-1}^{\gamma _{1}}(\{ z_{1},x\};z_{2},\dots ,z_{k})\bullet _{0}z_{k+1}\bullet _{0}
\dots \bullet _{0}z_{q}+\hfill \cr
z_{1}\bullet _{0}B_{p-2}^{\gamma _{2}}(\{ z_{2},x\};z_{3},\dots ,z_{k})\bullet _{0}z_{k+1}\bullet _{0}\dots \bullet _{0}z_{q}+\dots +\cr
z_{1}\bullet _{0}\dots \bullet _{0}z_{j-2}\bullet _{0}B_{1}^{\gamma _{j-1}}(\{z_{j-1},x\};z_{j},\dots ,z_{k})\bullet _{0}z_{k+1}
\bullet _{0}\dots \bullet _{0}z_{q}+\cr
\hfill z_{1}\bullet _{0}\dots \bullet _{0}z_{j-1}\bullet _{0}B_{k-j+1}^{\gamma _{j}}(x;z_{j},\dots ,z_{k})\bullet _{0}z_{k+1}\bullet _{0}\dots 
\bullet _{0}z_{q},\cr }$$
where the permutations $\gamma _{i}$ are defined in Proposition \ref{lemsh}. \hfill $\diamondsuit$
\bs

Let $X$ be a set, we denote by $$\{X\}_{0}:= \{ \{\{\{\{x_{1},x_{2}\},x_{3}\},\dots \},x_{n}\}\ : x_{i}\in X,\ 1\leq i\leq n\ {\rm and}\ n\geq 1\}$$ which 
is a subset of the free  ${\mathcal Prim}_{sh}$ algebra spanned by $X$. Let $\{X\}={\displaystyle \bigcup _{n\geq 0}}\{X\}_{n}$ be the set  
defined recursively as follows:
$$\displaylines {\{ X\}_{1}:=\{ X\}_{0}\bigcup \{ B_{1}^{\gamma}(x;y),/\ x,y\in \{ X\}_{0}\},\cr
\{X\}_{2}:=\{ X\}_{1}\bigcup \{ B_{1}^{\gamma}(x;y),/\ x\in \{X\}_{0}\ {\rm and}\ y\in \{ X\}_{1}\} \bigcup \hfill\cr
\hfill \{ B_{2}^{\gamma}(x;y_{1},y_{2}),/\ x\in \{X\}_{0}\ {\rm and}\ y_{j}\in \{ X\}_{1},\ {\rm for}\ j=1,2\},\cr
\{X\}_{n}:= \{X\}_{n-1}\bigcup \hfill \cr
 (\bigcup _{m\geq 1}\{ B_{m}^{\gamma}(x;y_{1},\dots ,y_{m}),/\ x\in \{X\}_{0}\ {\rm and}\ 
y_{j}\in \{X\} _{n-1},\ {\rm for}\ 1\leq j\leq m\}).\cr}$$

\noi Note that $\{X\}_{n}\subseteq \{X\}_{n+1}$, and that Proposition \ref{lemsh} and the Remark \ref{remsh} imply that the free ${\mathcal Prim}_{sh}$ 
algebra over $X$ is the vector space spanned by $\{X\}$.
 \ms

\begin{prop} \label{pshuff}Let $X$ be a positively graded set. The subspace ${\mathcal Prim}_{sh}(X)$ is the subspace of primitive elements of $K[{\mathcal P}_{\infty},X]$. 
Moreover, it is the free ${\mathcal Prim}_{sh}$ algebra spanned by $X$.
\end{prop}
\ms

\P Note first that it suffices to prove the result for the case where $X={\displaystyle \bigcup _{n\geq 1}}X_{n}$ with $X_{n}$ finite, for all $n\geq 1$. For the general case,
 $X$ is a limit of graded sets which verify this condition, and the result follows.
 \ms

 \noi  Proposition \ref{elsrim} states that ${\mathcal Prim}_{sh}(X)\subseteq {\rm Prim}(K[{\mathcal P}_{\infty},X])$, while Proposition \ref{todosh} implies that 
$K[{\mathcal P}_{\infty},X]={\overline T}({\mathcal Prim}_{sh}(X))$ as a vector space. From Theorem \ref{cofree} one has that 
$K[{\mathcal P}_{\infty},X]={\overline T}({\rm Prim}(K[{\mathcal P}_{\infty},X]))$, so ${\mathcal Prim}_{sh}(X)={\rm Prim}
(K[{\mathcal P}_{\infty},X])$.
\bs

\noi For the second assertion note that, since $K[{\mathcal P}_{\infty},X]$ is the free associative algebra over the 
set  ${\mbox {\it Irr}_{{\mathcal P},X}}={|displaystyle \bigcup _{n\geq 1}}{\mbox {\it Irr}_{{\mathcal P}_{n,X}}}$ of irreducible elements of ${\mathcal P}_{\infty}$ coloured with elements of $X$, the previous assertion states that
${\mathcal Prim}_{sh}(X)$ is linearly spanned by the set ${\mbox {\it Irr}_{{\mathcal P},X}}$. 
From Proposition \ref{lemsh} we know that ${\mathcal Prim}_{sh}(X)$ is a ${\mathcal Prim}_{sh}$ algebra which contains $e(X)$. To see that it is free, it 
suffices to define a bijective map $\beta :{\mbox {\it Irr}_{{\mathcal P},X}}\longrightarrow \{X\}$, where $\{ X\}$ is 
the set defined above. On $X\subset {\mbox {\it Ir}_{{\mathcal P},X}}$, 

\noi $\beta$ coincides with the identity map. 
Let $y=(\xi _{n},x)\bullet _{\gamma }y_{1}\in {\mbox {\it Irr}_{{\mathcal P},X}}$, with $x\in X_{n}$ and $n\geq 1$. We define $\beta (y)$ as follows:

\noi If $y_{1}\in {\mbox {\it Irr}_{{\mathcal P},X}}$ and $\gamma \neq \epsilon _{nm_{1}}$, then $\beta (y):=B_{1}^{\gamma }(x;\beta (y_{1}))$.

\noi If $y_{1}\in {\mbox {\it Irr}_{{\mathcal P},X}}$ and $\gamma = \epsilon _{nm_{1}}$, then 
$$\beta (y):=\begin{cases}\{ \beta (y_{1}),x\},&\ {\rm for}\ \beta (y_{1})\in \{X\}_{0} ,\\
B_{q}^{\overline \tau}(\{w,x\};z_{1},\dots ,z_{q}),&\ {\rm for}\ \beta (y_{1})=B_{q}^{\tau}(w;z_{1},\dots ,z_{q}),\end{cases}$$
where $\vert y_{1}\vert =m_{1}$, $w\in \{X\}_{0}$, $\vert w\vert =s$, $z_{j}\in \{ X\}$ for $1\leq j\leq q$, $\sum _j\vert z_{j}\vert =r$ and 
${\overline \tau}:=(\tau\times 1_n)\cdot (1_s\times \epsilon _{nr})$.
\ms

\noi Suppose that $y_{1}=t_{1}\bullet _{0}\dots \bullet _{0}t_{p}$, with $t_{i}\in {\mbox {\it Irr}_{{\mathcal P},X}}$ and $p>1$. The fact that $y$ is irreducible, implies that
  $\gamma (n +h_1 +\dots +h_{p-1}+1)<\gamma (n)$, for $\vert t_i\vert =h_i$. 
  
\noi If $\gamma (1)<\gamma (n+h_1)$, then $\beta (y):=B_{p}^{\gamma } (x;\beta (t_1),\dots ,\beta (t_{p}))$.
 \ms
  
\noi If $\gamma (n+h_{1} )<\gamma (1)-1$, then 
  $$\beta (y):=\begin{cases}B_{p-1}^{\overline \gamma}(\{\beta (t_{1}) ,x\};\beta (t_{2}),\dots ,\beta (t_{p})),&\ {\rm for}\ \beta (t_{1})\in \{X\}_{0}\\
 B_{q+p-1}^{\overline \tau}(\{w,x\};z_{1},\dots ,z_{q},\beta  (t_{2}),\dots ,\beta (t_{p})),&\ {\rm for}\ \beta (t_{1})=B_{q}^{\tau }(w;z_{1},\dots ,
 z_{q}),\end{cases}$$
 where ${\overline {\gamma}}=\gamma \cdot (\epsilon _{h_{1}n}\times 1_{h_{2}+\dots +h_{p}})$ and ${\overline \tau}:=
 \gamma \cdot (1_{n}\times \tau \times 1_{h_{2}+\dots +h_{r}})\cdot (\epsilon _{sr}\times 1_{h-s}),$ where $h:={\displaystyle \sum _{i=1}^{p}}h_{i}$.
\ms
  
\noi It is not difficult to check that $\beta $ is bijective, which implies that  ${\mathcal Prim}_{sh}(X)$ is isomorphic to the free ${\mathcal Prim}_{sh}$ algebra spanned 
by $X$ \hfill $\diamondsuit $
\bs

The following result is a straightforward consequence of Theorem \ref{cofree} and the previous results.
\ms

\begin{prop}\label{shfr} Let $X$ be a positively graded set, such that $K[X]$ is equipped with a coassociative graded coproduct $\Theta$. 
The unital infinitesimal bialgebra  $K[{\mathcal P}_{\infty},X]_{+}$ is isomorphic to $T^{fc}({\mathcal Prim}_{sh}(X))$, where 
${\mathcal Prim}_{sh}(X)$ is the free ${\mathcal Prim}_{sh}$ algebra spanned by $X$.
\end{prop}
\ms

We want to prove the equivalence between the categories of connected shuffle bialgebras and ${\mathcal Prim}_{sh}$ algebras.  More precisely, 
given a ${\mathcal Prim}_{sh}$ algebra $(V, B^{\gamma}, [-,-])$ and an homogeneous basis $X$ of the underlying vector space $V$, let 
${\mathcal U}_{Sh}(V)$ be the shuffle bialgebra obtained by taking the quotient of the free shuffle algebra $K[{\mathcal P},X]$ by the ideal 
(as a shuffle algebra) spanned by the set:
$$\{ B^{\gamma}(x;y) -{\overline B}^{\gamma}(x;y),\ \{ x,z\} -
[ x, z]\},$$ 
with $x\in X_{n}$, $y\in X_{m}$, $z\in X_{r}$ and $\gamma \in Sh(n ,m)$ such that $\gamma (1)<\gamma (n+m_{1} )$ 
and $\gamma (n+m_{1} +\dots +m_{n-1}+1)<\gamma (n )$, where 
$B^{\gamma }$ and $\{ -,-\}$ denote the operations associated to the shuffle algebra $K[{\mathcal P}_{\infty},X]$. 
\ms

\begin{thm}\label{MMsh} {\bf a)} Let $(H,\circ _{\gamma},\Delta )$ be a connected shuffle bialgebra, then $H$ is isomorphic to ${\mathcal U}_{Sh}({\rm Prim} (H))$, where 
${\rm Prim}(H)$ is the ${\mathcal Prim}_{sh}$ algebra of primitive elements of $H$.

\noi {\bf b)} Let $(V,{\overline B}_{n}^{\gamma },\{ -,-\})$ be a ${\mathcal Prim}_{sh}$ algebra, then $V$ is isomorphic to ${\rm Prim}({\mathcal U}_{Sh}(V))$.
\end{thm}
\ms

\P {\bf a)} Let $X$ be a basis of the vector space Prim$(H)$. Define $\varphi : Sh(X)\longrightarrow H$ as follows:
$$\varphi (x_{1}\bullet _{\gamma _{1}}(x_{2}\bullet _{\gamma _{2}}(\dots (x_{n-1}\bullet _{\gamma _{n-1}x_{n}})))):=
x_{1}\circ _{\gamma _{1}}(x_{2}\circ _{\gamma _{2}}(\dots (x_{n-1}\circ _{\gamma _{n-1}x_{n}}))),$$
where $x_{i}\in X$ for $1\leq i\leq n$.
Note that $\varphi (B^{\gamma}(x;y))=\varphi (x_{\gamma}y)=x\circ _{\gamma}y={\overline B}^{\gamma}(x;y)$, and $\varphi (\{x,y\})=
\varphi (x\bullet_{top} y - x\bullet _{0}y)=[x;y]$, so $\varphi $ factorizes through ${\mathcal U}_{Sh}({\rm Prim} (H))$. Moreover, it is immediate to check 
that $\varphi$ is a bialgebra homomorphism.
Applying Theorem \ref{cofree}, the inverse morphism is given by $$\varphi ^{-1}(x)=cl(e(x))+\sum cl(x_{(1)}\bullet _{0}x_{2}))+\dots 
+\sum cl(x_{1}\bullet _{0}\dots \bullet _{0}x_{n}),$$
where $cl$ denotes the class of the element in ${\mathcal U}_{Sh}(V)$
\ms

\noi {\bf b)} It is clear that $V\subseteq {\rm Prim}({\mathcal U}_{Sh}(V))$. If $X$ is a basis of $V$, then Proposition \ref{shfr} implies that ${\rm Prim}({\mbox {\it Sh}(X)})=
{\mathcal Prim}_{sh}(X)$. So, the primitive elements of ${\mathcal U}_{Sh}(V)$ are generated by $X$ under the operations $B^{\gamma}$ and $\{ ,\}$, which 
are precisely the elements of $V$ in the quotient.\hfill $\diamondsuit$

\bs

\bs

\bs

\section{Boundary map on the free shuffle bialgebra}
\ms

Let $X$ be a positively graded set, and let $\Theta :K[X]\longrightarrow K[X]\otimes K[X]$ be a coassociative graded coproduct on $K[X]$.

\noi Given an element $x=x_{1}\bullet _{\gamma _{1}}(x_{2}\bullet _{\gamma _{2}}(\dots (x_{r-1}\bullet _{\gamma _{n-1}x_{n}})))$ in the free shuffle 
algebra $Sh(X)$, with $x_{i}\in X$ for $1\leq i\leq r$, the {\it weight} of $x$ is $r$. We denote the weight of an element $x$ by $w(x)$. Note that the elements of 
weight $r$ in ${\mbox {\it Sh}(X)}$ correspond to the elements of the subspace ${\displaystyle \bigoplus _{n\geq 1}}K[{\mathcal P}_{n,X}^{r}]$.

\noi Define the linear map $\partial _{\Theta}:{\mbox {\it Sh}(X)}\longrightarrow {\mbox {\it Sh}(X)}$ as the unique linear homomorphism such that:\begin{enumerate}
\item For $x\in X_{n}$, 
$$\partial _{\Theta}(x):=\sum (\sum _{\sigma \in Sh(n_{1},n-n_{1})}(-1)^{n_{1}}sgn(\sigma)\ x_{(1)}\bullet _{\sigma}x_{(2)}),$$
where $\Theta (x)=\sum x_{(1)}\otimes x_{(2)}$, and $\vert x_{(1)}\vert =n_{1}$.
\item For $x,y\in Sh(X)$, with $w(x)=r$,
$$\partial _{\Theta}(x\bullet _{\gamma}y):=\partial _{\Theta}(x)\bullet _{\gamma}y -(-1)^{r}x\bullet _{\gamma }\partial_{\Theta}(y).$$
\end{enumerate}
Note that $\partial _{\Theta}({\mathcal P}_{n,X}^{r})\subseteq {\mathcal P}_{n,X}^{r-1}$.

\begin{lem}\label{difsh} The homomorphism $\partial _{\Theta}$ is a boundary map, that is $\partial_{\Theta}^{2}=\partial_{\Theta}\circ \partial
_{\Theta}=0$.
\end{lem}

\P Let $x\in X_{n}$, we have that:
$$\displaylines {
\partial_{\Theta}^{2}(x)=\hfill\cr
\sum _{{0<n_{1}<n}\atop {0<n_{2}<n-n_{1}}}\bigl ((\sum _{{\gamma \in Sh(n_{1}+n_{2},n-n_{1}-n_{2})}\atop {\sigma \in Sh(n_{1},n_{2})}}
(-1)^{2n_{1}+n_{2}}sgn(\gamma)sgn(\sigma)(x_{(1)}\bullet _{\sigma}x_{(2)})\bullet _{\gamma}x_{3})+\hfill\cr
\hfill (\sum _{{\delta \in Sh(n_{1},n-n_{1})}\atop {\tau \in Sh(n_{2},n-n_{1}-n_{2})}}(-1)^{2n_{1}+n_{2}+1}sgn(\delta )sgn(\tau)
x_{(1)}\bullet _{\delta}(x_{(2)}\bullet _{\tau}x_{(3)})\bigr),\cr }$$
where $\Theta ^{3}(x)=\sum x_{(1)}\otimes x_{(2)}\otimes x_{(3)}$ with $\mid x_{(i)}\mid =n_{i}$ for $i=1,2$.
So, 
$$\displaylines {
\partial_{\Theta}^{2}(x)=\hfill\cr
\sum _{n_{1}+n_{2}+n_{3}=n}\bigl(\sum _{\alpha \in Sh(n_{1},n_{2},n_{3})}(\sum (-1)^{n_{2}}sgn(\gamma)sgn(\sigma )(x_{(1)}\bullet _{\sigma}x_{(2)})
\bullet _{\gamma}x_{(3)})+\hfill\cr
\hfill(\sum (-1)^{n_{2}+1}sgn(\delta)sgn(\tau)x_{(1)}\bullet _{\delta}(x_{(2)}\bullet _{\tau}x_{(3)}))\bigr)=\cr
\sum _{n_{1}+n_{2}+n_{3}=n}\bigl(\sum _{\alpha \in Sh(n_{1},n_{2},n_{3})}\sum sgn(\alpha)( (-1)^{n_{2}}+(-1)^{n_{2}+1})
(x_{(1)}\bullet _{\sigma}x_{(2)})\bullet _{\gamma}x_{(3)})\bigr )=0,\cr }$$
where $\alpha =(\sigma \times 1_{n_{3}})\cdot \gamma=(1_{n_{1}}\times \tau )\cdot \delta$, for $\alpha \in {\mbox {\it Sh}(n_{1},n_{2},n_{3})}$ and 
$\vert x_{(i)}\vert =n_{i}$.
\ms

Suppose that $x=y\bullet _{\gamma}z$, with $y\in X$ and $z\in K[{\mathcal P}_{n-1,X}^{r}]$. Applying the definition of $\partial _{\Theta}$ 
and a recursive argument, we have that:
$$\displaylines {
\partial_{\Theta}^{2}(x)=
\partial_{\Theta}(\partial_{\Theta}(y)\bullet _{\gamma}z-(-1)^{1}y\bullet _{\gamma}\partial_{\Theta}(z))=\cr
\hfill (-1)^{w(\partial_{\Theta}(y))}\partial_{\Theta}(y)\bullet _{\gamma}\partial_{\Theta}(z)-(-1)^{1}
\partial_{\Theta}(y)\bullet _{\gamma}\partial_{\Theta}(z)=0.\hfill \diamondsuit \cr}$$
\ms

In the case that $X=\{ \xi _{n}\}_{n\geq 1}$, we have that ${\mbox {\it Sh}(X)}=K[{\mathcal P}_{\infty}]$. If we consider the coassociative coproduct $\Theta$ on $K[X]$ given by:
$$\Theta (\xi _{n}):=\sum _{i=1}^{n-1}\xi _{i}\otimes \xi _{n-i},$$
then the boundary map $\partial_{\Theta }$ coincides with the boundary map of the permutahedron (see \cite{Mil}). 
\bs

\bs

\bs

\section{Relations with dendriform and $2$-associative algebras}
\ms

In the previous sections we have defined functors:\begin{enumerate}
\item $H_{\mbox {\it inf-gr}}$ from the category of graded unital infinitesimal bialgebras to the category of shuffle bialgebras, 
\item two functors from the category of shuffle algebras into the category of associative algebras, in such a way that the restrictions of both functors to the 
category of shuffle bialgebras, $H_{0}$ and $H_{L}$, have their images contained in the subcategory of unital infinitesimal bialgebras,
\item from the category of shuffle algebras to the category of associative algebras, in such a way that the image under this functor of a shuffle bialgebra gives an 
associative bialgebra. We denote the restriction of this functor to the category of shuffle bialgebras  by $H_{\mbox {\it sh-as}}$.
\end{enumerate}

 In \cite{Ag}, M. Aguiar constructs functors relating infinitesimal bialgebras (see \cite{JR}) to dendriform algebras (see \cite {Lo}) and brace algebras (see \cite{Kad} 
 and \cite{Get}). We want to include shuffle bialgebras in his framework.
 \ms

\begin{defn} {\rm A {\it dendriform algebra} over $K$ is a vector space $D$ equipped with two bilinear maps $\succ , \prec :D\otimes D\longrightarrow D$
which verify the following relations:\begin{enumerate}
\item $x\succ (y\succ z)=(x\succ y+x\prec y)\succ z$,
\item $x\succ (y\prec z)=(x\succ y)\prec z$,
\item $x\prec (y\succ z + y\prec z)=(x\prec y)\prec z$,
\end{enumerate}
for $x,y,z\in D$.}
\end{defn}

Note that any dendriform algebra $D$ has a natural structure of associative algebra with the product $*$, defined by: $x*y=x\succ y + x\prec y$.
\ms

For nonnegative integers $n,m$, let ${\mbox {\it Sh}^{\succ }(n,m)}$ and ${\mbox {\it Sh}^{\prec }(n,m})$ be the following subsets of $\Sh$:

\noindent a) ${\mbox {\it Sh}^{\succ }(n,m)}:=\{ \sigma \in \Sh \mid \sigma (n+m)=n+m\},$

\noindent b) ${\mbox {\it Sh}^{\prec }(n,m}):=\{ \sigma \in \Sh \mid \sigma (n+m)=n\}.$

It is immediate to check that $\Sh$ is the disjoint union of ${|mbox {\it Sh}^{\succ}(n,m)}$ and ${\mbox {\it Sh}^{\prec }(n,m)}$. 

Let $(A,\bullet _{\gamma})$ be a shuffle algebra, Define on $A$ the operations $\succ $ and $\prec $ as follows:\begin{enumerate}
\item $x\succ y:={\displaystyle \sum _{\gamma \in {\mbox {\it Sh}^{\succ }(n,m)}}}x\bullet _{\gamma }y,$
\item $x\prec y:={\displaystyle \sum _{\gamma \in {\mbox {\it Sh}^{\prec }(n,m)}}}x\bullet _{\gamma }y,$
\end{enumerate} 
for $x\in A_{n}$ and $y\in A_{m}$. 
Note that the associative product $*$ defined in Lemma \ref{shuffprod} is the sum of $\succ $ and $\prec $. 
\ms

\begin{lem} Let $(A,\bullet _{\gamma})$ be a shuffle algebra, then $(A,\succ ,\prec )$ is a dendriform bialgebra. Moreover, if $(A,\bullet _{\gamma},\Delta )$ 
is a shuffle bialgebra, then $\succ $, $\prec $ and $\Delta $ verify the following equalities:
$$\displaylines {
\Delta (x\succ y)= \sum( x_{(1)}*y_{(1)})\otimes (x_{(2)}\succ y_{(2)})+\sum x_{(1)}\otimes (x_{(2)}\succ y)+\cr
\sum y_{(1)}\otimes (x\succ y_{(2)})+\sum (x*y_{(1)})\otimes y_{(2)}+x\otimes y,\cr
\Delta (x\prec y)=\sum (x_{(1)}*y_{(1)})\otimes (x_{(2)}\prec y_{(2)})+\sum y_{(1)}\otimes (x\prec y_{(2)})+\cr
\sum (x*y_{(1)})\otimes y_{(2)}+\sum x_{(1)}\otimes (x_{(2)}\prec y)+y\otimes x,\cr }$$
for all $x,y\in A$.
\end{lem}
\ms

\P To prove that $(A,\succ ,\prec )$ is a dendriform algebra, it suffices to note that the equality:
$$ (1_{n}\times {\mbox {\it Sh}(m,r)})\cdot {\mbox {\it Sh}(n,m+r)}=({\mbox {\it Sh}(n,m)}\times 1_{r})\cdot {\mbox {\it Sh}(n+m,r)},$$
is the sum of the following three equalities:
$$\displaylines {
(1_{n}\times {\mbox {\it Sh}^{\succ }(m,r)})\cdot {\mbox {\it Sh}^{\succ }(n,m+r)}=(\Sh\times 1_{r})\cdot {\mbox {\it Sh}^{\succ }(n+m,r)}\cr
(1_{n}\times {\mbox {\it Sh}^{\prec }(m,r)})\cdot {\mbox {\it Sh}^{\succ }(n,m+r)}=({\mbox {\it Sh}^{\succ}(n,m)}\times 1_{r})\cdot {\mbox {\it Sh}^{\prec }(n+m,r)},\cr
(1_{n}\times {\mbox {\it Sh}(m,r)})\cdot {\mbox {\it Sh}^{\prec}(n,m+r)}= ({\mbox {\it Sh}^{\prec}(n,m)}\times 1_{r})\cdot {\mbox {\it Sh}^{\prec}(n+m,r)}.\cr }$$

On the other hand given $1\leq r\leq n+m$, if $\gamma = (1_{n_{1}}\times \epsilon _{n-n_{1},m_{1}}\times 1_{m-m_1})\cdot (\gamma _{(1)}^{r}\times \gamma _{2}^{n+m-r})\in \Sh$ is the decomposition given in Remark \ref{diagonal} with $\gamma _{(1)}^{r}\in {\mbox {\it Sh}(n_1,m_1)}$,
$\gamma _{(2)}^{n+m-r}\in {\mbox {\it Sh}(n-n_1,m-m_1)}$ and $r=n_1+m_1$, then it is easy to prove that
$\gamma \in Sh^{\succ }(n,m)$ if, and only if $\gamma _{(2)}^{n+m-r}\in Sh^{\succ }(n-r_{1},m+r_{1}-r)$, which implies the formulas for the coproduct.\hfill $\diamondsuit$

The Lemma above states that any shuffle bialgebra has a natural structure of dendriform bialgebra, as defined in \cite{Ro}.
\ms

\begin{defn}{\rm (see \cite{LR3}) Let $X$ be a $K$-vector space equipped with two associative products $*$ and $\cdot $, and a coassociative coproduct $\Delta $, such that:
\begin{enumerate}
\item $(X,*,\Delta )$ is a bialgebra over $K$,
\item $(X,\cdot ,\Delta )$ is an infinitesimal unital bialgebra.
\end{enumerate}
Then $(X,*,\cdot ,\Delta )$ is called a $2${\rm -associative bialgebra}.}
\end{defn}
\ms

The following statement is a consequence of Corollary \ref{hs}.

If $(A,\bullet _{\gamma},\Delta )$ is a shuffle bialgebra, then $A_{+}=K \oplus A$ with the products:\begin{enumerate}
\item $x*y=\begin{cases}
{\displaystyle \sum _{\gamma \in Sh(n,m)}}x\bullet _{\gamma}y,&\ {\rm for}\ x\in A_{n},\ y\in A_{m},\\
x,&\ {\rm for}\ y=1_{K},\\
y,&\ {\rm for }\ x=1_{K},\end{cases} $ and
\item $x\cdot y=\begin{cases}x\bullet _{1_{n+m}}y,&\ {\rm for}\ x\in A_{n},\ y\in A_{m},\\
x,&\ {\rm for}\ y=1_{K},\\
y,&\ {\rm for }\ x=1_{K},\end{cases} $
\end{enumerate} 
is a $2$-associative bialgebra.
\bs

In previous work, we prove that:\begin{enumerate}
\item The subspace of primitive elements of a dendriform bialgebra has natural structure of brace algebra. Moreover, there exists an equivalence between the category of connected 
dendriform bialgebras and the category of brace algebras (see \cite{Ro}).
\item The subspace of primitive elements of a $2$-associative bialgebra has a natural structure of non-differential $B_{\infty }$ algebra (see \cite{LR3}). As in the previous case, 
the category of connected $2$-associative bialgebras is equivalent to the category of non-differential $B_{\infty}$ algebras.
\end{enumerate}

A non-differential $B_{\infty}$ algebra $B$ is a $K$-vector space equipped with linear maps $B_{n,m}:B^{\otimes (n+m)}\longrightarrow B$, verifying certain 
relations (see \cite{Vo}). 

\noi In particular, any brace algebra $(B, M_{1,n})$ is a non-differential $B_{\infty}$ algebra such that $B_{1,m}=M_{1,m}$, for $m\geq 1$, and $B_{n,m}=0$, 
for $n\neq 1$. 

\noi  The previous results show that the associative bialgebra associated to any shuffle bialgebra has structures of both dendriform bialgebra and $2$-associative bialgebra. 
The subspace of primitive elements of a shuffle bialgebra is a brace algebra, as well as a non-differential $B_{\infty}$ algebra. However, these structures 
are not always isomorphic as non-differential $B_{\infty}$ algebras, even if they give the same associative bialgebra structure. For instance, using the formulas given in 
\cite{LR3}, the subspace of primitive elements of the Malvenuto-Reutenauer Hopf algebra $K[S_{\infty}]$ has a structure of non-differential $B_{\infty}$ algebra 
which is not isomorphic to the brace algebra structure given in \cite{Ro}.
\bs 

\bs

\bs

\section{Preshuffle and grafting bialgebras} 
\ms

We introduce the notions of preshuffle algebras, related to leveled trees, and of a particular type of preshuffle algebras called grafting algebras, related to trees. 

\noi Shuffle algebras are related to monoids in the category $({\mbox{\it $\mathbb S$--Mod}}, \otimes _{\mathcal S})$, where we do no ask for a compatibility relation between 
the operations $\bullet _{\gamma}$ and the action of the symmetric group. In a similar way, grafting algebras are related to non-symmetric operads (see {Ma}); 
a grafting algebra structure on $A=\bigoplus A_{n}$ is equivalent to a non-symmetric operad ${\mathbb P}$ with ${\mathbb P}(n)=A_{n-1}$.

\begin{defn} \label{preshgraft} {\rm \begin{enumerate}\item A {\it preshuffle } algebra over $K$ is a graded vector space $A={\displaystyle \bigoplus _{n\geq 0}}A_{n}$ 
equipped with linear maps
$$\bullet _{i}:A\otimes A_{m}\rightarrow A,\ {\rm for}\ 0\leq i\leq m\ {\rm and}\ m\geq 0,$$
verifying:
 $$(x\bullet _i y)\bullet _j z=x\bullet _{i+j}(y\bullet _jz),\ {\rm for}\ 0\leq i\leq \vert y\vert \ {\rm and}\ 0\leq j\leq \vert z\vert.$$ 
 \item A {\it grafting algebra} is a preshuffle algebra $(A,\bullet _{i})$ such that the operations $\bullet _{i}$ 
verify the following additional conditions:
$$x\bullet _{i}(y\bullet _{j}z)= y\bullet _{j+\vert x\vert }(x\bullet _{i}z),\ {\rm for}\ 0\leq i<j,$$
for any elements $x,y,z\in A$.
\end{enumerate}}
\end{defn}

The relations verified by a preshuffle algebra may be pictured as follows:
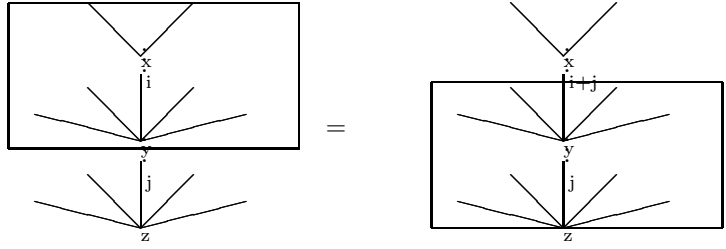
\begin{center}
\begin{picture}(400,100)
\put(80,0){\scriptsize z}
\put(80,33){\scriptsize y}
\put(80,66){\scriptsize x}
\put(82,20){\scriptsize j}
\put(82,58){\scriptsize i}
\put(60,25){\line(1,-1){20}}
\put(60,90){\line(1,-1){20}}
\put(60,90){\line(1,-1){20}}
\put(40,15){\line(4,-1){40}}
\put(80,5){\line(0,1){25}}
\put(80,5){\line(1,1){20}}
\put(80,70){\line(1,1){20}}
\put(80,70){\line(1,1){20}}
\put(80,5){\line(4,1){40}}
\put(80,30){\vdots}
\put(80,64){\vdots}
\put(40,48){\line(4,-1){40}}
\put(60,58){\line(1,-1){20}}
\put(80,38){\line(1,1){20}}
\put(80,38){\line(4,1){40}}
\put(80,38){\line(0,1){25}}
\put(60,25){\line(1,-1){20}}
\put(40,15){\line(4,-1){40}}
\put(80,5){\line(0,1){25}}
\put(80,5){\line(1,1){20}}
\put(80,5){\line(4,1){40}}
\put(40,48){\line(4,-1){40}}
\put(60,58){\line(1,-1){20}}
\put(80,38){\line(1,1){20}}
\put(80,38){\line(4,1){40}}
\put(80,38){\line(0,1){25}}
\put(30,35){\line(1,0){110}}
\put(30,35){\line(0,1){55}}
\put(30,90){\line(1,0){110}}
\put(140,35){\line(0,1){55}}
\put(150,40){=}
\put(240,0){\scriptsize z}
\put(240,33){\scriptsize y}
\put(240,66){\scriptsize x}
\put(242,20){\scriptsize j}
\put(242,58){\scriptsize i+j}
\put(220,25){\line(1,-1){20}}
\put(220,90){\line(1,-1){20}}
\put(220,90){\line(1,-1){20}}
\put(200,15){\line(4,-1){40}}
\put(240,5){\line(0,1){25}}
\put(240,5){\line(1,1){20}}
\put(240,70){\line(1,1){20}}
\put(240,70){\line(1,1){20}}
\put(240,5){\line(4,1){40}}
\put(240,30){\vdots}
\put(240,64){\vdots}
\put(200,48){\line(4,-1){40}}
\put(220,58){\line(1,-1){20}}
\put(240,38){\line(1,1){20}}
\put(240,38){\line(4,1){40}}
\put(240,38){\line(0,1){25}}
\put(220,25){\line(1,-1){20}}
\put(200,15){\line(4,-1){40}}
\put(240,5){\line(0,1){25}}
\put(240,5){\line(1,1){20}}
\put(240,5){\line(4,1){40}}
\put(200,48){\line(4,-1){40}}
\put(220,58){\line(1,-1){20}}
\put(240,38){\line(1,1){20}}
\put(240,38){\line(4,1){40}}
\put(240,38){\line(0,1){25}}
\put(190,5){\line(1,0){110}}
\put(190,5){\line(0,1){55}}
\put(190,60){\line(1,0){110}}
\put(300,5){\line(0,1){55}}
\end{picture}
\end{center}
\ms

For non-negative integers $n, m$ and $0\leq i\leq m$, let 
$\omega _{i}^{n,m}$ be the permutation $$\omega _{i}^{n,m}:=\epsilon_{n,i}\times 1_{m-i}\in \Sh.$$
It is immediate to see that the following equality holds:
$$ (1_{n}\times \omega _{j}^{m,r})\cdot \omega _{i}^{n,m+r}= (\omega _{i-j}^{n,m}\times 1_{r})\cdot \omega _{j}^{n+m,r}.$$
Moreover, note that the element $\omega _{m}^{n,m}=\epsilon_{n,m}$
\ms 

Note that, since the permutation $\omega _{i}^{n,m}$ belongs to $\Sh$, for $0\leq i\leq m$, any 
shuffle algebra is a preshuffle algebra, with the operations 
$$x\bullet _{i}y:=x\bullet _{\omega _{i}^{n,m}}y,\ {\rm for}\ x\in A_{n},\ y\in A_{m}\ {\rm and}\ 0\leq i\leq m.$$ 

\begin{defn}\label{defbia} {\rm {\bf (1)} Let $(A,\bullet _i)$ be a positively graded preshuffle algebra, such that $A$ is equipped with a graded coassociative 
coproduct $\Delta $. We say that $(A,\bullet _i , \Delta  )$ is a {\it preshuffle bialgebra} if it verifies:}
$$\displaylines { (1)\ \Delta (x\bullet _{0} y)=\sum x_{(1)}\otimes (x_{(2)}\bullet _{0}y)+ x\otimes y +\sum (x\bullet _{0}y_{(1)})\otimes y_{(2)}.
\hfill\cr
(2)\ \Delta (x\bullet _i y)=\sum _{\vert y_{(1)}\vert \leq i} y_{(1)}\otimes (x\bullet _{i-\vert y_{(1)}\vert }y_{(2)})+\hfill\cr
\hfill \sum _{\vert y_{(1)}\vert =i} (x_{(1)}\bullet _{i}y_{(1)})\otimes (x_{(2)}\bullet _{0} y_{(2)})+\sum _{\vert y_{(1)}\vert \geq i}(x\bullet _{i}y_{(1)})
\otimes y_{(2)},\cr 
{\rm for }\ 1\leq i\leq \vert y\vert .\hfill \cr
(3)\ \Delta (x\bullet _{\vert y\vert }y)=\sum y_{(1)}\otimes (x\bullet _{\vert y_{(2)}\vert }y_{(2)}) + y\otimes x + \sum (x_{(1)}\bullet _{\vert y\vert }
y)\otimes x_{(2)}.\hfill \cr }$$
{\bf (2)} {\rm A {\it grafting bialgebra} is a preshuffle bialgebra $(A,\bullet _{i},\Delta )$ such that $(A,\bullet _{i})$ is a grafting algebra.}
\end{defn}
\ms 

It is immediate to check that any shuffle bialgebra is a preshuffle bialgebra.
\bs

\begin{exmpls} \label{infunitbialg} \noi{\rm {\bf a) The free preshuffle algebra.} Let $X={\displaystyle \bigcup _{n\geq 1}X_{n}}$ be a positively graded set. 
Let $K [{\mathcal F}_{\infty},X]
:=\bigoplus _{n\geq 1} K [{\mathcal F}_{n,X}]$, where ${\mathcal F}_{n,X}$ is the set of pairs $(f;x_{1},\dots ,x_{r})$ with $f$ a map
from $\{ 1,\dots ,n\}$ to  $\{ 1,\dots ,r\}$, and $x_{i}\in X_{n_{i}}$, for $n_{i}:=\vert f^{-1}(i)\vert$.  The operations $\bullet _{i}$ are defined, using the shuffle algebra structure of $K [{\mathcal F}_{\infty},X]$ defined in example {\bf b)} of \ref{bish}, as follows:}
$$\displaylines {
(f ; x_{1}, \dots , x_{r})\bullet _{i}(g ; y_{1}, \dots , y_{k}):=\hfill \cr
\hfill ((g(1)+r,\dots ,g(i)+r,f(1),\dots ,f(n),g(i+1)+r,\dots ,g(m)+r); x_{1},\dots ,x_{r},y_{1},\dots, y_{k}).\cr } $$

{\rm The subspaces $K[{\mathcal K}_{\infty},X]:={\displaystyle \bigoplus _{n\geq 1}}K [{\mathcal K}_{n,X}]$ and $K
[{\mathcal P}_{\infty},X]:={\displaystyle \bigoplus _{n\geq 1}} K [{\mathcal P}_{n,X}]$ are closed under the operations $\bullet _{i}$, so they are sub-preshuffle algebras of }
$K [{\mathcal F}_{\infty},X]$.
\ms

\subsubsection {\bf Theorem.} For any graded set $X={\displaystyle \bigcup _{n\geq 1}}X_{n}$, the space $K [{\mathcal K}_{\infty},X]$ with the operations $\bullet _{i}$ 
described above is the free preshuffle algebra spanned by $X$.
\ms

\P {\rm Any element of $x\in X_{n}$ is identified with the pair $(\xi _{n};x)$. The result follows easily using that any element $z$ in the free preshuffle algebra 
spanned by $X$ is of the form $x\bullet _{i}y$, with $x\in X_{r}$ for some $1\leq r$ and $y$ an element of the free preshuffle algebra such that $\vert y\vert <\vert z 
\vert $. \hfill $\diamondsuit$}

{\rm Moreover, given a coassociative coproduct $\Theta :K[X]\longrightarrow K[X]\otimes K[X]$, the coproduct $\Delta _{\theta}:K [{\mathcal P}_{\infty},X]
\longrightarrow K [{\mathcal P}_{\infty},X]\otimes K [{\mathcal P}_{\infty},X]$, defined in Example {\bf c)} of 
\ref{exbial}, restricts to $K [{\mathcal K}_{\infty},X]$. So, any free preshuffle algebra is a preshuffle bialgebra.}
\ms

{\rm Note that the free preshuffle algebra spanned by one element of degree one is just the space  $K[S_{\infty }]:={\displaystyle \bigoplus _{n\geq 1}}K [S_{n}]$, equipped 
with the products $$\sigma \bullet _{i}\tau :=(\tau (1)+n,\dots ,\tau (i)+n,\sigma (1),\dots ,\sigma (n),\tau (i+1)+n,\dots ,\tau (m)+n),$$
{\rm for} $\sigma \in S_{n}$, $\tau \in S_{m}$ and $0\leq i\leq m$. Its associated Hopf algebra is the Malvenuto-Reutenauer bialgebra.}
\bs

\noi {\rm {\bf b) Infinitesimal bialgebras} Given a graded nonunital infinitesimal bialgebra $(A,\cdot ,\Delta )$, example {\bf b)} of \ref{exbial} shows that there exists a natural way to define a shuffle bialgebra structure on 
$A$, where the coproduct is $\Delta$ and the operations $\bullet _{\gamma }$ are constructed using $\cdot$ and $\Delta$. It is easy to see that the preshuffle algebra 
structure on $A$, given by $\bullet _{i}=\bullet _{\omega _{i}^{n,m}}$ is in fact a grafting algebra. So, any graded nonunital infinitesimal bialgebra gives rise to a grafting 
bialgebra.}
\bs

\noi {\rm {\bf c) The algebra of planar trees.} The graded vector space $K[T_{\infty}]$ spanned by the set of planar trees $T_{\infty}:={\displaystyle \bigcup _{n\geq 1}}T_{n}$, 
with the products $\circ _{i}$ described in Definition \ref{injerto} is a grafting algebra. Moreover, the subspace $K[Y_{\infty}]$,
spanned by the set of planar binary trees, is a grafting subalgebra of $K[T_{\infty}]$. 
\ms

For any graded set $X={\displaystyle \bigcup _{n\geq 1}}X_{n}$, let $T_{n,X}$ be the set of planar rooted trees with $n+1$ leaves and the internal vertices coloured by
 the elements of $X$, in such a way that a vertex with $r+1$ inputs is coloured by an element $x$ of $X_{r}$. The grafting structure of $K[T_{\infty}]$ induces 
a grafting algebra structure on the space $K[T_{\infty,X}];={\displaystyle \bigoplus _{n\geq 1}}K[T_{n,X}]$ in an obvious way.} 
\ms

{\rm To describe the free grafting algebra spanned by a graded set $X$, we need the following result. Its proof is straightforward.}
 
\begin{lem}\label{freegraft} Let $A$ be a grafting algebra. Let $x_{0},\dots ,x_{q},y,z$ be elements of of $A$ with $\vert x_{i}\vert =n_{i}$, $\vert y\vert =m$ and 
$\vert z\vert =r$. For any family of integers $0\leq i_0<\dots
<i_r\leq r $ and $0\leq j\leq n_{q} $ we have that:
$$x_{0}\bullet _{i_0}(x_{1}\dots (x_{q-1}\bullet _{i_{q-1}}((y\bullet _jx_{q})\bullet _{i_q}z))) = 
(-1)^{m (n_{0} +\dots +n_{q-1} )}y\bullet _i (x_0\bullet _{i_0}(\dots (x_{q}\bullet _{i_q}z))),$$
where $i = j + {\displaystyle \sum _{l=0}^{q-1}} n_{l} + i_q1$.
\end{lem}
\ms

\bs

\begin{thm}  For any graded set $X$, the vector space  $K[T_{\infty ,X}]$, equipped with the linear maps $\circ _i$,
for $i\geq 0$, is the free grafting algebra spanned by $X$.
\end{thm}
\ms

\P {\rm Let ${\mbox {\it Graft}(X)}$ denote the free grafting algebra spanned by $X$. For any element $x\in X_{k}$, the tree $\cc _k$ with its vertex coloured by $x$ is denoted 
by $(\cc _k,x)$. Let $\iota : X\rightarrow K[T_{\infty ,X}]$ be the map which sends an element $x\in X_k$ to $(\cc _k,x)$. 
\ms

\noi From the definition of grafting it is immediate to check that for any $z$ in 
${\mbox {\it Graft}(X)}\setminus X$ there exist elements $z_{1},\dots ,z_{r}$ in ${\mbox {\it Graft}(X)}$ and $x\in X_{n}$ such that:
$$z=z_{1}\bullet _{i_{1}}(z_{2}\bullet _{i_{2}}(\dots (z_{r}\bullet _{i_{r}}x))).$$
 Moreover the integer $r$, the elements $z_{1},\dots ,z_{r}$ and the collection 
$i_{1},\dots ,i_{r}$ are unique. 

\noi The homomorphism $\kappa : {\mbox {\it Graft}(X)}\longrightarrow K[T_{\infty ,X}]$ is defined by induction on the number of elements of $X$ which appear in $z$, this number is denoted by $o(z)$. If $o(z)=1$, then 
$\kappa (z):=\iota (z)$. If $z=z_{1}\bullet _{i_{1}}(z_{2}\bullet _{i_{2}}(\dots (z_{r}\bullet _{i_{r}}x)))$, with $r>0$, then $o(z_{i})<
o(z)$ for all $1\leq i\leq r$. Define
$$\kappa (z):=(t^0,\dots ,t^{n})\circ \iota(x),$$
where $t^{j}=\begin{cases}\kappa (z_{i_{l}})&\ {\rm if}\ j=i_{l},\\
\vert&\ {\rm if}\ j\notin\{i_{1},\dots ,i_{r}\}.\end{cases}$

\noi Lemma \ref{freegraft} implies that $\kappa$ is a homomorphism of grafting algebras.
\ms

Conversely, since any planar rooted tree $t$ with the vertices coloured with the elements of $X$, may be written in a unique way as 
$(t^{0},\dots ,t^{n})\circ (\cc _{n},x)$, with $x\in X_{n}$, the inverse of $\kappa$ is defined by the conditions:\begin{enumerate}
\item $\kappa ^{-1}(\cc _{n},x):=x$, 
\item $\kappa ^{-1}(t):=\kappa ^{-1}(t^{i_{1}})\bullet _{i_{1}}(\dots (\kappa ^{-1}(t^{i_{k}}\bullet _{i_{k}}x))),$
where $t^{i_{1}},\dots ,t^{i_{k}}$ are the trees in $\{ t^{0},\dots ,t^{n}\}$ which are different from $\vert$.\hfill $\diamondsuit$
\end{enumerate} }
\ms
  
\begin{cor} The free grafting algebra on one generator is $(K[Y_{\infty}],\circ _{i})$, while $(K[T_{\infty}],\circ _{i})$ is the free grafting algebra on the graded 
set $\{\cc_{n}\}_{n\geq 1}$, which has exactly one element of degree $n$: the tree $\cc _{n}$, with $n+1$ leaves and a unique vertex.
\end{cor}

{\rm On the graded vector space $K[Y_{\infty}]$, spanned by all planar binary rooted trees, the coproduct 
$\Delta_{PR}$ is defined as the unique counital coproduct such that:}
$$\displaylines {
\Delta _{PR}(\arboluno)= 0,\hfill \cr
\Delta _{PR}(t\vee w):=\sum t_{(1)}\otimes (t_{2}\vee w)+t\otimes (\vert \vee w)+(t\vee \vert )\otimes w+\sum (t\vee w_{(1)})\otimes w_{(2)}.\cr }$$
{\rm The vector space $K[Y_{\infty}]$, with the products $\circ _{i}$ and $\Delta _{PR}$ is a  grafting bialgebra.}
\bs

{\rm Let $X={\displaystyle \bigcup _{n\geq 1}}X_{n}$ be a graded set. If $\Theta $ is a graded coassociative coproduct on $K[X]$, then there exists a coproduct $\Delta _{\theta }$
 on $K[T_{\infty ,X}]$ given by:}
$$\displaylines {
\Delta _{\theta }(\cc _{n},x):=\sum _{1\leq i\leq n-1}(\sum _{\vert x_{(1)}\vert =i}(\cc_{i},x_{(1)})\otimes (\cc_{n-i},x_{(2)})),\ {\rm for}\ x\in X_{n}\ 
{\rm and}\ n\geq 1,\hfill\cr
\Delta _{\theta }((t^0,\dots ,t^{n})\circ (\cc _n,x)):=\hfill \cr 
\sum _{0\leq i\leq n}\bigl( \sum _{\vert x_{(1)}\vert=\vert t^{0}\vert +\dots +\vert t_{(1)}^{i}\vert}
(t^{0},\dots ,t_{(1)}^{i})\circ (\cc _{n_{j}},x_{(1)})\otimes (t_{(2)}^{i},\dots ,t^{n})\circ (\cc_{n-n_{j}},x_{(2)})+\cr
\sum _{\vert x_{(1)}\vert =\vert t^{0}\vert +\dots +\vert t^{i}\vert}(t^{0},\dots ,t^{i})\circ (\cc _{n_{j}},x_{(1)})\otimes 
(t^{i+1},\dots ,t^{n})\circ (\cc_{n-n_{j}},x_{(2)})\bigr ),\cr}$$
{\rm where} $\vert x_{(1)}\vert =n_{j}$, $\Delta _{\theta }(t^{i})=\sum t_{(1)}^{i}\otimes t_{(2)}^{i}${\rm , and} $\Theta (x ):=\sum x_{(1)}\otimes x_{(2)}$.
\ms

\noi {\rm It is not difficult to check that, for any $\Theta $, the space $K[T_{\infty}]$ equipped with the operations $\circ _{i}$ and $\Delta _{\theta }$ is a 
grafting bialgebra.}
\bs

\noi {\rm {\bf d) The space of Hochschild cochains.} (see \cite{Ger} ) Let $A$ be a unital $K$-algebra, and let 
$C^*(A):=\bigoplus _{n\geq 0}{\rm Hom}_{K}(A^{\otimes n},A)$ be the space of Hochschild cochains on $A$. 

\noi The space $C^*(A)[1]:={\displaystyle \bigoplus _{n\geq 0}}{\rm Hom}_{K}(A^{\otimes (n+1)},A)$ is a grafting algebra with the operations $\bullet _i$ defined as follows:
$$g\bullet _i f:=f\circ (id_{A}^{\otimes (i-1)}\times g\times id _{A}^{\otimes (n-i)}),$$
for $g\in C^m(A,A)$ and $f\in C^n(A,A)$.}
\ms

\noi {\rm Consider on $C^*(A)[1]$ the following coproduct:}
$$\Delta (f):=\sum _{i=1}^{n-1}f_{(1)}^{i}\otimes f_{(2)}^{n-i+1},\ {\rm for}\ f\in C^n(A,A),$$
{\rm where} \begin{enumerate} 
\item $f_{(1)}^{i}(x_{1},\dots ,x_{i}):=f(x_{1},\dots ,x_{i},1_{A},\dots ,1_{A})\in C^i(A,A)$
\item $f_{(2)}^{n-i+1}(x_{1},\dots ,x_{n+1-i}):=f(1_{A},
\dots ,1_{A},x_{1},\dots ,x_{n+1-i})\in C^{n+1-i}(A,A)$.\end{enumerate}
\ms

{\rm It is easy to see that }$(C^*(A)[1],\bullet _{i},\Delta )$ {\rm is a grafting bialgebra.}
\bs

\noi {\rm {\bf e) The underlying space of an algebraic operad.}  Let $K$ be a field of characteristic $0$, and let $\mathbb P$ be a $K$-linear operad as described in 
\cite{GK}. Consider the graded $K$-vector space ${\mathbb P}[1]:={\displaystyle \bigoplus _{n\geq 0}}{\mathbb P}(n+1)$ equipped with the maps:
$$\lambda \bullet _i \nu :=\gamma _{1,\dots, 1,n,1,\dots ,1}(\nu \otimes 1\otimes \dots \otimes 1\otimes \lambda \otimes
1\otimes \dots \otimes 1),$$
where $1\in {\mathbb P}(1)={\mathbb P}[1]_0$ is the identity operation, and $\lambda \in {\mathbb P}(m)$ is at the $i+1$-th place.
It is easy to check that ${\mathbb P}[1]$ with these products is a grafting algebra over $K$.}
\bs

{\rm As an example of grafting bialgebra consider the grafting algebra associated to the operad $As$.

The grafting structure of $As[1]={\displaystyle \bigoplus _{n\geq 0}}K[S_{n+1}]$ is given by the operations:}
$$(\sigma \bullet _{i}\tau )= (\tau_{(1)}^{i}\times \sigma\times \tau _{(2)}^{m-i-1})\cdot \delta _{i}^{n},$$
{\rm where} $\tau=({\tilde \tau}_{(1)}^{i}\times 1_{1}\times {\tilde \tau }_{(2)}^{m-i-1})\cdot \delta $ with $\delta \in {\mbox {\it Sh}(i,1.m-i-1)}$, 
${\tilde \tau}_{(1)}^{i}\in S_{i}$, 

\noi ${\tilde \tau }_{(2)}^{m-i-1}\in S_{m-i-1}${\rm , and }
$$\delta _{i}^{n}(k):=\begin{cases}\delta (k),&\ {\rm for}\ \delta (k)\leq i\ {\rm and}\ k<\delta ^{-1}(i+1),\\
\delta (k)+n-1,&\ {\rm for}\ \delta (k)> i\ {\rm and}\ k<\delta ^{-1}(i+1),\\
i+r+1,&\ {\rm for}\ k=\delta ^{-1}(i+1)+r\ {\rm and}\ 0\leq r<n,\\
\delta (k-n+1),&\ {\rm for}\ \delta (k)\leq i\ {\rm and}\ k>\delta ^{-1}(i+1)+n-1,\\
\delta (k-n+1)+n-1,&\ {\rm for}\ \delta (k)> i\ {\rm and}\ k>\delta ^{-1}(i+1)+n-1.\end{cases}$$
{\rm In fact, $\sigma \bullet _{i}\tau $ is obtained by replacing $i+1$ in the image of $\tau$ by $(\sigma (1)+ i,\dots ,\sigma (n)+i)$. For instance,}
$$(2,4,1,3)\bullet _{1}(1,3,2,5,4)=(1,6,3,5,2,4,8,7).$$
\ms

{\rm To define a coproduct on $As[1]$, let $\gamma \in S_{m+1}$ be a permutation, for an integer $0\leq i\leq m$, there exists unique decompositions:}
 $$\gamma = ({\tilde \gamma} _{(1)}^{i+1}\times {\tilde \gamma} _{(2)}^{m-i})\cdot \delta = 
 ({\tilde \gamma} _{(1)}^{i}\times {\tilde \gamma} _{(2)}^{m+1-i})\cdot \epsilon ,$$
 {\rm where} ${\tilde \gamma} _{(i)}^{j}\in S_{j}$, for $i=1,2$, $\delta^{-1}\in {\mbox {\it Sh}(i+1,m-i)}$ and $\epsilon \in {\mbox {\it Sh}(i,m-i+1)}$. 
 {\rm Define }
 $$\Delta _{As} (\gamma ):=\sum _{i=0}^{m}{\tilde \gamma} _{(1)}^{i+1}\otimes {\tilde \gamma} _{(2)}^{m+i-i}.$$
 \ms
 
 \subsubsection {\bf Proposition } The space $As[1]={\displaystyle \bigoplus _{n\geq 1}}K[S_{n+1}]$, equipped with the operations $\bullet _{i}$ and the coproduct 
 $\Delta _{As}$ is a grafting bialgebra.
 \ms
 
 \P {\rm We know that $(As[1], \bullet _{i})$ is a grafting algebra. To check that $\Delta $ is coassociative, it suffices to note that, for $\gamma \in S_{m+1}$, 
we have that:}
 $$(\Delta _{As}\otimes id_{As[1]})\circ \Delta _{As}(\gamma )=\sum _{i+j+k=m}{\tilde \gamma} _{(1)}^{i+1}\otimes {\tilde \gamma} _{(2)}^{j+1}\otimes 
{\tilde \gamma} _{(3)}^{k+1}= (id_{As[1]}\otimes \Delta _{As})\circ \Delta _{As} (\gamma ),$$
 {\rm where, for each compositions $(i,j,k)$ of $m$, the following equalities hold:}
 $$\gamma = ({\tilde \gamma}_{(1)}^{i+1}\times {\tilde \gamma} _{(2)}^{j}\times {\tilde \gamma} _{(3)}^{k})\cdot \delta _{1}=
  ({\tilde \gamma} _{(1)}^{i}\times {\tilde \gamma}_{(2)}^{j+1}\times {\tilde \gamma} _{(3)}^{k})\cdot \delta _{2}=
  ({\tilde \gamma} _{(1)}^{i}\times {\tilde \gamma} _{(2)}^{j}\times {\tilde \gamma} _{(3)}^{k+1})\cdot \delta _{3},$$ 
 {\rm with } ${\tilde \gamma}_{(l)}^{p}\in S_{p}$, for $l=1,2,3$, $\delta _{1}\in {\mbox {\it Sh}(i+1,j,k)}$, $\delta _{2}\in {\mbox {\it Sh}(i,j+1,k)}$ 
 {\rm and}  $\delta _{3}\in {\mbox {\it Sh}(i,j,k+1)}$. 
\ms 
 
\noi  {\rm To prove the relationship between $\Delta _{As}$ and the operations $\bullet _{i}$, note that for any $\gamma \in S_{n}$ and any $0\leq i\leq n$, there 
 exist unique order preserving bijections 
  $\varphi _{(1)}^{i}:\{ 1,\dots ,i\}\longrightarrow \gamma ^{-1}(\{1,\dots ,i\})$ and 
 $\varphi _{(2)}^{n-i}:
 \{ 1,\dots ,n-i\}\longrightarrow \gamma ^{-1}(\{ i+1,\dots ,n\})$. The permutations ${\tilde \gamma}_{(1)}^{i}$ and ${\tilde \gamma }_{(2)}^{n-i}$ are 
 given by the formulas:
 $$\displaylines {
{\tilde \gamma}_{(1)}^{i}=(\gamma (\varphi _{(1)}^{i}(1)),\dots ,\gamma (\varphi _{(1)}^{i}(i)))\hfill \cr
{\tilde \gamma }_{(2)}^{n-i}=(\gamma (\varphi _{(2)}^{n-i}(1)),\dots ,\gamma (\varphi _{(2)}^{n-i}(n-i))).\hfill \cr  }$$

Using the formulas above, it is easily seen that, for $\sigma \in S_{n+1}$, $\tau \in S_{m+1}$ and $0\leq j\leq n+m$, we have that:
$$({\tilde {\sigma \bullet _{i}\tau} })_{(1)}^{j+1}=\begin{cases}{\tilde \tau} _{(1)}^{j+1},& {\rm for}\ 0\leq j<i\\
{\tilde \sigma} _{(1)}^{j-i+1}\bullet _{i}{\tilde \tau }_{(1)}^{i+1},& {\rm for}\ i\leq j\leq i+n\\
\sigma \bullet _{i}{\tilde \tau }_{(1)}^{j-n+1},&{\rm for}\ i+n<j\leq n+m.\end{cases}$$
$$({\tilde {\sigma \bullet _{i}\tau} })_{(2)}^{n+m-j+1}=\begin{cases}\sigma\bullet _{i-j}{\tilde \tau} _{(2)}^{m-j+1},& {\rm for}\ 0\leq j<i\\
{\tilde \sigma} _{(2)}^{n+i-j+1}\bullet _{0}{\tilde \tau }_{(2)}^{m-i+1},& {\rm for}\ i\leq j\leq i+n\\
{\tilde \tau }_{(2)}^{m+n-j+1},&{\rm for}\ i+n<j\leq n+m,\end{cases}$$
which ends the proof.}\hfill $\diamondsuit $
\end{exmpls}
\bs

The following result is an extension of Lemma \ref{hadsh} to preshuffle algebras, its proof is straightforward.. 

\begin{lem}\label{hadpre}  Let $(A, \bullet _{i})$ and $(B,\circ _{j})$ be preshuffle (respectively grafting) algebras.\begin{enumerate}
\item The Hadamard product $A\underset {H}{\otimes}B$ has a natural structure of preshuffle (respectively grafting) algebra, given by the operations:
$$(x \otimes  y)\bullet _{i}(x\rq \otimes y\rq ):= (x\bullet _{i}x\rq ) \otimes(y\circ _{i}y\rq),$$ 
for $x\in A_{n}, y\in B_{n}, x\rq \in A_{m}, y\rq \in B_{m}$ and $0\leq i\leq m$.
\item The tensor product $A\otimes B$  has a natural structure of preshuffle (respectively grafting) algebra, given by the operations:
$$(x\otimes y)\bullet _{i}(x\rq \otimes y\rq ):=\begin{cases}
(x\bullet _{i}x\rq )\otimes (y\circ _{i-\vert x\rq \vert }y\rq ),& \ {\rm for}\ 
\ i=\vert x\rq \vert,\\
0,& {\rm otherwise.}\end{cases}$$
\end{enumerate}
\end{lem}
\ms

\begin{prop} Let $(A,\bullet _{i}, \Delta _{A})$ and $(B,\circ _{j}, \Delta_{B} )$ be two preshuffle (respectively grafting) bialgebras. The Hadamard product 
$A\underset {H}{\otimes}B$ with the operations $\bullet _{i}$ given in Definition \ref{hadpre} and the coproduct given by:
$$\Delta _{A\underset {H}{\otimes}B}(x\otimes y)=\sum _{\vert x_{(1)}\vert =\vert y_{(1)}\vert }(x_{(1)}\otimes y_{(1)})\otimes (x_{(2)}\otimes y_{(2)}),$$
where $\Delta _{A}(x)=\sum x_{(1)}\otimes x_{(2)}$ and $\Delta _{B}(y)=\sum y_{(1)}\otimes y_{(2)}$,
is a preshuffle (respectively grafting) bialgebra.
\end{prop}
\ms

\P  Let $x\in A_{n}$, $y\in B_{n}$, $x\rq \in A_{m}$ and $y\rq \in B_{m}$. Since $\vert x_{(1)}\vert <\vert y\vert <\vert x\bullet _{0}x\rq _{(1)}\vert $ and 
$\vert y_{(1)}\vert <\vert x\vert <\vert y\circ _{0}y\rq _{(1)}\vert $, for all $x_{(1)}$ and $y_{(1)}$, we have that:
$$\displaylines {
\Delta _{A\underset {H}{\otimes}B}((x\bullet _{0}x\rq )\otimes (y\circ _{0}y\rq ))=\hfill\cr
\sum _{\vert x_{(1)}\vert =\vert y_{(1)}\vert }(x_{(1)}\otimes y_{(1)})
\otimes ((x_{(2)}\bullet _{0}x\rq )\otimes (y_{(2)}\circ _{0}y\rq ))+\hfill \cr
\hfill (x\otimes y)\otimes (x\rq \otimes y\rq )+\sum _{\vert x\rq _{(1)}\vert =\vert y\rq _{(1)}\vert }((x\bullet _{0}x\rq _{(1)})\otimes (y\circ _{0}y
\rq _{(1)}))\otimes (x\rq _{(2)}\otimes y\rq _{(2)})=\cr
\sum ((x\otimes y)_{(1)}\otimes ((x\otimes y)_{(2)}\bullet _{0}(x\rq \otimes y\rq ))_{(2)}+(x\otimes y)\otimes (x\rq \otimes y\rq )+\hfill \cr
\hfill +\sum ((x\otimes y)\bullet _{(0)}(x\rq \otimes y\rq )_{(1)})\otimes (x\rq \otimes y\rq )_{(2)}.\cr }$$

For $0< i< m$, note that $\vert ( x\bullet _{i}x\rq) _{(1)}\vert  =\vert ( y\circ _{i}y\rq) _{(1)}\vert $ only in the following cases:
\begin{enumerate}\item  $(x\bullet _{i}x\rq) _{(1)}= x\rq _{(1)}$,\  $ (y\circ _{i}y\rq) _{(1)}=y\rq _{(1)}$ and $\vert x\rq _{(1)}\vert =\vert y\rq _{(1)}\vert 
\leq i$,
\item $(x\bullet _{i}x\rq) _{(1)}= x_{(1)}\bullet _{i}x\rq _{(1)}$, $ (y\circ _{i}y\rq) _{(1)}=y_{(1)}\circ _{i}y\rq _{(1)}$, $\vert x\rq _{(1)}\vert =
\vert y\rq _{(1)}\vert =i$, and $\vert x\rq _{(1)}\vert =\vert y\rq _{(1)}\vert $,
\item $(x\bullet _{i}x\rq) _{(1)}= x\bullet _{i}x\rq _{(1)}$, $(y\circ _{i}y\rq) _{(1)}=y\circ _{i}y\rq _{(1)}$, and $\vert x\rq _{(1)}\vert =
\vert y\rq _{(1)}\vert \geq i$.
\end{enumerate}
The equalities above imply that $\Delta _{A\underset {H}{\otimes}B}((x\bullet _{i}x\rq )\otimes (y\circ _{i}y\rq ))=$
$$\displaylines {
\sum _{\vert (x\rq \otimes y\rq )_{(1)}\vert \leq i}(x\rq \otimes y\rq )_{(1)}\otimes ((x\otimes y)\circ _{i}(x\rq \otimes y\rq)_{(2)}) +\hfill\cr
\sum _{\vert (x\rq \otimes y\rq )_{(1)}\vert = i}((x\otimes y)_{(1)}\circ _{i}(x\rq \otimes y\rq)_{(1)})\otimes ((x\otimes y)_{(2)}\bullet _{0} 
(x\rq \otimes y\rq )_{(2)})+ \cr
\hfill \sum _{\vert (x\rq \otimes y\rq )_{(1)}\vert \geq i}((x\otimes y)\circ _{i}(x\rq \otimes y\rq)_{(1)})\otimes 
(x\rq \otimes y\rq )_{(2)},\cr}$$
which ends the proof for $0<i<m$. The result for $i=m$ is obtained in the same way.
\hfill $\diamondsuit $
\bs

\bs

Let $(A, \bullet _{i},\Delta )$ be a grafting bialgebra, we want to show that there exist a natural way of defining operations $\bullet _{\gamma}$ on $A$ is such a way that 
$(A,\bullet _{\gamma},\Delta )$ is a shuffle bialgebra.

\begin{enumerate}\item  Given a composition $(n_{1},\dots ,n_{p})$ of $n$, we denote by 
$\Delta _{n_{1},\dots ,n_{p}}$ the composition $\pi _{n_{1},\dots ,n_{p}}\circ \Delta ^{p}$,
where $\pi _{n_{1},\dots ,n_{p}}$ is the projection from $A^{\otimes p}$ to 

\noi $A_{n_{1}}\otimes \dots \otimes A_{n_{p}}$.
\item  Let $\gamma $ be an {\mbox {\it (n,m)}}-shuffle. There exist unique compositions $(n_1,\dots ,n_r)$ of $n$ and $(m_1,\dots ,m_{r+1})$ of $m$ such that 
$m_1\geq 0$, $m_{r+1}\geq 0$, $m_i\geq 1$ for $2\leq i\leq r$, and $n_j\geq 1$ for $1\leq j\leq r$, such that 
$$\gamma  =(n+1,\dots n+m_{1},1,\dots ,n_1,n+m_1+1,\dots ,n+m_1+m_2,n_1+1,\dots ,n_1+n_2,\dots , n+m),$$
that is $$\gamma (j)=\begin{cases}j+n-\sum_{i=1}^kn_i,&\ {\rm for}\ 0<j-\sum _{i=1}^kn_i+m_i\leq m_{k+1},\ {\rm with}\ 0\leq k\leq r\\
j-\sum_{i=1}^km_i,&\ {\rm for}\ m_{k+1}<j-\sum _{i=1}^kn_i+m_i\leq m_{k+1}+n_{k+1},\ {\rm with}\ 1\leq k\leq r.\end{cases}$$

\noi For instance, if $\gamma =(1,3,4,2,5,6)\in {\mbox {\it Sh}(2,4)}$, then $(m_1,m_2,m_3)=(0,2,2)$, and $(n_1,n_2)=(1,1)$.

\noi Given elements $x\in A_{n}$ and $y\in A_{m}$, define the element $x\bullet _{\gamma}y\in A_{n+m}$ as follows:
$$x\bullet _{\gamma }y:=\sum x_{(1)}^{n_{1}}\bullet _{m_1}(\dots (x_{(r-1)}^{n_{r-1}}\bullet _{m_1+\dots +m_{r-1}}(x_{(r)}^{n_{r}}\bullet _{m_1+\dots +m_{r}}y)))),$$
where $\Delta _{n_{1},\dots ,n_{p}}(x)=\sum x_{(1)}^{n_{1}}\otimes \dots \otimes x_{(p)}^{n_{p}}$.
\end{enumerate}
\ms

\begin{thm}\label{grftbi} Let $(A, \bullet _{i},\Delta )$ be a grafting bialgebra. The graded space $A$ equipped with the operations $\bullet _{\gamma}$ defined above for any 
shuffle $\gamma $, is a shuffle bialgebra.
\end{thm}

\P  Let $x\in A_{n}$, $y\in A_{m}$ and $z\in A_{r}$ be homogeneous elements of $A$, and let $\gamma \in {\mbox {\it Sh}(n,m+r)}$, $\delta \in {\mbox {\it Sh}(m,r)}$, 
$\lambda \in {\mbox {\it Sh}(n+m,r)}$ and $\sigma \in {\mbox {\it Sh}(n,m)}$ be such that 
$$(1_{n}\times \delta )\cdot \gamma=(\sigma \times 1_{r})\cdot \lambda.$$
We want to verify that $x\bullet _{\gamma} (y\bullet _{\delta }z)=(x\bullet _{\sigma }y)\bullet _{\lambda }z$.

\noi Let $\gamma $ be the permutation given by the integers
$(n_1,\dots ,n_p)\vdash n$ and $(h_1,\dots ,h_{p+1})\vdash m+r$. We proceed by a recursive argument on $p$.

\noi If $p=1$, then $\gamma =\omega _{h}^{n,m+r}$. 

\noi Suppose that 
$\delta \in {\mbox {\it Sh}(m,r)}$ is given by integers $(m_1,\dots ,m_q)\vdash m$ and $(r_1,\dots ,r_{q+1})\vdash r$, we have to consider two different cases.

\noi {\bf a)} If there exists $0\leq k\leq q$ such that $0< h-{\displaystyle \sum _{i=1}^k}r_i+m_i<r_{k+1}$,
then 

\noi $(1_n\times \delta )\cdot \gamma = (\sigma\times 1_r)\cdot \lambda$, where $\sigma =\omega _{m_1+\dots +m_k}^{n,m}$ and $\lambda $ is the $(n+m,r)$ shuffle 
associated to the compositions $(m_1,\dots ,m_k,n,m_{k+1},\dots ,m_q)$ 
of $n+m$ and $(r_1,\dots ,r_k,{\underline h},r_{k+1}-{\underline h},r_{k+2},\dots ,r_{q+1})$ of $r$, with ${\underline h}:=h-{\displaystyle \sum _{i=1}^k}r_i+m_i$.
\medskip 

\noi Applying the properties of a grafting algebra, we get that:
$$\displaylines { 
x\bullet _{\gamma }(y\bullet _{\delta}z)=\hfill \cr
x\bullet _h(y_{(1)}^{m_1}\bullet _{r_1}(\dots (y_{(q)}^{m_q}\bullet _{r_1+\dots +r_q}z)))=\cr
y_{(1)}^{m_1}\bullet _{r_1}(\dots y_{(k)}^{m_k}\bullet _{r_1+\dots +r_k}(x\bullet _{ h-\sum_{1\leq i\leq k}m_i}(y_{(k+1)}^{m_{k+1}}\bullet _{r_1+\dots +r_{k+1}}(\dots (y_{(q)}^{m_q}\bullet _{r_1+\dots +r_q}z))))).\cr} $$
Since $(x\bullet _{m_1+\dots +m_k}y)_{(j)}^{m_j}=\begin{cases}y_{(j)}^{m_j},&\ {\rm for}\ 1\leq j\leq k,\\
x,&\ {\rm for}\ j=k+1\\
y_{(j-1)}^{m_{j-1}},&\  {\rm for}\ j<k+1,\end{cases}$,

\noi we get the result.
\medskip

\noi {\bf b)} If there exists $0\leq k\leq q$ such that $0< h-{\displaystyle \sum _{i=1}^{k-1}}(r_i+m_i) +r_k<m_{k}$, then 

\noi $\sigma =\omega _{m_1+\dots +m_{k-1}+{\underline h}}
^{n,m}$, with  ${\underline h}:=h-{\displaystyle \sum _{i=1}^{k-1}}(r_i+m_i)-r_k$, and $\lambda $ is the 
$(n+m, r)$-shuffle associated to the compositions $(m_1,\dots ,m_{k-1},m_k+n,\dots ,m_p)$ of $n+m$ and $(r_1,\dots ,r_{p+1})$ of $r$.

\noi We have that \begin{enumerate}
\item $(x\bullet _{m_1+\dots +m_{k-1}+{\underline h}}y)_{(j)}^{m_j}=y_{(j)}^{m_j}$, for $j\neq k$, 
\item  $(x\bullet _{m_1+\dots +m_{k-1}+{\underline h}}y)_{(k)}^{n+m_k}=
x\bullet _{\underline h}y_{(k)}^{m_k}$.
\end{enumerate}
\ms

\noi  In this case, using the properties of grafting algebras, it is immediate to check that:
$$x\bullet _{\gamma }(y\bullet _{\delta}z)=(x\bullet _{m_1+\dots +m_{k-1}+{\underline h}}y)\bullet _{\lambda }z.$$
\medskip

\noi For $p>1$, note that if $\gamma$ is the $(n,m+r)$-shuffle associated to the compositions $(n_1,\dots ,n_p)\vdash n$ and $(h_1,\dots , h_{p+1})$ of $m+r$, then
$$x\bullet _{\gamma }(y\bullet _{\delta }z)=x_{(1)}^{n_1}\bullet _{h_1}(x_{(2)}^{n-n_1}\bullet _{\tilde {\gamma}}(y\bullet _{\delta }z),$$
where ${\tilde {\gamma}}$ is the $(n-n_1,m+r)$-shuffle associated to the compositions 
$(n_2,\dots ,n_p)$ of $n-n_1$ and $(h_1+h_2,\dots , h_{p+1})$ of $m+r$.

\noi We get that:
$$
x\bullet _{\gamma }(y\bullet _{\delta }z)=x_{(1)}^{n_1}\bullet _{h_1}((x_{(2)}^{n-n_1}\bullet _{\tilde {\sigma}}y)\bullet _{\tilde {\lambda} }z)=
(x_{(1)}^{n_1}\bullet _{k_1}(x_{(2)}^{n-n_1}\bullet _{\tilde {\sigma}}y))\bullet _{\lambda }z,$$
where $(1_{n-n_1}\times \delta )\cdot {\tilde {\gamma}}=({\tilde {\sigma}}\times 1_r)\cdot {\tilde {\lambda}}$, and  

\noi $(1_{n+m-n_1}\times {\tilde {\lambda} })\cdot \omega _{h_1}^{n_1,n+m+r-n_1}=(\omega _{k_1}^{n_1,n+m-n_1}\times 1_r)\cdot \lambda$. 

\noi So, we have that $$x\bullet _{\gamma }(y\bullet _{\delta }z)=(x\bullet _{\sigma}y)\bullet _{\lambda}z,$$
with $\sigma :=(1_{n_1}\times {\tilde {\sigma}})\cdot \omega _{k_1}^{n_1,n+m-n_1}$ and $(1_n\times \delta )\cdot \gamma =(\sigma \times 1_r)\cdot \lambda $, which ends the proof.\hfill
 $\diamondsuit $
\bs

\bs

\noi{\bf Primitive elements of preshuffle bialgebras.} 
\ms

Since any shuffle bialgebra is a preshuffle algebra, we look for the operations obtained by compositions and linear combinations of the primitive operations
 $\{-,-\}$ and $B^{\gamma}$, introduced in Section 4, which can be defined in terms of the multiplications $\bullet _{i}$ of a preshuffle algebra.
 
 \noi Let $(A,\bullet _{i})$ be a preshuffle algebra, and let $x\in A_{n}$, $y\in A_{m}$ and $z\in A_{r}$ be elements of $A$. 
 Note that $\{x,y\}=x\bullet _{top}y-x\bullet _{0}y$ and $B^{\omega _{i}^{n,m}}(x;y)=x\bullet _{i}y$ are defined for all $1\leq i\leq m-1$. But also the element 
 $$B^{1_{n}\times \omega _{i}^{m,r}}(\{x,z\};y)=z\bullet _{i+n}(x\bullet _{0}y)-x\bullet _{0}(z\bullet _{i}y)$$
 may be defined in $A$ for $1\leq i\leq m$. In a similar way, the element
 $$\displaylines {
 B_{q}^{1_{n_{1}}\times \omega _{n_{\geq 2}+i}^{m,n_{\geq 2}+r}}(\{x_{1},z\};x_{2},\dots ,x_{q};y)=\hfill\cr
 z\bullet _{n+i}(x_{1}\bullet _{0}x_{2}\bullet_{0}\dots x_{q}\bullet _{0}y)-x_{1}\bullet _{0}(z\bullet _{n_{\geq 2}+i}
 (x_{2}\bullet _{0}\dots \bullet _{0}x_{q}\bullet_{0}y)),\cr }$$
where $\vert x_{i}\vert =n_{i}$, $n={\displaystyle \sum _{i=1}^q}n_i$ and $n_{\geq k}:={\displaystyle \sum _{i=k}^q}n_i$, may be defined on $A$.

\begin{defn}\label{primpre}{\rm Let $(A,\bullet _{i})$ be a preshuffle algebra over $K$. For $q\geq 0$ and $1\leq p\leq n_{r}$, the $q+2$-ary operation 
$L_{q}^{p}: A^{\otimes q}\otimes A\longrightarrow A$ is defined by the following formulas:
$$\displaylines {
L_{q}^{p}(x_{1},\dots ,x_{q};y;z):=\hfill \cr
z\bullet _{p}y,\ \ \ \ \ \ \ \ \ \ \ {\rm for}\ q=0\ {\rm and}\ 0<p<\vert y\vert,\cr
\{ y,z\}=y\bullet _{top}z -y\bullet _{0}z,\ \ \ \ \ \ \ \ \ \ \ {\rm for}\ q=0\ {\rm and}\ p=\vert y\vert ,\cr
z\bullet _{p+n}(x_{1}\bullet _{0}\dots \bullet _{0}x_{q}\bullet _{0}y) - x_{1}\bullet _{0}(z\bullet _{p+n_{\geq 2}}(x_{2}\bullet _{0}\dots 
\bullet _{0}x_{q}\bullet _{0}y)),\ \ \ \ \ {\rm for}\ q\geq 1,\cr }$$
where $(x_{1},\dots ,x_{q};y;z)$ denotes the element $x_{1}\otimes \dots \otimes x_{q}\otimes y\otimes z\in A^{\otimes (q+2)}$, 

\noi $n_{k} :=\vert x_{k}\vert $, $n_{\geq k}:={\displaystyle \sum _{i=k}^qn_i}$ and $n=n_{\geq 1}$, for $1\leq k\leq q$.}
\end{defn}
\bs 

The following result implies that the subspace of primitive elements of a preshuffle bialgebra is closed under the ${\mbox q+1}$-ary operations $L_{q}^{p}$, its proof is similar to 
Proposition \ref{elsrim} one.
\ms

\begin{prop}\label{elprim} Let $(A={\displaystyle\bigoplus _{k\geq 1}}A_{k},\bullet _{i}, \Delta )$ be a preshuffle bialgebra. If the elements $x_{1},\dots x_{q},y,z$ belong to 
${\mbox {\rm Prim}(A)}$, then $L_{q}^{p}(x_{1},\dots ,x_{q};y;z)$ belongs to ${\mbox {\rm Prim}}(A)$, for any $1\leq p\leq \vert x_{q}\vert$.
\end{prop}
\bs

\noi Note that, in fact, $\Delta (L_{0}^{i}(x;y))=0$ for all $y\in A$ and $x\in {\mbox {\rm Prim}(A)}$ and $0<i<\vert x\vert $.
\bs

In order to describe the relationship verified by the new operations $L_{q}^{p}$, we need the following Lemma.

\begin{lem}\label{lem*} Let $x_{1},\dots ,x_{q}, y,z$ be elements of a preshuffle algebra $A$. With the same notations that in Definition \ref{primpre}, the product $\bullet _{0}$ and the operations $L_{q}^{j}$ defined above verify the following equalities:

\noi {\bf 1)} For $j<\vert y\vert$,
$$L_{q}^{j}(x_{1},\dots ,x_{q};y;z\bullet _{0}w)=
\sum _{k=0}^{q}L_{k}^{j+n_{\geq k}}(x_{1},\dots ,x_{k};L_{q-k}^{j}(x_{k+1},\dots ,x_{q};y;w);z),$$ and
$$\displaylines{L_{q}^{\vert y\vert}(x_{1},\dots ,x_{q};y;z\bullet _{0}w)=\hfill \cr
\sum _{k=0}^{q}L_{k}^{\vert y\vert+n_{\geq k+1}}(x_{1},\dots ,x_{k},L_{q-k}^{\vert y\vert}(x_{k+1},\dots ,
x_{q};y;w);z)+L_{q}^{\vert y\vert }(x_{1},\dots ,x_{q};y;z)\bullet _{0}w,\cr }$$
\ms

\noi {\bf 2)} For $1\leq j\leq\vert z\vert $, $L_{q}^{j}(x_{1},\dots ,x_{q},z\bullet _{0}y;w)=L_{q}^{j}(x_{1},\dots ,x_{q};z;w)\bullet _{0}y$, and
$$L_{q}^{j+\vert z\vert }(x_{1},\dots ,x_{q};z\bullet _{0}y;w)=\begin{cases}L_{1}^{j}(z,y;w)+z\bullet _{0}L_{0}^{j}(y;z),&{\rm for}\ q=0\\
L_{q+1}^{j}(x_{1},\dots ,x_{q},z;y;w),&{\rm for}\ q\geq  1,\end{cases}$$
\ms

\noi{\bf 3)} $$L_{q}^{j}(x_{1},\dots ,x_{q-1},z\bullet _{0}x_{q};y;w)=\begin{cases}L_{q+1}^{j}(x_{1},\dots ,x_{q-1},z,x_{q};y;w),&{\rm for}\ q\geq 2,\\
L_{2}^{j}(z,x_{1};y;w)+z\bullet _{0}L_{1}^{j}(x_{1};y;w),&{\rm for}\ q=1.\end{cases}$$
\end{lem}
\ms 

\P The formulas are straightforward to check. We prove for instance the last one, the other ones may be obtained similarly.

\noi For $q\geq 2$, the result is obvious. 

\noi For $q=1$, we have that:
$$\displaylines {
L_{1}^{j}(z\bullet _{0}x;y;w)=w\bullet _{j+n+r}(z\bullet _{0}x\bullet _{0}y)-(z\bullet _{0}x)\bullet _{0}(w\bullet _{j}y)=\hfill\cr
(w\bullet _{j+n+r}(z\bullet _{0}x\bullet _{0}y)-z\bullet _{0}(w\bullet _{j+n }(x\bullet _{0}y))) +
(z\bullet _{0}(w\bullet _{j+n }(x\bullet _{0}y))-(z\bullet _{0}x)\bullet _{0}(w\bullet _{j}y))=\cr
\hfill L_{2}^{j}(z,x;y;w)+z\bullet _{0}L_{1}^{j}(x;y;w).\qquad \diamondsuit\cr }$$
\bs

We introduce some notation, in order to prove the relations satisfied by the operations $L_{q}^{p}$.

\begin{notn} {\rm Let $(A,\bullet _{i})$ be a preshuffle algebra over $K$ and let  ${\bold x}=(x_{1},\dots ,x_{n})$, $y$, ${\bold z}=(z_{1},\dots ,z_{m})$, $t$ and $w$ be a collection of elements in $A$. Given nonnegative integers $0\leq j\leq \vert y\vert$, $0\leq k\leq \vert w\vert $ and $1\leq l\leq m$ we define:

\noi {\bf a)} for a partition
${\underline p}=\{ p_{1}, \dots , p_{m}\}$ of $m$ with $p_i\geq 0$ for $1\leq i\leq m$, the element $L_{\underline p}^j({\bold x},y,{\bold z})$ as follows:
\begin{enumerate}\item if $m=1$, then $L_{\underline p}^j({\bold x},y,z_1):=L_n^j(x_1,\dots ,x_n;y;z)$.
\item it $m>1$, then $$L_{\underline p}^j({\bold x},y,{\bold z}):= L_{p_1}^{j+n_{>p_1}}(x_1,\dots ,x_{p_1};L_{{\underline p}_1}^{j}((x_{k_1+1},\dots ,x_n);y;(z_2,\dots ,z_m));z_1),$$
where ${{\underline p}_1}:=(p_2,\dots ,p_m)$ and $n_{>j}:={\displaystyle \sum _{i=j+1}^{q}}\vert x_{i}\vert $.
\end{enumerate}
\medskip

\noi {\bf b)} let $m_1:={\displaystyle \sum _{i=1}^m}\vert z_i\vert $, for a partition ${\underline q}=(q_1,\dots ,q_{m+1})$ of $n$, and an integer $0<k\leq n_{>q_1}+\vert y\vert +m_1$, the element 
$$L_{\underline q}^{jk}({\bold x},y,{\bold z},t):=L_{q_1}^{j+k+n_{>q_1}+m_1}(x_1,\dots ,x_{q_1};L_{{\underline q}_1}^j({\bold x}^{q_1} ,y, {\bold z});t),$$
where ${{\underline q}_1}:=(q_2,\dots ,q_{m+1})$ and ${\bold x}^{q_1}:=(x_{q_1+1},\dots ,x_n)$.
\medskip

\noi {\bf c)} for a partition ${\underline r}=(r_1,\dots ,r_{l+1})$ of $n$, the element 
$$ L_{l{\underline r}}^k({\bold x},y,{\bold z},w,t):=L_{r_1+m-l+1}^k(x_1,\dots ,x_{r_1},L_{{\underline r}_1}^{\vert y\vert}({\bold x}^{r_1},y,{\bold z}^{\leq l}),z_{l+1},\dots ,z_m;w;t),$$
where ${\underline r}_1:=(r_2,\dots ,r_l)$, ${\bold x}^{r_1}:=(x_{r_1+1},\dots ,x_n)$ and ${\bold z}^{\leq l}:=(z_1,\dots ,z_l)$.}
\end{notn}
\ms

\begin{thm}\label{primpersh}  Let $(A,\bullet _{i})$ be a preshuffle algebra over $K$. Given elements $x_{1},\dots ,x_{n}$, $y,z_{1},\dots ,z_{m},w,t$ of $A$, the operations $L_{n}^{j}$ verify the following relations:
\ms
 $$\displaylines { 
{\bf a)}\ L_{n}^{j}(x_{1},\dots ,x_{n};y;L_{0}^{k}(w;t))=\sum _{r=0}^{n}L_r^{j+k+n_{>r}}(x_1,\dots ,x_r;L_{n-r}^j(x_{r+1},\dots ,x_n;y;w);t)+\cr
\delta _{j\vert y\vert}L_{n+1}^k(x_1,\dots ,x_n,y;w;t)-\delta _{k\vert w\vert}\sum _{r=0}^{n}L_r^{j+n_{>r}}(L_r^{j+k+n_{>r}}(x_1,\dots ,x_r;L_{n-r}^{j}(x_{r+1},\dots ,x_n;y;t);w),\cr }$$
where $\delta _{pq}:=\begin{cases} 1,& {\rm for}\ p=q,\\
0,&{\rm otherwise}.\end{cases}$
\bs

\noi {\bf b)} For $m\geq 1$ , 
$$\displaylines {
L_{n}^{j}(x_{1},\dots ,x_{m};y;L_{m}^{k}(z_{1},\dots ,z_{m};w;t))= \sum _{\underline p} \bigl ( L_{\underline p}^{jk}({\bold x},y,{\bold z}\rq ,t)-\hfill\cr
 L_{p_1}^{j+n_{>p_1}}(x_1,\dots ,x_{p_1};L_{{\underline p}_1}^{jk}({\bold x}^{p_1} ,y,{{\bold z}\rq }^{1},t);z_1)\bigr )+
 \delta _{j\vert y\vert} \bigl ( \sum _{l=1}^m \sum _{\underline q}L_{l{\underline q}}^k({\bold x},y,{\bold z},w,t)-\cr
\hfill  \sum _{l=2}^m(\sum _{\underline r}L_{r_1}^{\vert y\vert +n_{>r_1}}(x_1,\dots ,x_{r_1};L_{(l-1){\underline r}_1}^k({\bold x}^{r_1},y,{\bold z}^1,y,w,t);z_1)+
L_{n+m+1}^{k}(x_1,\dots ,x_n,y,z_1,\dots ,z_m;w;t)-\cr
\sum _{s=0}^nL_s^{\vert y\vert + n_{>s}}(x_1,\dots ,x_s;L_{n+m-s}^k(\dots ,x_n,y,z_2,\dots ,z_m;w;t);z_1)\bigr ),\cr }$$
where the first sum is taken over all partitions ${\underline p}=(p_1,\dots ,p_{m+1})$, the second one is taken over all partitions ${\underline q}=(q_1,\dots ,q_{l+1})$, and the third one over all ${\underline r}=(r_1,\dots ,r_{l+1})$ of $n$, ${\bold z}\rq :=(z_1,\dots ,z_m,w)$ and ${\bold y}^{s}:=(y_{s+1},\dots ,y_h)$ for any partition ${\bold y}=(y_1,\dots ,y_h)$.
\end{thm}
\bs

\P For $n,m=0,1$ the formulas may be checked by a straightforward calculation. The other cases are obtained by recursive arguments on $n$ and $m$, applying Lemma 
\ref{lem*} and the following formula:
$$L_{n+1}^{i}(x_{1},\dots ,x_{n+1};y;w)=L_{n}^{i+\vert x_{n+1}\vert  }(x_{1},\dots ,x_{n};x_{n+1}\bullet _{0}y;w),$$
for $q\geq 1. \hfill \diamondsuit $
\bs

\begin{defn}{\rm A {\it ${\mbox {\it Prim}_{psh}}$ algebra} is a graded vector space $V$ equipped with operations $L_{n}^{j}:V^{\otimes n}\otimes V_{m}\otimes V
\longrightarrow V$, for $n\geq 0$ and $1\leq j\leq m$, which verify the relations of the Theorem \ref{primpersh}. }
\end{defn}

Note that the Theorem \ref{primpersh} implies that the free ${\mbox {\it Prim}_{psh}}$ algebra over a set $X$ is linearly spanned by elements $x=L_{n}^{p}(x_{1},\dots ,x_{n};y;z)$, 
where $x_{1},\dots ,x_{n},y$ are elements in the free algebra, and $z\in X$.
\ms

Theorem \ref{primpersh} states that there exists a functor from the category ${\mbox {\it Presh}}$ of preshuffle algebras to the category of ${\mbox {\it Prim}_{psh}}$ algebras. Given a  
preshuffle bialgebra $H$, the subspace of  primitive elements ${\mbox {\rm Prim}(H)}$ is a ${\mbox {\it Prim}_{psh}}$ subalgebra of $H$.
\bs

Applying the same arguments that in Section 4 for shuffle algebras, we study the structure of free preshuffle algebras..
\ms

Let $X$ be a positively graded set, since the triple $(K[{\mathcal K}_{\infty},X]_{+},\bullet _{0}, 
\Delta _{\Theta +})$ is a connected unital infinitesimal bialgebra, the map $$e(\xi _{n},x)\mapsto (\xi _{n};x)-\sum (\xi _{n_{1},n_{2}};x_{(1)},x_{(2)})+\dots 
+(-1)^{r+1}\sum (\xi _{n_{1},\dots ,n_{r}};x_{(1)},\dots ,x_{(r)})+\dots $$
 gives a bijection between $X$ and the subset $e(X)$ of ${\mbox {\rm Prim}(K[{\mathcal K}_{\infty},X])}$. 
 \ms

\noi We denote by  ${\mathcal Prim}_{psh}(X)$ the subspace of $K[{\mathcal K}_{\infty},X]$ spanned by the set $e(X)$ with the operations $L_{n}^{i}$, and 
by  ${\mathcal Prim}_{psh}(X)^{\bullet _{0}n}$ the space spanned by all the elements of the form $z_{1}\bullet _{0}z_{2}\bullet _{0}\dots 
\bullet _{0}z_{n}$, with each $z_{j}\in {\mathcal Prim}_{psh}(X)$, for $1\leq j\leq n$. Note that Theorem \ref{primpersh} states that any element $w$ in 
${\mathcal Prim}_{psh}(X)$ is a sum of elements of type $L_{i}^{n}(x_1,\dots ,x_n;y;t)$, with $x_1,\dots ,x_n,y\in {\mathcal Prim}_{psh}(X)$ 
and $t\in e(X)$.
\ms

\begin{prop}\label{todos} Let $X$ be a positively graded set, equipped with a coassociative graded coproduct $\Theta $ on $K[X]$. 
Any element $z$ in 
$K[{\mathcal K}_{\infty},X]$ belongs to ${\displaystyle \bigoplus _{n\geq 1}}{\mathcal Prim}_{psh}(X)^{\bullet _{0}n}$.
\end{prop}
\ms

\P We only need to check that an element $$z=e(\xi _{n};x)\bullet _{j}(z_{1}\bullet _{0}z_{2}\bullet _{0}\dots \bullet _{0}z_{r}),$$ with 
$x\in X_{n}$ and $z_{i}\in {\mathcal Prim}_{psh}(X)$, belongs to ${\displaystyle \bigoplus _{n\geq 1}}{\mathcal Prim}_{psh}(X)^{\bullet _{0}n}$. 
\ms

\noi We show it applying a recursive argument on $r$.

\noi If $r=0$, then $z=e(\xi _{n};x)$ belongs to $e(X)$, and the result is obvious.

\noi If $r=1$ and $0<j<\vert z_{1}\vert$, then $z=L_{0}^{j}(z_{1};e(\xi _{n};x))$ belongs to ${\mathcal Prim}_{psh}(X)$.

\noi If $r=1$ and $j=\vert z_{1}\vert $, then $z=L_{0}^{\vert z_{1}\vert }(z_{1};e(\xi _{n};x))+z_{1}\bullet _{0}e(\xi _{n};x)$ belongs to 
${\mathcal Prim}_{psh}(X)\bigoplus {\mathcal Prim}_{psh}(X)^{\bullet _{0}2}$.
\ms

\noi Suppose that $r\geq 2$. If $0<j\leq \vert z_{1}\vert +\dots +\vert z_{r}\vert$, then there exists $1\leq k\leq r$ such that $\vert z_{1}\vert +\dots +\vert z_{k-1}
\vert <j\leq \vert z_{1}\vert +\dots +\vert z_{k}\vert$, and
$$z=(e(\xi _{n};x)\bullet _{j}(z_{1}\bullet _{0}\dots \bullet _{0}z_{k}))\bullet _{0}z_{k+1}\bullet _{0}\dots \bullet _{0}z_{r}.$$
Clearly, if $k<r$ the result follows immediately by recursive hypothesis.

\noi If $k=r$, then
$$z=L_{r-1}^{j-m}(z_{1},\dots ,z_{r-1};z_{r};e(\xi _{n};x))+z_{1}\bullet _{0}(e(\xi _{n};x)\bullet _{i-\mid z_{1}\mid }(z_{2}\bullet _{0}\dots 
\bullet _{0}z_{r})),$$
where $m=\vert z_{1}\vert +\dots +\vert z_{r-1}\vert $. But $L_{r-1}^{j-m }(z_{1},\dots ,z_{r-1};z_{r};e(\xi _{n};x))\in {\mathcal Prim}_{psh}(X)$ and,
 by recursive hypothesis, $e(\xi _{n};x)\bullet _{j-\vert z_{1}\vert }(z_{2}\bullet _{0}\dots \bullet _{0}z_{r})\in {\displaystyle \bigoplus _{n\geq 1}}
 {\mathcal Prim}_{psh}(X)^{\bullet _{0}n}$. So, 
$z\in {\displaystyle \bigoplus _{n\geq 1}}{\mathcal Prim}_{psh}(X)^{\bullet _{0}n}$. \hfill $\diamondsuit$
\ms

\begin{prop}\label{fund} Let $X$ be a positively graded set. The subspace ${\mathcal Prim}_{psh}(X)$ is the subspace of primitive elements of 
$K[{\mathcal K}_{\infty},X]$. Moreover, it is the free ${\mathcal Prim}_{psh}$ algebra spanned by $X$.
\end{prop}
\ms

\P As for shuffle algebras, it suffices to prove the result for the case where $X={\displaystyle \bigcup _{n\geq 1}}X_{n}$ with $X_{n}$ finite, for all $n\geq 1$. 
 \ms
 
 \noi  Proposition \ref{elprim} states that ${\mathcal Prim}_{psh}(X)\subseteq {\mbox {\rm Prim}(K[{\mathcal K}_{\infty},X])}$, while Proposition \ref{todos} implies that 
$K[{\mathcal K}_{\infty},X]={\overline T}({\mathcal Prim}_{psh}(X))$ as a vector space. From Theorem \ref{cofree} one has that 
$K[{\mathcal K}_{\infty},X]={\overline T}({\mbox {\rm Prim}(K[{\mathcal K}_{\infty},X]))}$, so ${\mathcal Prim}_{psh}(X)={\mbox {\rm Prim}
(K[{\mathcal K}_{\infty},X])}$.
\ms

\noi For the second point, we know that the dimension of the subspace 
${\mathcal Prim}_{psh}(X)_{n}$ of homogeneous elements of degree $n$ of ${\mathcal Prim}_{psh}(X)$ is $\vert  {\mbox {\it Irr}_{{\mathcal K}_{n,X}}}\vert$.
\ms

\noi An element $(f;x_{1},\dots ,x_{r})\in {\mathcal K}_{n,X}$ is irreducible if, and only if,  $f=\xi _{n_{1},\dots ,n_{r}}\cdot \sigma $, with $\sigma \in 
Sh^{-1}(n_{1},\dots ,n_{r})\bigcap {\mbox {\it Irr}_{S_{n}}}$ and $\vert x_{i}\vert =n_{i}$. 
\ms

\noi It is easily seen that if $\sigma \in {\mbox {\it Irr}_{S_n}}$, $\tau \in S_{m}$ and $1\leq i\leq n$, then 
$\tau \bullet _{i}\sigma \in {\mbox {\it Irr}_{S_{n+m}}}$. 
\ms

\noi Moreover, let $\sigma =\sigma _{1}\times\dots \times\sigma _{k}$, with $\sigma _{i}\in {\mbox {\it Irr}_{S_{n_i}}}$ and $k\geq 2$. The permutation 
$\tau \bullet _{j}\sigma \in {\mbox {\it Irr}_{S_{n+m}}}$, for  any ${\displaystyle \sum _{i=1}^{k-1} n_{i}< j\leq n=\sum _{i=1}^{k}n_{i}}$.
\bs

\noi We need to check that the dimension of the homogeneous subspace of degree $n$ of the free ${\mathcal Prim}_{psh}$ algebra spanned by $X$ 
is precisely $\vert {\mbox {\it Irr}_{{\mathcal K}_{n,X}}}\vert$.
\ms

\noi Let $Ppsh(X)_{n }$ denotes the subspace of degree $n$ of the free ${\mathcal Prim}_{psh}$ algebra spanned by $X$.
Note that  $Ppsh(X)_{n}$ is the vector space with basis 
$$\displaylines {
{\mathbb B}_{n}:=X_n\bigcup \hfill\cr
\hfill\{ L_{q}^{j}(x_{1},\dots ,x_{q};y;z)\ \vert \ q\geq 0, 1\leq j\leq \vert y\vert, \ x_{1},\dots ,x_{q},y\in 
\bigcup _{k=1}^{n-1}{\mathbb B}_{k}\ , \ z\in X\ {\rm and} \sum _{i=1}^q\vert x_i\vert +\vert y\vert +\vert z\vert =n\}.\cr }$$
\ms

\noi The map $\alpha :{\displaystyle \bigcup _{n\geq 1}{\mathbb B}_{n}\longrightarrow \bigcup _{n\geq 1}{\mbox {\it Irr}_{{\mathcal K}_{n,X}}}}$ is defined as follows:
\begin{enumerate}\item $\alpha (x):= (\xi _{n};x)$, for $x\in X_{n}$,
\item $\alpha (L_{0}^{j}(y;z)):=(\xi _{\vert z\vert};z)\bullet _{j}\alpha (y)$, for $1\leq j\leq \vert y\vert$,
\item $\alpha (L_{q}^{j}(x_{1},\dots ,x_{q};y;z)):=$

\noi \qquad \qquad \qquad \qquad $(\xi _{\vert z\vert};z)\bullet _{j+m}(\alpha (x_{1})
\bullet _{0}\dots \bullet _{0}\alpha (x_{q})\bullet _{0}\alpha _{\vert y\vert}(y))$,

\noi for $1\leq j\leq \vert y\vert$, where $m={\displaystyle \sum _{i=1}^q}\vert x_i\vert $.
\end{enumerate}
\bs

Conversely, let $(f;x_{1},\dots ,x_{r})\in {\mathcal K}_{n,X}$ be an irreducible element, where 

\noi $f=\xi _{n_{1},\dots ,n_{r}}\cdot \sigma $. One has that
$\sigma ^{-1}(i)=\sigma ^{-1}(1)+i-1$, for $1\leq i\leq n_{1}$, with $\sigma ^{-1}(1)> 1$ and $f=\xi_{n_{1}}\bullet _{\sigma ^{-1}(1)}
(\xi _{n_{2},\dots ,n_{r}}\cdot \sigma \rq )$, for a unique $\sigma \rq \in S_{n-n_{1}}$. 
\ms

\noi There exists a unique decomposition $$(\xi _{n_{2},\dots ,n_{r}}\cdot \sigma \rq ;x_{2},\dots ,x_{r})=(g_{1};x_{2},\dots ,x_{j_{1}})\bullet _{0}
\dots \bullet _{0}(g_{m};x_{j_{m-1}+1},\dots ,x_{r}),$$ with $(g_{i};x_{j_{i-1}+1},\dots ,x_{j_{i}})$ irreducible. 

\noi By a recursive argument we suppose that 
$\alpha ^{-1}(g_{i};x_{j_{i-1}+1},\dots ,x_{j_{i}})$ is a well-defined element in ${\displaystyle \bigcup _{l=1}^{n-1}}{\mathbb B}_{l}$, for $1\leq i\leq m$. 
Since $(f;x_{1},\dots ,x_{r})\in {\mbox {\it Irr}_{{\mathcal K}_{n,X}}}$, one has that $\sum _{i=1}^{m-1}\vert g_{i}\vert <\sigma ^{-1}(1)$. So, we define 
$$\displaylines {\alpha ^{-1}(f;x_{1},\dots ,x_{r})=\hfill \cr
L_{m-1}^{s}(\alpha ^{-1}(g_{1};x_{2},\dots ,x_{j_{1}}),\dots ,\alpha ^{-1}(g_{m-1};\dots ,x_{j_{m-1}});
\alpha ^{-1}(g_{m};x_{j_{m-1}+1},\dots ,x_{r}); x_{1}),\cr }$$
where $s=\sigma ^{-1}(1)-\sum _{i=1}^{m-1}\vert g_{i}\vert $.
\ms

\noi Clearly, the map $\alpha ^{-1}$ is the inverse of $\alpha$. \hfill $\diamondsuit $
\bs

The following result is a straightforward consequence of Theorem \ref{cofree} and the previous results.
\ms

\begin{prop} Let $X$ be a positively graded set, such that $K[X]$ is equipped with a coassociative graded coproduct $\Theta$. 
The unital infinitesimal bialgebra  $K[{\mathcal K}_{\infty},X]_{+}$ is isomorphic to $T^{fc}({\mathcal Prim}_{psh}(X))$, where 
${\mathcal Prim}_{psh}(X)$ is the free ${\mathcal Prim}_{psh}$ algebra spanned by $X$.
\end{prop}
\ms

We want to prove the equivalence between the categories of connected preshuffle bialgebras and ${\mathcal Prim}_{psh}$ algebras.  More precisely, 
given a ${\mathcal Prim}_{psh}$ algebra $(V, {\overline L}_{n}^{i})$ and an homogeneous basis $X$ of the underlying vector space $V$, let 
${\mathcal U}_{Psh}(V)$ be the preshuffle bialgebra obtained by taking the quotient of the free preshuffle algebra $K[{\mathcal K},X]$ by the ideal 
(as a preshuffle algebra) spanned all the elements:
$$L_{q}^{i}(x_{1},\dots ,x_{q};y;z) -{\overline L}_{q}^{i}(x_{1},\dots ,x_{q};y;z),$$ with $x_{1},\dots ,x_{q},y,z\in X$, $ q\geq 0$ and 
$1\leq i\leq \mid y\mid$, where $L_{q}^{i}$ denotes the operations associated to the preshuffle algebra $K[{\mathcal K},X]$. 
\ms

The proof of the following result is similar to the proof of Theorem \ref{MMsh}.

\begin{thm}\label{fin} {\bf a)} Let $(H,\circ _{i},\Delta )$ be a connected preshuffle bialgebra, then $H$ is isomorphic to ${\mathcal U}_{Psh}({\rm Prim} (H))$, where 
${\mbox {\rm Prim}(H)}$ is the ${\mathcal Prim}_{psh}$ algebra of primitive elements of $H$.

\noi {\bf b)} Let $(V,{\overline L}_{q}^{i})$ be a ${\mathcal Prim}_{psh}$ algebra, then $V$ is isomorphic to ${\mbox {\rm Prim}({\mathcal U}_{Psh}(V))}$.
\end{thm}
\bs

\bs

\bs

\noi{\bf Primitive elements of grafting bialgebras.} 
\ms

\noi Let $(A,\bullet _{i},\Delta )$ be a grafting bialgebra. 

\noi By Proposition \ref{elprim}, the elements $L_{n}^{p}(x_{1},\dots ,x_{n};y;z)$ are primitive, for $1\leq p<
\vert y\vert$, whenever the elements $x_{1},\dots ,x_{n},y,z$ belong to ${\mbox {\rm Prim}(A)}$. But an easy calculation shows that $L_{n}^{p}(x_{1},\dots ,x_{n};y;z)=0$ for any 
$x_{1},\dots ,x_{n},y,z\in A$ and $n\geq 1$. Using this fact, we introduce the following definition.  

\begin{defn}\label{primgraft}{\rm A ${\mathcal Prim}_{gr}$ algebra over $K$ is a graded vector space $V$ equipped with a family of binary operations $\{-,-\}:V\otimes V
\longrightarrow V$ and $\bullet _{p}:V\otimes V_{n}\longrightarrow V$, for $1\leq p < n$, such that:}
\begin{enumerate}\item $\{\{x,y\},z\}=\{x,\{y,z\}\} + y\bullet _{\vert x\vert}\{x,z\}$, for $x,y,z\in V$.
\item $\{x\bullet _{p}y,z\}=x\bullet _{p}\{y,z\}$,
\item $\{x, y\bullet _{p}z\}=y\bullet _{\vert x\vert +p}\{ x,z\}$,
\item $\{ x,y\}\bullet _{p}z=y\bullet _{\vert x\vert +p}(x\bullet _{p}z)-x\bullet _{p}(y\bullet _{p}z)$, for $1\leq p<\vert z\vert$,
\item $(x\bullet _{p}y)\bullet _{q}z=x\bullet _{p+q}(y\bullet _{q}z)$,
\item $x\bullet _{p}(y\bullet _{q}z)=y\bullet _{\vert x\vert +q}(x\bullet _{p}z)$, if $1\leq p<q<\vert z\vert $,
\end{enumerate}
{\rm for} $x,y,z\in V$.\end{defn}

Clearly, any grafting bialgebra $(A,\bullet _{i},\Delta )$ has a natural structure of ${\mathcal Prim}_{gr}$ algebra, such that ${\mbox {\rm Prim}(A)}$ is a ${\mathcal Prim}_{gr}$ 
subalgebra of $A$.

For any positively graded set $X$ and any coassociative coproduct $\Theta $ on $K[X]$, there exists a natural extension of the coproduct to a coassociative coproduct 
$\Delta _{\Theta}$ such that $(Graf(X),\bullet _{i},\Delta _{\Theta })$ is a grafting bialgebra. Moreover, the vector space $K[T_{\infty},X]_{+}$ equipped with the 
associative product $\circ _{0}$ and $\Delta _{\Theta +}$ is a unital infinitesimal connected bialgebra, so it is isomorphic to $T^{fc}({\mbox {\rm Prim} (K[T_{\infty},X]_{+}))}$.
 \ms
 
Let ${\mathcal Prim}_{gr}(X)$ the subspace of $K[T_{\infty},X]$ spanned by $e(X)$ with the operations $\{-,-\}$ and $\circ _{p}$.

\begin{prop}\label{todosgraft} Let $X$ be a positively graded set, equipped with a coassociative graded coproduct $\Theta $ on $K[X]$. Any element $z$ in 
$K[T_{\infty},X]$ may be written as a sum  $z=\sum _{k}z_{1}^{k}\circ _{0}z_{2}^{k}\circ _{0}\dots \circ _{0}z_{r_{k}}^{k}$, 
with $z_{i}^{k}\in {\mathcal Prim}_{gr}(X)$.
\end{prop}
\ms

\P The space $K[T_{\infty},X]$ is a quotient of $K[{\mathcal K}_{\infty},X]$, let 

\noi $\Pi : K[{\mathcal K}_{\infty},X]\longrightarrow 
K[T_{\infty},X]$. Let $e_{\mathcal K}$ (respectively, $e_{T}$) denotes the projection of $K[{\mathcal K}_{\infty},X]$ (respectively, $K[T_{\infty},X]$) 
into its primitive part. 

\noi The set $\Pi ^{-1}(e_{\mathcal K}(x))$ has a unique element, for any $x\in X$; which implies that the restriction of $\Pi$ to 
$e_{\mathcal K}(X)$ is an injective map, whose image is $e_{T}(X)$. 

\noi Moreover, since $\Pi$ sends the product $\bullet _{p}$ to $\circ _{p}$, for $p\geq 0$, we have that 

\noi $\Pi ({\mathcal Prim}_{psh}(X))\subseteq 
\Pi ({\mathcal Prim}_{gr}(X))$.

\noi Let $z\in K[T_{\infty},X]$, there exist at least one element ${\tilde z}\in K[{\mathcal K}_{\infty},X]$ such that 

\noi $\Pi ({\tilde z})=z.$

\noi We know that ${\tilde z}=\sum _{k}{\tilde z}_{1}^{k}\bullet _{0}{\tilde z}_{2}^{k}\bullet _{0}\dots \bullet _{0}{\tilde z}_{r_{k}}^{k}$, with
${\tilde z}_{i}^{j}\in {\mathcal Prim}_{psh}(X)$. 

\noi So, $z=\sum_{k}\Pi ({\tilde z}_{1}^{k})\circ _{0}\Pi ({\tilde z}_{2}^{k})\circ _{0}\dots 
\circ _{0}\Pi ({\tilde z}_{r_{k}}^{k})$, with $\Pi ({\tilde z}_{1}^{k})\in {\mathcal Prim}_{gr}(X)$.\hfill $\diamondsuit$
\bs

\begin{prop}\label{fundgraf} Let $X$ be a positively graded set, equipped with an associative coproduct $\Theta$ on $K[X]$. The subspace ${\mathcal Prim}_{gr}(X)$
 is the subspace of primitive elements of $K[T_{\infty},X]$. Moreover, it is the free ${\mathcal Prim}_{gr}$ algebra spanned by $X$.
\end{prop}
\ms

\P The proof of the first assertion is identical to the ones given for preshuffle and shuffle algebras. The unique point to see is that ${\mathcal Prim}_{gr}(X)$ is the free 
${\mathcal Prim}_{gr}$ algebra spanned by $X$.
\ms

To prove the second one, we may suppose that the set $X_{n}$ is finite, for all $n\geq 1$. 

\noi Since the algebra $({\mbox {\rm Prim}(K[T_{\infty},X])},\bullet _{0})$ is free on the set $\{ t=\bigvee _{x}(\vert ,
t^{1},\dots ,t^{r})\}$, the dimension of the subspace of homogeneous elements of degree 
$n$ of ${\mbox {\rm Prim}(K[T_{\infty},X])}$ is the number of trees of the form $t=\bigvee _{x}(\vert ,
t^{1},\dots ,t^{r})$, where $x\in X_{r}$ and 

\noi $\vert t\vert =\sum _{1\leq j\leq r}\vert t^{j}\vert +r-1$.

\noi Let $\{X\}$ be the set of all elements of the form $z=\{ x_{1},\{\dots ,\{x_{k-1},x_{k}\}\}\}$, with $k\geq 1$ and $x_{i}\in X$. 
From Definition \ref{primgraft} we have that the elements of $X$ and the elements of type: 
$$z=x_{1}\bullet _{i_{1}}(\dots \bullet _{i_{r-2}}(x_{r-1}\bullet _{i_{r-1}}(x_{r}\bullet _{i_{r}}w))),$$
with $i_{1}>\dots>i_{r},\ 1\leq i_{j}<
\vert x_{j+1}\vert +\dots +\vert x_{r}\vert +\vert w\vert,\ x_{j}\in X$, and $w\in \{ X\} ,$
are a basis of the free ${\mathcal Prim}_{gr}$ algebra over $X$ as a vector space.

\noi Define a map $\gamma$ from the basis described above to the set $\{ t=\bigvee _{x}(\vert ,t^{1},\dots ,t^{r}) \mid \ {\rm with}\ x\in X_{r}\ {\rm and}\ \vert t\vert =
\sum _{1\leq j\leq r}\vert t^{j}\vert +r-1\} $, as follows:
$$\displaylines {
\gamma (x):=(\cc _{n},x),\ {\rm for}\ x\in X\hfill\cr
\gamma (\{ x_{1},\{\dots ,\{x_{k-1},x_{k}\}\}\}):=\bigvee _{x_{1}}(\vert, \gamma (\{ x_{2},\{\dots ,\{x_{k-1},x_{k}\}\}\})),\ {\rm for}
\ x_{1},\dots ,x_{k}\in X\hfill\cr
\gamma (x_{1}\bullet _{i_{1}}(\dots \bullet _{i_{r-2}}(x_{r-1}\bullet _{i_{r-1}}(x_{r}\bullet _{i_{r}}w)))):=
(\cc_{n_{1}},x_{1})\circ _{i_{1}}(\dots ((\cc _{n_{r}},x_{r})\circ _{i_{r}}\gamma (w))),\hfill \cr}$$
for $\vert x_{i}\vert =n_{i}.$ 

\noi Clearly, $\gamma$ is a graded bijection, which sends elements of degree $n$ of the basis to trees of the same degree. So, ${\mathcal Prim}_{gr}(X)$ is a quotient of the 
free ${\mathcal Prim}_{gr}$ algebra over $X$ such that both space have the same dimension on each degree, which implies they are isomorphic. \hfill $\diamondsuit$
\ms

Again, we have a natural equivalence between the categories of connected grafting bialgebras and ${\mathcal P}_{Gr}$ algebras.  As in previous cases we define, for any
 ${\mathcal P}_{Gr}$ algebra $(V,  [-,-],\circ _{p})$ and an homogeneous basis $X$ of $V$, the universal grafting envelopping algebra 
${\mathcal U}_{Gr}(V)$ as the quotient of the free grafting algebra $K[T_{\infty},X]$ by the ideal 
 spanned by the elements:

\noi $\{x,y\}-[x,y]$ and $x\bullet _{p}y-x\circ _{p}y$, for $x,y\in X$ and $1\leq p<\vert y\vert$
, where $\{-,-\}$ and $\bullet _{p}$ denote the operations associated to the grafting algebra $K[T_{\infty},X]$. 
\ms

The proof of the following result is similar to the proof given for shuffle and preshuffle bialgebras.

\begin{thm}\label{fin} {\bf a)} Let $(H,\circ _{i},\Delta )$ be a connected grafting bialgebra, then $H$ is isomorphic to ${\mathcal U}_{gr}({\mbox {\rm Prim} (H))}$, where 
${\mbox {\rm Prim}(H)}$ is the ${\mathcal P}_{Gr}$ algebra of primitive elements of $H$.

\noi {\bf b)} Let $(V,\{-,-\},\bullet _{p})$ be a ${\mathcal P}_{Gr}$ algebra, then $V$ is isomorphic to ${\mbox {\rm Prim}({\mathcal U}_{Gr}(V))}$.
\end{thm}

\bs

\bs

\bs

\section{Applications.}
\ms

\noi{\bf I. Basis of primitive elements for the Malvenuto-Reutenauer algebra and for the bialgebra $K[{\mathcal P}_{\infty}]$.}
\ms

In \cite {DHT} and \cite{AS1} the authors describe different basis for the subspace of primitive elements of the Malvenuto-Reutenauer bialgebra. We construct another 
one using about preshuffle bialgebras. We may extend this basis to a basis of the subspace of primitive elements of $K[{\mathcal P}_{\infty}]$.
\bs

{\bf 1) The Malvenuto-Reutenauer bialgebra.} 
\ms

We know that the dimension of the subspace of primitive elements of degree $n$ of $K[S_{\infty }]$ is the number $\vert {\mbox {\it Irr}_{S_n}}\vert $ of irreducible permutations of 
$S_{n}$. Using Proposition \ref{pshuff}, we associate to any $\sigma \in {\mbox {\it Irr}_{S_n}}$ a primitive element $E_{\sigma }$ in $K[S_{\infty }]$. 

\noi If $\sigma =(1)$, then $E_{(1)}:=(1)$.

\noi Let $\sigma \in {\mbox {\it Irr}_{S_n}}$ with $n\geq 1$, there exist a family of irreducible permutations $\sigma _1,\dots ,\sigma _r$ and  a shuffle $\delta \in {\mbox {\it Sh}(1,n-1)}$ such that $\sigma =((1)\times 
\sigma_1\times \dots \times \sigma _r )\cdot \delta $. The integer $r$ and the permutations $\sigma _1,\dots ,\sigma _r,\delta$ are unique.

\noi For example, $\sigma =(5,1,4,6,3,2)=((1)\times (4,3,5,2,1))\cdot (2,1,3,4,5,6)$, and 
$\sigma \rq =(3,4,2,5,1)=((1)\times (2,3,1)\times (1))\cdot (2,3,4,5,1)$.

\noi Suppose that $\vert \sigma _i\vert = n_i$, with $\displaystyle \sum _{i=1}^r n_i=n-1$. If $\sigma ^{-1}(1)\leq n_1+\dots +n_{r-1}+1$, then $\sigma = (((1)\times \sigma _1\times \dots 
\times \sigma _{r-1})\cdot \delta \rq )\cdot \sigma _r$, which is impossible because $\sigma $ is irreducible. So, $\sigma ^{-1}(1)-1-n_1-\dots -n_{r-1}>0$.
\medskip

\noi Define $E_{\sigma }$ as the following primitive element,
$$E_{\sigma }:=L_{r-1}^{\sigma ^{-1}(1)-1-n_1-\dots -n_{r-1}}(E_{\sigma _1},\dots ,E_{\sigma _{r-1}};E_{\sigma _r};(1)),$$
where the operations $L_i^j$ are the operations introduced of Definition \ref{primpre}.

From Proposition \ref{fund} we get that the set $\{ E_{\sigma }\}_{\sigma\in \bigcup {\mbox {\it Irr}_{S_n}}}$ is a basis of the subspace of primitive elements of the Malvenuto-Reutenauer bialgebra. For instance, one has that:
$$\displaylines {
E_{(2,1)}= (2,1)-(1,2),\hfill\cr
E_{(3,1,2)}= L_0^1(E_{(2,1)};(1))=(3,1,2)-(2,1,3),\hfill\cr
E_{(3,4,2,5,7,1,6)}=L_2^1(E_{(2,3,1)},E_{(1)};E_{(2,1)};(1))=\hfill\cr
(3,4,2,5,7,1,6) -(2,3,1,5,7,4,6)-(3,4,2,5,6,1,7)+(2,3,1,5,6,4,7)-\cr
\hfill (2,4,3,5,7,1,6) +(1,3,2,5,7,4,6)+(2,4,3,5,6,1,7)-(1,3,2,5,6,4,7).\cr }$$
\bs

{\bf 2) The bialgebra $K[{\mathcal P}_{\infty}]$.} As a shuffle algebra $K[{\mathcal P}_{\infty}]$ is the free shuffle algebra spanned by the family $\{ \xi _n
\}_{n\geq 1}$. So, any coassociative coproduct  $\Theta $ on the vector space spanned by   $\{ \xi _n\}_{n\geq 1}$ gives a shuffle bialgebra structure on 
$K[{\mathcal P}_{\infty}]$. Note that $$\Theta ^r(\xi _n)=\sum _{n_1+\dots +n_r}c_{n_1\dots n_r}\xi _{n_1}\otimes \dots \otimes \xi _{n_r},$$
where $(n_1,\dots ,n_r)$ is a composition of $n$, and $c_{n_1\dots n_r}\in K$.

\noi So, $$E_{\xi _n}:=\sum _{r=1}^n(-1)^{r-1}(\sum _{n_1+\dots +n_r=n}c_{n_1\dots n_r}\xi _{n_1}\bullet _0\dots \bullet _0 \xi _{n_r}).$$

\noi Given an irreducible element $f\in {\mbox {\it Irr}_{{\mathcal P}_{n}}}$, we associate to it a primitive element $E_{\theta}(f)$ in 
$(K[{\mathcal P}_{\infty}], \Delta _{\theta })$ as follows. Let $n_{1}=\vert f^{-1}(1)\vert $, 

\noi If $n_{1}=n$, then $f=\xi _{n}$ and $E_{\theta}(f)$ is defined above.

\noi If $n_{1}<n$, then $f=\xi_{n_{1}}\bullet _{\gamma}f_{1}$, with $\gamma \neq 1_{n}$.There exist a unique family $g_1,\dots ,g_r$ of irreducible elements, such that 
$f_1=g_1\bullet _0\dots \bullet _0g_r$. Let $\vert g_j\vert =m_j$, for $1\leq j\leq r$. Since $f\in {\mbox {\it Irr}_{{\mathcal P}_{n}}}$, it is immediate that
 $\gamma (n_{1})> \gamma (n_{1}+m_{1}+\dots +m_{r-1}+1) $.

\noi There exists  $0\leq k\leq r-1$ such that $\gamma (n_{1}+m_{1}+\dots +m_{k}+1)<\gamma (1)$. \begin{enumerate}
\item If $k=0$, then $E_{\theta }(f):=E_{\theta}(\xi _{n_1})\bullet _{\gamma }(E_{\theta }(g_1)\bullet _0\dots \bullet _0 E_{\theta }(g_r)).$
\item If $k\geq 1$, then $$E_{\theta }(f):=\{E_{\theta }(g_1),E_{\theta }(\xi _n)\}\bullet _{\tilde {\gamma }}(E_{\theta }(g_2)\bullet _0\dots \bullet _0 E_{\theta}(g_r),$$
where $\gamma =(\epsilon _{n_1,m_1}\times 1_{m_2+\dots +m_r})\cdot {\tilde {\gamma }}$, with ${\tilde {\gamma }}=1_{m_1}\times \alpha $ and 

\noi $\alpha (n_1)>\alpha (n_1+m_2+\dots +m_{r-1}+1$.
\end{enumerate}

Applying Proposition \ref{fund}, we get that the family $\{E_{\theta}(f)\}_{f\in {\mbox {\it Irr}_{{\mathcal P}_{\infty}}}}$ is a basis of the space of primitive elements of 
$(K[{\mathcal P}_{\infty}], \Delta _{\theta})$.
\medskip

\noi For example, let $\Theta (\xi _n)={\displaystyle \sum _{i=1}^{n-1}}\xi _i\otimes \xi _{n-i}$, for $n\geq 1$, and let 
$$f=(3,2,4,1,6,4,1,5,5)=\xi_2\bullet _{(3,4,5,1,6,7,2,8,9)}((2,1)\bullet _0(1,3,1,2,2)).$$

\noi We get that $$\displaylines {
E_{\theta }(\xi _2)=(1,1)-(1,2),\hfill\cr
E_{\theta}(2,1)=(2,1)-(1,2),\hfill\cr
E_{\theta}(2,1,1)=\{ E_{\theta}(\xi _1),E_{\theta}(\xi _2)\}=(2,1,1)-(3,1,2)-(1,2,2)+(1,2,3),\hfill\cr
E_{\theta}(1,3,1,2,2)=E_{\theta}(\xi _2)\bullet _{(1,3,2,4,5)}E_{\theta}(2,1,1)=(1,3,1,2,2)-(1,4,1,2,3)-\hfill\cr
\hfill (1,2,1,3,3)+(1,2,1,3,4)-(1,4,2,3,3)+(1,5,2,3,4)+(1,3,2,4,4)-(1,3,2,4,5),\cr
E_{\theta}(3,2,4,1,6,4,1,5,5)=\{E_{\theta}(2,1),E_{\theta}(\xi _2)\}\bullet _{(1,2,5,3,6,7,4,8,9)}E_{\theta}(1,3,1,2,2).\hfill\cr }$$

\bs

\noi {\bf II. Some triples of operads}
\ms

Note that preshuffle algebras, shuffle algebras and grafting algebras are not described by classical linear operads. This Section is devoted to describe somme {\it good triples} 
of $K$-linear operads, following the definition of \cite{Lo1}, where the co-operad is always the co-associative operad. 

{\bf 1) Duplicial bialgebras.} 
\ms

This example is studied in \cite{Lo1}, we include it since the results may be obtained easily from our computation of primitive elements.

\begin{defn} {\rm A duplicial algebra over $K$ is a vector space $A$ equipped with two bilinear maps $/$, $\backslash :A\otimes A\longrightarrow A$, verifying 
the following relations:}
$$\displaylines {
x/(y/z)=(x/y)/z\hfill \cr
x/(y\backslash z)=(x/y)\backslash z\hfill\cr
x\backslash (y\backslash z)=(x\backslash y)\backslash z,\hfill \cr}$$
{\rm for $x,y,z\in A$}.\end{defn}

Note that for any grafting algebra $(A,\bullet _{i})$, the space $A$ with the products:
$$x/y:=x\bullet _{0}y\ {\rm and}\ x\backslash y:=y\bullet _{\mid y\mid}x$$
is a duplicial algebra.

It is not difficult to verify (see \cite{Pi} or \cite{Lo1}) that the free duplicial algebra spanned by a set $E$ is the space of planar binary rooted trees $K[Y_{\infty},E]$, with the 
vertices coloured by the elements of $E$. We denote it by ${\mbox {\it Dup}(E)}$.
 
 \noi Let $(A,/,\backslash )$  be a duplicial algebra, a coassociative coproduct on $A$ is admissible for the duplicial structure if
 $$\displaylines{
 \Delta (x/y)=\sum x_{(1)}\otimes (x_{(2)}/y)+\sum (x/y_{1})\otimes y_{(2)}+x\otimes y,\hfill \cr
 \Delta (x\backslash y)=\sum x_{(1)}\otimes (x_{(2)}\backslash y)+\sum (x\backslash y_{1})\otimes y_{(2)}+x\otimes y,\hfill \cr}$$
 for $x,y\in A$.
 
 An duplicial bialgebra is an duplicial algebra $A$ equipped with an admissible coproduct. 
 
 \noi Note that any grafting bialgebra is a duplicial bialgebra. In particular, for any set $E$, the free duplicial algebra ${\mbox {\it Dup}(E)}$ is a duplicial bialgebra.
 
 \noi Clearly, the unique operation of ${\mathcal Prim}_{gr}$ which may be defined in any duplicial algebra is the product $\{ -,-\}$, which does not verify any relation.
 
 \begin{defn} {A magmatic algebra over $K$ is a vector space $M$ equipped with a bilinear map $M\otimes M\longrightarrow M$.}
 \end{defn}
 
 There exists a functor $F_{Dup-Mag}$ from the category of duplicial algebras to the category of magmatic algebras, which maps $(A,/,\backslash )\mapsto (A,\{-,-\})$.
 If $(A,/,\backslash,\Delta)$  is a duplicial bialgebra, then $({\mbox {\rm Prim}(A)},\{-,-\})$ is a magmatic subalgebra of $(A,\{-,-\})$.
 
\ms

For any set $E$, let $\{ E,E\}$ denote the subspace of the free duplicial algebra ${\mbox {\it Dup}(E)}$ spanned by the elements of $E$ under the operation $\{-,-\}$ and let 
${\overline T}\{ E,E\}$ be the subspace of ${\mbox {\it Dup}(E)}$ spanned by the elements of the form $z=z_{1}/\dots ,/z_{n}$, with $n\geq 1$ and $z_{i}\in \{ E,E\}$ for 
$1\leq i\leq n$.

\begin{prop} The space ${\overline T}\{ E,E\}$ is isomorphic to ${\mbox {\it Dup}(E)}$.\end{prop}
\ms

\P We need to prove that any planar binary tree $t$ with the vertices coloured with elements of $E$ is a finite sum ${\displaystyle \sum _{j}}z_{1}^{j}/\dots /z_{n_{j}}^{j}$, with 
$z_{i}^{j}\in \{ E,E\}$. Note that it suffices to prove the result for the trees of the form $t=\vert \vee _{e}t\rq $, with $e\in E$ and $t\rq \in {\mbox {\it Dup}(E)}$.

\noi In order to simplify notation we denote $e$ the tree $(\cc _{1},e)$, for $e\in E$.

If $\vert t\vert =1$, then $t=e\in \{ E,E\}$.

If $\vert t\vert =2$, then $t=e\backslash f=\{ e,f\}+e/f\in {\overline T}\{ E,E\}$.

The proof follows from the following assertions:\begin{enumerate}
\item Let $t=w_{1}/\dots /w_{r}\backslash z$, with $w_{i},z\in \{ E,E\}$. We have that:
$$t=w_{1}/\dots /w_{r-1}/\{ w_{r},z\} + w_{1}/\dots /w_{r}/z\in {\overline T}\{ E,E\}.$$

\noi A recursive argument on $m$ proves that if $w_{1}/\dots /w_{r}\backslash z_{1}\backslash \dots \backslash z_{m-1}\in {\overline T}\{ E,E\}$, then 
$w_{1}/\dots /w_{r}\backslash z_{1}\backslash \dots \backslash z_{m}\in {\overline T}\{ E,E\}$, with the elements $w_{i}$ and $z_{j}$ in $\{ E,E\}$.
\item Let $z_{1},\dots ,z_{m}\in \{ E,E\}$, then
$$\displaylines {\{ z_{1},\{ z_{2},\dots ,\{ z_{m-1},z_{m}\}\dots \}\}=z_{1}\backslash \{z_{2},\dots ,\{z_{m-1},z_{m}\} \dots \}-\hfill \cr
\hfill -z_{1}/z_{2}\backslash \{z_{3},\dots ,\{ z_{m-1},z_{m}\}\dots \} +\dots +(-1)^{m-1}z_{1}/z_{2}/\dots /z_{m}. \cr }$$
\medskip

\noi For $m=2$, using the first formula, we get that:
$$w_{1}/\dots /w_{r}\backslash (z_{1}/z_{2})= -w_{1}/\dots w_{r-1}/w_{r}\backslash \{z_{1},z_{2}\} +w_{1}/\dots /w_{r}\backslash z_{1}
\backslash z_{2},$$
for $w_{1},\dots ,w_{r},z_{1},z_{2}\in \{ E,E\}$. So, $w_{1}/\dots /w_{r}\backslash (z_{1}/z_{2})\in {\overline T}\{ E,E\}$.

\noi For $m\geq 3$, a recursive argument and the formulas above imply that $w_{1}/\dots /w_{r}\backslash (z_{1}/z_{2}/\dots /z_{m})\in {\overline T}\{ E,E\}$.
\hfill $\diamondsuit$
\end{enumerate}
\ms

Using that for any duplicial bialgebra $(A,/,\backslash ,\Delta)$, the triple $(A_{+},/,\Delta _{+})$ is a unital infinitesimal bialgebra, we get that 
for any set $E$, the subspace $\{ E,E\}$ is equal to ${\mbox {\rm Prim} ({\mbox {\it Dup}(E)})}$, and that, as coalgebras $T^{fc}(\{ E,E\})$ and ${\mbox {\it Dup}(E)}$ are isomorphic.
\ms

\noi To end the example we only need to prove that $(\{ E,E\},\{-,-\})$ is the free magmatic algebra spanned by $E$, denoted by ${\mbox {\it Mag}(E)}$. 
\ms

\begin{prop} For any set $E$, the subspace of primitive elements of ${\mbox {\it Dup}(E)}$, equiped with the binary product $\{-,-\}$ is the free magmatic algebra 
${\mbox {\it Mag}(E)}$ over $E$.
\end{prop}
 \ms
 
 \P Again, we prove the result for any finite set $E$, showing that the dimension of the subspaces of homogeneous elements of degree $n$ of 
 ${\mbox {\rm Prim}({\mbox {\it Dup}(E)})}$ and ${\mbox {\it Mag}(E)}$ are the same, for all $n$.
 
 \noi We know (see \ref{Hol}) that the dimension of the space of homogeneous elements of degree $n$ of ${\mbox {\it Mag}(E)}$ is $c_{n-1}\vert E\vert ^{n}$, where $c_{n-1}$ is the 
 Catalan number, which counts the number of planar binary rooted trees with $n$ leaves.
 
 \noi On the other hand, we know that $({\mbox {\it Dup}(E)}, /)$ is the free associative algebra spanned by the trees of type $t=\vert \vee _{e}t\rq$, with $e\in E$ and 
 $\vert t\rq \vert =\vert t\vert -1$. Since, we also have that ${\mbox {\it Dup}(E)}$ is isomorphic as a coalgebra to $T^{fc}({\mbox {\rm Prim}({\mbox {\it Dup}(E)})})$, 
 we may assert that the dimension of the subspace of homogeneous elements of ${\mbox {\rm Prim}({\mbox {\it Dup}(E)})}$ of degree $n$ is the number ot trees of type $\vert \vee _{e}t\rq$, where $t\rq $ is a planar binary rooted tree with $n$ 
  leaves and its $n-1$ internal vertices coloured with elements of $E$. So, we have that $dim _{K}({\mbox {\rm Prim}({\mbox {\it Dup}(E)})})_{n}=c_{n-1}\vert E\vert ^{n}$, which ends the proof.
  \hfill $\diamondsuit$
  \ms
  
  As in the previous cases, given a magmatic algebra $(M, \cdot )$ we consider the free duplicial algebra ${\mbox {\it Dup}(E)}$ over the underlying vector space $M$. We define the universal  duplicial enveloping algebra ${\mathcal U}_{dup}(M)$ of $(M, \cdot )$ as the quotient of ${\mbox {\it Dup}(E)}$ by the ideal spanned by the elements of the form 
  $\{ x, y\}-x\cdot y$, for $x,y\in M$; where $\{-,-\}$ is the magmatic product defined on ${\mbox {\it Dup}(E)}$.
  
  The proof of a Cartier-Milnor-Moore type theorem for connected duplicial bialgebras follows using the previous result in the same way that we proved it for the cases of preshuffle 
  and shuffle bialgebras. We just state it.
  \ms
  
  \begin{thm}\label{dupli} {\bf a)} Let $(A,/,\backslash ,\Delta )$ be a connected duplicial bialgebra, then $A$ is isomorphic to ${\mathcal U}_{dup}({\mbox {\rm Prim} (A)})$, 
  where ${\mbox {\rm Prim}(A)}$ is the magmatic algebra of primitive elements of $A$.

\noi {\bf b)} Let $(M,\cdot)$ be amagmatic algebra, then $M$ is isomorphic to ${\mbox {\rm Prim}({\mathcal U}_{dup}(M))}$.
\end{thm}

\bs

\bs
 
 \bs
 
 {\bf 2) The $2$-infinitesimal nonunital bialgebra} 
 \ms
 
Recall that a ${\mbox {\it 2-ass}}$ algebra is simply a vector space equipped with two associative products $\cdot$ and $\circ $. 
Let $E$ be a set, define ${\tilde T}_{n,E}$ as the set of all planar rooted trees with $n$ leaves, with the leaves coloured by the elements of $E$. 
Consider the vector space $K[{\tilde T}(E)]$ spanned by the graded set ${\tilde T}(E):={\displaystyle \bigcup _{n\geq 1}}{\tilde T}_{n,E}$. In \cite{LR3}, we 
equipped the tensor space $T(K[{\tilde T}(E)])$ with two associative products $\cdot$ and $\circ $, and proved that $(T(K[{\tilde T}(E)]),\cdot ,\circ )$ is the free 
${\mbox {\it 2-ass}}$ algebra spanned by $E$, we denote it by ${\mbox {\it 2-ass}(E)}$. We just give a brief description of the ${\mbox {\it 2-ass}}$ structure of ${\mbox {\it 2-ass}(E)}$:\begin{enumerate}
\item For $t=\bigvee (t^{1},\dots ,t^{r})\in {\tilde T}_{n,E}$ and $w=\bigvee (w^{1},\dots ,w^{k})\in {\tilde T}_{m,E}$, 
$$t\cdot w:=\bigvee (t^{1},\dots ,t^{r},w^{1},\dots ,w^{k}),\hfill $$
and $t\circ w:=t\otimes w$.
\item For $x=t^{1}\otimes \dots \otimes t^{r}$ and $y=w^{1}\otimes \dots \otimes w^{k}$, with $t_{1},\dots ,t_{r},w_{1},\dots ,w_{k}\in {\tilde T}(X)$,
$$x\cdot y:=(\bigvee (t^{1},\dots ,t^{r}))\vee (\bigvee (w^{1},\dots ,w^{k})),$$
and $x\circ y:=t^{1}\otimes \dots \otimes t^{r}\otimes w^{1}\otimes \dots \otimes w^{k}$.
\end{enumerate}
\ms

\begin{defn} {\rm A $2$-infinitesimal nonunital bialgebra is a ${\mbox {\it 2-ass}}$ algebra $(A,\cdot ,\circ )$ equipped with a coassociative coproduct $\Delta $, such that the triples 
$(A_{+},\cdot ,\Delta _{+})$ and $(A_{+},\circ ,\Delta _{+})$ are infinitesimal unital bialgebras. That means that $\Delta $ verifies the following relations:}
$$\displaylines {
\Delta (x\cdot y)=\sum (x\cdot y_{(1)})\otimes y_{(2)}+\sum x_{(1)}\otimes (x_{(2)}\cdot y)+x\otimes y,\hfill\cr
\Delta (x\circ y)=\sum (x\circ y_{(1)})\otimes y_{(2)}+\sum x_{(1)}\otimes (x_{(2)}\circ y)+x\otimes y,\hfill\cr}$$
{\rm for $x,y\in A$, where $\Delta (x)=\sum x_{(1)}\otimes x_{(2)}$.}\end{defn}
\ms

For any preshuffle algebra $(A,\bullet _{i})$, the triple $(A, \bullet _{0},\bullet _{L})$ is a ${\mbox {\it 2-ass}}$ algebra; and if $(A,\bullet _{i},\Delta )$ is a preshuffle bialgebra, 
then $(A,\bullet _{0},\bullet _{L},\Delta )$ is a $2$-infinitesimal nonunital bialgebra. 
\ms

Let us describe the $2$-infinitesimal nonunital bialgebra structure of ${\mbox {\it 2-ass}(E)}$, for any set $E$. We define $\Delta $ recursively as follows:\begin{enumerate}
\item $\Delta (\vert ,e):= 0$, for all $e\in E$,
\item $$\displaylines {\Delta (\bigvee (t^{1},\dots ,t^{r})):= \sum _{i=1}^{r}\bigvee (t^{1},\dots ,t^{i-1},\tilde {t}_{(1)}^{i})\otimes \bigvee ({\tilde t}_{(2)}^{i},t^{i+1},
\dots ,t^{r})+\cr
{\tilde t}^{1}\otimes (\bigvee (t^{2},\dots ,t^{r}))+ \sum _{j=2}^{r-2}\bigvee (t^{1},\dots ,t^{j})\otimes \bigvee (t^{j+1},\dots ,t^{r})+\bigvee (t^{1},
\dots ,t^{r-1})\otimes {\tilde t}^{r},\cr }$$
where ${\tilde w}:=\begin{cases}w,&\ {\rm for}\ w=(\vert ,e)\\
w^{1}\otimes \dots \otimes w^{p},&\ {\rm for}\ w=\bigvee (w^{1},\dots ,w^{p}).\end{cases}$
\item $$\displaylines {
\Delta (t^{1}\otimes \dots \otimes t^{r}):=\sum _{i=1}^{r}(t^{1}\otimes \dots \otimes t^{i-1}\otimes t_{(1)}^{i})\otimes (t_{(2)}^{i}\otimes t^{i+1}
\otimes \dots \otimes t^{r})+\hfill\cr
\hfill \sum _{j=1}^{r-1}(t^{1}\otimes \dots \otimes t^{j})\otimes (t^{j+1}\otimes \dots \otimes t^{r}).\cr }$$
\end{enumerate}
\bs

\noi If we look at the structure of ${\mathcal Prim}_{psh}$ algebra of $A$, given in Definition \ref{primpre}, 
the operations which are defined using only the products $\bullet _{0}$ and $\bullet _{L}$ are the  $n+2$-ary products 
$L_{n}^{\vert y\vert}(x_{1},\dots ,x_{n};y;z)$, for $n\geq 1$. Note that they do not verify any relationship, which leads us to the following definition.
\ms

\begin{defn}{\rm A ${\mbox {\it Mag}(\infty )}$ algebra over $K$ is a vector space $M$, equipped with $n$-linear maps 
$\mu _{n}:M^{\otimes n}\longrightarrow M$, for $n\geq 2$.}\end{defn}
\ms

Let $(A,\cdot ,\circ )$ be a ${\mbox {\it 2-ass}}$ algebra, define $\mu _{n}: A^{\otimes n}\longrightarrow A$ be the operations defined as follows:
$$\displaylines {
\mu _{2}(x_{1},x_{2}):=x_{1}\cdot x_{2}-x_{1}\circ x_{2},\hfill\cr
\mu _{n}(x_{1},\dots ,x_{n}):= (x_{1}\cdot (\dots \cdot (x_{n-2}\cdot x_{n-1})))\circ x_{n}-x_{1}\cdot ((x_{2}\cdot (\dots \cdot (x_{n-2}\cdot 
x_{n-1})))\circ x_{n}),\hfill\cr }$$
for $x_{1},\dots ,x_{n}\in A$ and $n\geq 2$.

\noi Clearly, $(A, \mu _{n})$ is a ${\mbox {\it Mag}(\infty )}$ algebra.
\ms

\begin{prop} Let $(A,\cdot ,\circ ,\Delta )$ be a $2$-infinitesimal nonunital bialgebra, the subspace  ${\mbox {\rm Prim}(A)}$ of primitive elements of $A$ is closed 
under the products $\mu _{n}$, for $n\geq 2$.\end{prop}
\ms

\P The result is a straigthforward consequence of Proposition \ref{elprim}, it suffices to note that $\mu _{n}(x_{1},\dots ,x_{n})$ coincides with 
$L_{n-2}^{\mid x_{n-1}\mid}(x_{1},\dots ,x_{n-2};x_{n-1};x_{n})$,

\noi for $x\bullet _{0} y:=x\cdot y$ and $x\bullet _{\mid y\mid}y=y\circ x$.\hfill 
$\diamondsuit$
\ms

For a set $E$, let ${\mbox {\it 2-ass}(E)}$ be the free ${\mbox {\it 2-ass}}$ algebra spanned by $E$, let $M(E)$ be the subspace of ${\mbox {\it 2-ass}(E)}$ spanned by the elements of $E$ under the operations 
$\mu _{n}$ defined above, and let $T(M(E))={\displaystyle \bigoplus _{n\geq 1}}M(E)^{\otimes n}$ the tensor space over $M(E)$ equipped with the deconcatenation coproduct. 

\begin{prop} Given a set $E$, the coalgebra ${\mbox {\it 2-ass}(E)}$ is isomorphic to $T(M(E))$.
\end{prop}
\ms

\P It suffices to prove that any homogeneous element of degree $n$ of ${\mbox {\it 2-ass}(E)}$ belongs to $T(M(E))$ for all $n$. We proceed by induction on the degree of $n$. 
For $n=1$ and $n=2$ the result is immediate to check.

For a tree $t=\bigvee (t^{1},\dots ,t^{r})={\tilde t}^{1}\cdot \dots \cdot {\tilde t}^{r}$ the result is immediate because, a recursive argument states that any 
${\tilde t}^{j}$ belongs to $T(M(E))$.

Suppose then that $t=t^{1}\otimes \dots \otimes t^{r} $, since $\vert t^{1}\otimes \dots \otimes t^{r-1} \vert < \vert t\vert $ and $\vert t^{r}\vert <\vert t\vert$, we know that 
both elements belong to $T(M(E))$. So, we may restrict ourselves to prove that an element of the form $x\otimes (z^{1}\cdot \dots \cdot z^{m})$ is in $T(M(E))$, for 
$x\in T(M(E))$ and $z^{j}\in M(E)$. 

\noi Again, we proceed by a recursive argument on $m$. Note first that:
$$\displaylines { 
\mu _{2}(z^{1},\mu _{2}(z^{2},(\dots ,\mu _{2}(z_{m-1},z_{m}))))= z^{1}\cdot \dots \cdot z^{m}-\hfill\cr
z^{1}\cdot \dots \cdot z^{m-2}\cdot (z^{m-1}\circ z^{m})-z^{1}\cdot \dots \cdot z^{m-3}\cdot(z^{m-2}\circ \mu _{2}(z^{m-1},z^{m}))-\cr
\hfill \dots -z^{1}\circ \mu _{2}(z^{2},(\dots ,\mu _{2}(z_{m-1},z_{m})))),\cr }$$

\noi But, using the recursive hypothesis and the arguments above, we get that:
$$x\otimes (z^{1}\cdot \dots \cdot z^{k}\cdot (z^{k+1}\circ \mu _{2}(z^{k+2},(\dots ,\mu _{2}(z_{m-1},z_{m})))\in T(M(E)),$$
for $0\leq k< m-2,$ which implies that $x\otimes (z^{1}\cdot \dots \cdot z^{m})$ is in $T(M(E))$.\hfill $\diamondsuit$
\ms

For any $2$-infinitesimal nonunital bialgebra $(A,\cdot ,\circ ,\Delta)$, the triple $(A_{+},\cdot ,\Delta _{+})$ is a unital infinitesimal bialgebra. So, we have that
for any set $E$, the subspace $M(E)$ is the subspace of primitive elements of ${\mbox {\it 2-ass}(E)}$.  But, since as a vector space $2-ass(E)$ is just $T(K[{\tilde T}(E)]$, the dimension 
of the homogeneous elements of degree $n$ of ${\mbox {\rm Prim}({\mbox {\it 2-ass}(E)})}$ is just $C_{n-1}\vert E\vert ^{n}$, where $C_{n-1}$ is the super Catalan number. So, for any finite set $E$,
the dimension of $M(E)_{n}$ is $C_{n-1}\vert E\vert ^{n}$.
\ms

\begin{prop} For any set $E$, the subspace of primitive elements of ${\mbox {\it 2-ass}(E)}$, equiped with the $n$-ary products $\mu _{n}$ is the free ${\mbox {\it Mag}(\infty)}$ algebra over $E$.
\end{prop}
 \ms
 
 \P It suffices to note that for any finite set $E$, the subspace of homogeneous elements of degree $n$ of the free $Mag(\infty)$ algebra spanned by $E$ is the number 
 of planar rooted trees with $n$ leaves, whith the leaves coloured by the elements of $E$, which is precisely $C_{n-1}\mid E\mid ^{n}$ (for a more detailed study of this algebra 
 we refer to \cite{Hol}).\hfill $\diamondsuit$
\ms
  
  As in the previous cases, given a ${\mbox {\it Mag}(\infty)}$ algebra $(M, \mu_{n} )$ we consider the free $2-ass $ algebra spanned by $M$. We define the universal 
$2-ass$ enveloping algebra ${\mathcal U}_{2-ass}(M)$ of $(M, \mu _{n} )$ as the quotient of ${\mbox {\it 2-ass}(M)}$ by the ideal spanned by the elements of the form 
  ${\overline {\mu}}_{n}(x_{1},\dots ,x_{n})-\mu _{n}(x_{1},\dots ,x_{n})$, for $x_{1},\dots ,x_{n}\in M$; where ${\overline {\mu}}_{n}$ are the magmatic 
  products defined on ${\mbox {\it 2-ass}(M)}$.
  
 The proof of a Cartier-Milnor-Moore type theorem for connected duplicial bialgebras follows as in previous cases.
  \ms
  
  \begin{thm}\label{biass} {\bf a)} Let $(A,\cdot,\circ ,\Delta )$ be a connected  $2$-infinitesimal nonunital bialgebra, then $A$ is isomorphic to 
  ${\mathcal U}_{2-ass}({\mbox {\rm Prim} (A)})$, where ${\mbox {\rm Prim}(A)}$ is the ${\mbox {\it Mag}(\infty)}$ algebra of primitive elements of $A$.

\noi {\bf b)} Let $(M,\cdot)$ be a ${\mbox {\it Mag}(\infty)}$ algebra, then $M$ is isomorphic to ${\rm Prim}({\mathcal U}_{2-ass}(M))$.
\end{thm}
\bs

\bs

\bs

\end{document}